\documentclass[10pt,draftcls, onecolumn]{IEEEtran}

\usepackage{braket,amsfonts}

\usepackage{array}

\usepackage{pgfplots}

\newtheorem{theorem}{Theorem}
\newtheorem{lemma}{Lemma}
\newtheorem{remark}{Remark}
\newtheorem{assumption}{Assumption}

\usepackage{graphicx,epstopdf,caption,subcaption}

\usepackage{amsopn}

\usepackage{subeqnarray}
\usepackage{amssymb,amsmath}
\usepackage{ifthen}
\usepackage{setspace}
\newcommand{\R}{\mathbb{R}}
\newcommand{\qedd}{\ \hfill{$\Box$}}
\newcommand{\one}{\mathbf{1}}
\newboolean{showcomments}
\setboolean{showcomments}{true}
\usepackage{cite}
\newboolean{isfullversion}
\setboolean{isfullversion}{true}

\newboolean{appendixin}
\setboolean{appendixin}{true}

\newcommand{\lina}[1]{  \ifthenelse{\boolean{showcomments}}
	{ \textcolor{green}{(  #1)}} {}  }
\newcommand{\guannan}[1]{\ifthenelse{\boolean{showcomments}}
	{ \textcolor{green}{( #1)} } {} }
\newcommand{\guannanrevise}[1]{\ifthenelse{\boolean{showcomments}}
	{ \textcolor{black}{#1} } {#1} }

\newcommand{\gquhide}[1]{}

\makeatletter
\newcommand\footnoteref[1]{\protected@xdef\@thefnmark{\ref{#1}}\@footnotemark}
\makeatother

\allowdisplaybreaks 

	\title{Accelerated Distributed Nesterov Gradient Descent}
	\author{Guannan Qu, Na Li
		\thanks{Guannan Qu and Na Li are affiliated with John A. Paulson School of Engineering and Applied Sciences at Harvard University. Email: gqu@g.harvard.edu, nali@seas.harvard.edu. This work is supported under NSF ECCS 1608509,
			NSF CAREER 1553407, and AFOSR YIP.} \thanks{Preliminary results have been accepted to 54th Annual Conference on Communication, Control and Computing as ``Accelerated Distributed Nesterov Gradient Descent for Smooth and Strongly Convex Functions", and to 56th IEEE Conference on Decision an Control as ``Accelerated Distributed Nesterov Gradient Descent for Convex and Smooth Functions''.}}

\begin{document}
	\maketitle
	

		\begin{abstract}
This paper considers the distributed optimization problem over a network, where the objective is to optimize a global function formed by a sum of local functions, using only local computation and communication. We develop an Accelerated Distributed Nesterov Gradient Descent (Acc-DNGD) method. When the objective function is convex and $L$-smooth, we show that it achieves a $O(\frac{1}{t^{1.4-\epsilon}})$ convergence rate for all $\epsilon\in(0,1.4)$. We also show the convergence rate can be improved to $O(\frac{1}{t^2})$ if the objective function is a composition of a linear map and a strongly-convex and smooth function. When the objective function is $\mu$-strongly convex and $L$-smooth, we show that it achieves a linear convergence rate of $O([ 1 - C  (\frac{\mu}{L})^{5/7} ]^t)$, where $\frac{L}{\mu}$ is the condition number of the objective, and $C>0$ is some constant that does not depend on $\frac{L}{\mu}$. 
		\end{abstract}
		
		

	
\section{Introduction}
Given a set of agents $\mathcal{N}=\{1,2,\ldots,n\}$, each of which has a local convex cost function $f_i(x):\R^N\rightarrow \R$, the objective of distributed optimization is to find $x$ that minimizes the average of all the functions,
\begin{equation*}
\min_{x\in \R^N} f(x) \triangleq  \frac{1}{n} \sum_{i=1}^n f_i(x)， 
\end{equation*}
using local communication and local computation. The local communication is defined through an undirected connected communication graph. 
 In the past decade, this problem has received much attention and has found various applications in multi-agent systems, sensor networks, machine learning, etc \cite{bazerque2010distributed,forero2010consensus,johansson2008distributed}. 

{\color{black} Many distributed optimization algorithms have been developed based on different local computation schemes for individual agents. For instance, in Alternating Direction Method of Multipliers (ADMM) based distributed methods, e.g., \cite{shi2014linear,wei2012distributed,chang2015multi,chatzipanagiotis2015augmented}, each agent needs to solve a sub-optimization problem at each iteration. Another example is the line of work \cite{scaman2017optimal} in which each agent needs to use the \emph{dual gradient}, i.e. gradient of the Fenchel conjugate of the local cost function. 
	These methods assume that the sub-optimization or computation problem associated with ADMM or dual gradient is easy to solve, limiting their applicability. When it is not the case, it is preferable to use only local (sub-)gradient evaluation of the local objective $f_i$. The focus of this paper is such primal-only \textit{gradient based} distributed optimization. Here ``primal-only gradient'' means that the gradient of functions $f_i$, as opposed to the ``dual gradient'' used in \cite{scaman2017optimal}.  }


There is also a large literature on such gradient based distributed optimization,   \cite{tsitsiklis1984distributed,bertsekas1989parallel,nedic2009distributed,lobel2008convergence, duchi2012dual,ram2010distributed,nedic2014stochastic,nedic2015distributed,matei2011performance, olshevsky2014linear,zhu2012distributed,lobel2011distributed,lan2017communication}, most of which are distributed (sub)gradient algorithms based on a consensus/averaging scheme \cite{olshevsky2009convergence}. At each iteration, each agent performs one or multiple consensus steps plus one local gradient evaluation step. These methods \cite{nedic2009distributed,lobel2008convergence, duchi2012dual,ram2010distributed,nedic2014stochastic,nedic2015distributed,matei2011performance, olshevsky2014linear,zhu2012distributed,lobel2011distributed} have achieved sublinear convergence rates for convex functions. When the convex functions are nonsmooth, the sublinear convergence rate matches the Centralized subGradient Descent (CGD) method.  {\color{black}More recent work \cite{shi2015extra,xu2015augmented,di2015distributed,di2016next,qu2016harnessing,nedich2016achieving,nedic2016geometrically,zeng2015extrapush,xi2015linear,xi2016add,pu2018push,xin2018linear,pu2018distributed}, have improved these results to achieve linear convergence rates for strongly convex and smooth functions, or $O(\frac{1}{t})$ convergence rates for convex and smooth functions, which match the centralized gradient descent method as well. For example, the EXTRA method in \cite{shi2015extra}, and the class of methods in \cite{xu2015augmented,qu2016harnessing,di2016next,nedich2016achieving}, can all achieve a linear convergence rate for strongly-convex and smooth functions.} Further, {\color{black}among these work, some have additional focus beyond convergence rate, like time-varying graph \cite{nedich2016achieving}, directed graph \cite{zeng2015extrapush,xi2015linear,xi2016add,pu2018push,xin2018linear}, uncoordinated step sizes \cite{xu2015augmented,nedic2016geometrically}, and stochastic gradient \cite{pu2018distributed}.}

It is known that among all centralized gradient based algorithms, centralized Nesterov Gradient Descent (CNGD) \cite{nesterov2013introductory} achieves the optimal convergence rate in terms of first-order oracle complexity. Specifically, for convex and $L$-smooth problems, the convergence rate is $O(\frac{1}{t^2})$, which improves over centralized gradient descent (CGD)'s $O(\frac{1}{t})$ rate; for $\mu$-strongly convex and $L$-smooth problems, the convergence rate is $O((1- \sqrt{\frac{\mu}{L}})^t)$ whose dependence on condition number $\frac{L}{\mu}$ improves over CGD's rate $O((1 - \frac{\mu}{L})^t)$ in the large $\frac{L}{\mu}$ regime. The nice convergence rates naturally lead to the following question: how to achieve similar rate as CNGD in the realm of gradient-based distributed optimization? To the best of our knowledge, for $\mu$-strongly convex and $L$-smooth functions, existing gradient-based distributed methods have the same or worse dependence on condition number compared to CGD. For convex and $L$-smooth functions, paper \cite{jakovetic2014fast} has developed Distributed Nesterov Gradient (D-NG) method and shown that
 it has a convergence rate of $O(\frac{\log t}{t})$.\footnote{The convergence rate for strongly convex and smooth functions is not studied in \cite{jakovetic2014fast}. From our simulation results in Section~\ref{sec:numerical}, the convergence rate for strongly convex and smooth functions is sublinear.} {\color{black} Paper \cite{jakovetic2014fast} also studies an algorithm that achieves almost $O(\frac{1}{t^2})$ convergence rate but it contains an inner loop of multiple consensus steps per iteration. The inner loop places a larger communication burden and needs extra coordination among agents (e.g., agreeing on when to terminate the inner loop), while communication burden has been recognized as the major bottleneck for distributed optimization \cite{lan2017communication}.\footnote{See Section \ref{subsec:innerloop} for a more detailed comparison between our methods and \cite{jakovetic2014fast}, as well as a discussion on why we avoid using the inner loop.} \textcolor{black}{In light of the above discussions, the aim of this paper is to design distributed Nesterov gradient methods that only use one communication step per gradient evaluation but achieve fast convergence rates to bridge the gap between distributed gradient methods and centralized ones, in particular reaching the same level of convergence rate as CNGD in the realm of distributed optimization, for both $\mu$-strongly convex, $L$-smooth functions and convex, $L$-smooth functions.} }

In this paper, we propose an Accelerated Distributed Nesterov Gradient Descent (Acc-DNGD) algorithm. Our algorithm has two versions. The first version Acc-DNGD-SC is designated for $\mu$-strongly convex and $L$-smooth functions. It achieves a convergence rate 
of $O((1-   C (\frac{\mu}{L})^{5/7} )^t)$, where $C$ is a constant that does not depend on $L$ or $\mu$ and the value of $C$ will be clear in Theorem~\ref{thm:str_cvx}.
We emphasize here that the dependence on the condition number $L/ \mu$ in our convergence rate ($[ 1- C (\frac{\mu}{L})^{5/7} ]^t$ ) is strictly better than that of CGD ($(1 - \frac{\mu}{L})^t$), and to the best of our knowledge, \textcolor{black}{the asymptotic dependence on condition number in the large $\frac{L}{\mu}$ regime is the best among all primal \emph{gradient-based} distributed algorithms proposed so far for the class of $\mu$-strongly convex and $L$-smooth functions.}\footnote{\color{black}If dual-gradients are available, the dependence on $ \frac{L}{\mu}$ can be better \cite{scaman2017optimal}. } 
The second version Acc-DNGD-NSC works for convex and $L$-smooth functions. It achieves a $O(\frac{1}{t^{1.4-\epsilon}})$ ($\forall \epsilon\in(0,1.4)$) convergence rate if a vanishing step size is used. We further show that the convergence rate can be improved to $O(\frac{1}{t^2})$ when we use a fixed step size and the objective function is a composition of a linear map and a strongly-convex and smooth function. Both rates are faster than CGD ($O(\frac{1}{t})$). {\color{black}To the best of our knowledge, the $O(\frac{1}{t^{1.4-\epsilon}})$ rate is also the fastest among all existing primal gradient-based distributed algorithms that use one or a constant number of consensus steps per iteration for the class of convex and $L$-smooth functions. \textcolor{black}{In summary, for both function classes, we have achieved rates that are strictly better than CGD but still slower than CNGD, hence we have partially bridged the gap between centralized gradient methods and distributed ones. How to fully close the gap remains an open question.}}

The major technique we use in our algorithm is a gradient estimation scheme \cite{xu2015augmented,di2015distributed,di2016next,qu2016harnessing,nedich2016achieving,xi2016add,nedic2016geometrically}, and we combine the scheme with CNGD.  
{\color{black}	
	Briefly speaking, the scheme tracks the average gradient in an efficient way, and it avoids the use of an inner loop as in \cite{jakovetic2014fast,chen2012thesis}. We will explain the intuition of the gradient estimation scheme in Section~\ref{subsec:innerloop}. } 

{\color{black} Finally, we note that our proof uses the framework of inexact Nesterov method \cite{devolder2014first,devolder2013first}. It is known that inexact Nesterov methods accumulate error and are prone to divergence \cite[Sec. 6]{devolder2014first} \cite[Sec 5]{devolder2013first}. One key step in our proof is to show that when the error has a special structure, we can still get exact convergence. This proof technique may be of independent interest in general inexact Nesterov methods. See Remark \ref{rem:inexact_contribution_sc} and \ref{rem:inexact_contribution_nsc} for more details.   }

\textbf{Notations.} 
In this paper, $n$ is the number of agents, and $N$ is the dimension of the domain of the $f_i$'s. 
Notation $\Vert\cdot\Vert$ denotes $2$-norm for vectors and Frobenius norm for matrices, while $\Vert \cdot\Vert_*$ denotes spectral norm for matrices, $\Vert\cdot\Vert_1$ denotes $1$-norm for vectors, and $\langle \cdot,\cdot\rangle$ denotes inner product for vectors. Notation $\rho(\cdot)$ denotes spectral radius for square matrices, and $\one$ denotes a $n$-dimensional all one column vector. All vectors, when having dimension $N$ (the dimension of the domain of the $f_i$'s), will all be regarded as row vectors. As a special case, all gradients, $\nabla f_i(x)$ and $\nabla f(x)$ are treated as $N$-dimensional row vectors. Notation ``$\leq$'', when applied to vectors of the same dimension, denotes element wise ``less than or equal to''. 

\section{Problem and Algorithm}\label{sec:preliminaries}

\subsection{Problem Formulation}
Consider $n$ agents, $\mathcal{N} = \{1,2,\ldots,n\}$, each of which has a convex function $f_i:\R^N \rightarrow \R$. The objective of distributed optimization is to find $x$ to minimize the average of all the functions, i.e.
\begin{equation}
\min_{x\in \R^N} f(x) \triangleq  \frac{1}{n} \sum_{i=1}^n f_i(x) \label{eq:problem}
\end{equation}
using local communication and local computation. The local communication is defined through a \textit{connected} and \textit{undirected} communication graph $\mathcal{G} = (\mathcal{N},\mathcal{E})$, where the nodes are the agents, and edges $\mathcal{E}\subset \mathcal{V}\times \mathcal{V}$. Agent $i$ and $j$ can send information to each other if and only if $(i,j)\in \mathcal{E}$. The local computation means that each agent can only make its decision based on the local function $f_i$ and the information obtained from its neighbors.  

\textbf{Global Assumptions. }Throughout the paper, we make the following assumptions without explicitly stating them. We assume that each $f_i$ is convex (hence $f$ is also convex). We assume each $f_i$ is $L$-smooth, that is, $f_i$ is differentiable and the gradient is $L$-Lipschitz continuous, i.e., $\forall x,y\in\R^N$, 
$$ \Vert \nabla f_i(x) - \nabla f_i(y)\Vert \leq L \Vert x - y\Vert.$$
As a result, $f$ is also $L$-smooth. We assume $f$ has at least one minimizer $x^*$ with $f(x^*) = f^*$. 

\textbf{Other Assumptions.} We will use the following assumptions in different parts of the paper, and those will be stated explicitly. 

\begin{assumption}\label{assump:str_cvx}
	$\forall i \in \mathcal{N}$, $f_i$ is $\mu$-strongly convex, i.e. $\forall x,y\in \R^N$, we have
	$$f_i(y) \geq f_i(x) + \langle \nabla f_i(x), y-x\rangle + \frac{\mu}{2} \Vert y-x\Vert^2.$$
	As a result, $f$ is also $\mu$-strongly convex.
\end{assumption}

\begin{assumption}\label{assump:compact}
	The set of minimizers of $f$ is compact. 
\end{assumption}

\subsection{Centralized Nesterov Gradient Descent (CNGD)}
We here briefly introduce two versions of centralized Nesterov Gradient Descent (CNGD) algorithm that is derived from \cite[Scheme (2.2.6)]{nesterov2013introductory}.\footnote{The two versions are essentially the same algorithm in \cite[Scheme (2.2.6)]{nesterov2013introductory} with two set of parameters. Describing the two versions in one form as in \cite[Scheme (2.2.6)]{nesterov2013introductory} will need more cumbersome notations. For this reason, we write down the two versions separately.   } The first version, which we term ``CNGD-SC'', is designated for $\mu$-strongly convex and $L$-smooth functions. Given step size $\eta$, let $\alpha = \sqrt{\mu\eta}$. CNGD-SC keeps updating three variables $x(t), v(t), y(t)\in\R^N$. Starting from an initial point $x(0) = v(0) = y(0)\in \R^N$, the updating rule is given by
\begin{subeqnarray}\label{eq:cngd_s}
	x(t+1) &= &y(t) - \eta \nabla f(y(t)) \slabel{eq:cngd_s_1}\\
	v(t+1) &= &(1-\alpha)v(t) + \alpha y(t) - \frac{\eta}{\alpha} \nabla f(y(t))\slabel{eq:cngd_s_2}\\
	y(t+1) &=& \frac{ x(t+1) + \alpha v(t+1)}{1+\alpha}.
\end{subeqnarray}
	
The following theorem (adapted from \cite[Theorem~2.2.1, Lemma 2.2.4]{nesterov2013introductory}) gives the convergence rate of CNGD-SC.
\begin{theorem}\label{thm:cngd}
	Under Assumption \ref{assump:str_cvx},  when $0<\eta\leq\frac{1}{L}$, in CNGD-SC (\ref{eq:cngd_s}) we have $f(x(t)) - f^* = O( (1 - \sqrt{\mu\eta})^t)$.
\end{theorem}

The second version, which we term ``CNGD-NSC'', is designated for convex (possibly not strongly-convex) and $L$-smooth functions. CNGD-NSC keeps updating three variables $x(t), v(t), y(t)\in\R^N$. Starting from an initial point $x(0) = v(0) = y(0)\in \R^N$, the updating rule is given by
\begin{subeqnarray}\label{eq:nsc:cngd}
	x(t+1) &= &y(t) - \eta \nabla f(y(t)) \slabel{eq:nsc:cngd_s_1}\\
	v(t+1) &= &v(t) - \frac{\eta}{\alpha_t} \nabla f(y(t))\slabel{eq:nsc:cngd_s_2}\\
	y(t+1) &=& (1-\alpha_{t+1})x(t+1) + \alpha_{t+1} v(t+1)
\end{subeqnarray}
where $(\alpha_t)_{t=0}^\infty$ is defined by an arbitrarily chosen $\alpha_0\in(0,1)$ and the update equation $\alpha_{t+1}^2 = (1-\alpha_{t+1})\alpha_t^2$. Here $\alpha_{t+1}$ always takes the unique solution in $(0,1)$. 
The following theorem (adapted from \cite[Theorem 2.2.1, Lemma 2.2.4]{nesterov2013introductory}) gives the convergence rate of CNGD-NSC.
\begin{theorem}\label{thm:nsc:cngd}
	In CNGD-NSC (\ref{eq:nsc:cngd}),  when $0<\eta\leq\frac{1}{L}$, we have $f(x(t)) - f^* = O(\frac{1}{t^2} )$.
\end{theorem}

\subsection{Our Algorithm: Accelerated Distributed Nesterov Gradient Descent (Acc-DNGD)} 
We design our algorithm based on a consensus matrix $W= [w_{ij}]\in\R^{n\times n}$. Here $w_{ij}$ stands for how much agent $i$ weighs its neighbor $j$'s information. Matrix $W$ satisfies the following conditions: 
\begin{enumerate}
	\item[(a)] For any $ (i,j)\in \mathcal{E}$, $w_{ij}>0$. For any $ i\in\mathcal{N}$, $w_{ii} > 0$. Elsewhere, $w_{ij} = 0$.
	\item[(b)] Matrix $W$ is doubly stochastic, i.e. $\sum_{i'=1}^n w_{i'j}=\sum_{j'=1}^n w_{ij'}=1$ for all $i, j \in \mathcal{N}$. 
\end{enumerate}
As a result, $\exists \sigma\in(0,1)$ being the second largest singular value of $W$, such that for any $\omega\in \R^{n\times 1}$, we have the ``averaging property'', $\Vert W \omega - \one \bar{\omega}\Vert \leq \sigma  \Vert  \omega - \one\bar{\omega}\Vert$ where $\bar{\omega} = \frac{1}{n} \one^T \omega$ (the average of the entries in $\omega$) \cite{olshevsky2009convergence}. 
How to select a consensus matrix to satisfy these properties has been intensely studied, e.g. \cite{consensus_richard, olshevsky2009convergence}.

Analogous to the centralized case, we present two versions of our distributed algorithm. The first version, Acc-DNGD-SC, is designated for $\mu$-strongly convex and $L$-smooth functions. Each agent keeps a copy of the three variables in CNGD-SC, $x_i(t)$, $v_i(t)$, $y_i(t)$, in addition to a new variable $s_i(t)$ which serves as a gradient estimator. The initial condition is $x_i(0) = v_i(0) = y_i(0)\in\R^{1\times N}$ and $s_i(0) = \nabla f_i(y_i(0))$, and the algorithm updates as follows:
\begin{subeqnarray}\label{eq:alg:update_elementwise}
	x_i(t+1) &= &\sum_{j=1}^n w_{ij} y_j(t) - \eta s_i(t)\slabel{eq:alg:update_elementwise_1}\\
	v_i(t+1) &= &(1-\alpha)\sum_{j=1}^n w_{ij} v_j(t)  + \alpha \sum_{j=1}^n w_{ij} y_j(t) - \frac{\eta}{\alpha} s_i(t) \nonumber\\
	\slabel{eq:alg:update_elementwise_2} \\
	y_i(t+1) &= &  \frac{x_i(t+1) + \alpha v_i(t+1)}{1+\alpha}\slabel{eq:alg:update_elementwise_3} \\
	s_i(t+1) &= & \sum_{j=1}^n w_{ij} s_j(t) + \nabla f_i(y_i(t+1)) - \nabla f_i(y_i(t)) \nonumber \\
	 \slabel{eq:alg:update_elementwise_4}
\end{subeqnarray}
where $[w_{ij}]_{n\times n}$ are the consensus weights, $\eta>0$ is a fixed step size and $\alpha = \sqrt{\mu\eta}$. 


The second version, Acc-DNGD-NSC, is designated for convex (not necessarily strongly convex) and $L$-smooth functions. Similarly as Acc-DNGD-SC, each agent keeps variable $x_i(t)$, $v_i(t)$, $y_i(t)$ and $s_i(t)$. The initial condition is $x_i(0) = v_i(0) = y_i(0)= 0$ and $s_i(0) = \nabla f(0)$,\footnote{We note that the initial condition $s_i(0)=\nabla f(0) = \frac{1}{n} \sum_{i=1}^n \nabla f_i(0)$ requires the agents to conduct an initial run of consensus averaging. We impose this initial condition for technical reasons, while we expect the results of this paper to hold for a relaxed initial condition, $s_i(0) = \nabla f_i(y_i(0))$ which does not need initial coordination. We use the relaxed condition in numerical simulations.}  and the algorithm updates as follows:

\begin{subeqnarray}\label{eq:nsc:update_elementwise}
	x_i(t+1) &= &\sum_{j=1}^n w_{ij} y_j(t) - \eta_t s_i(t)\slabel{eq:nsc:update_elementwise_1}\\
	v_i(t+1) &= &\sum_{j=1}^n w_{ij} v_j(t)  - \frac{\eta_t}{\alpha_t} s_i(t)\slabel{eq:nsc:update_elementwise_2} \\
	y_i(t+1) &= &  (1-\alpha_{t+1})x_i(t+1) + \alpha_{t+1} v_i(t+1)  \slabel{eq:nsc:update_elementwise_3} \\
	s_i(t+1) &= & \sum_{j=1}^n w_{ij} s_j(t) + \nabla f_i(y_i(t+1)) - \nabla f_i(y_i(t)) \nonumber\\
	  \slabel{eq:nsc:update_elementwise_4}
\end{subeqnarray}
where $[w_{ij}]_{n\times n}$ are the consensus weights and $\eta_t\in (0, \frac{1}{L} )$ are the step sizes. Sequence $(\alpha_t)_{t\geq 0}$ is generated as follows. First we let $\alpha_0 = \sqrt{\eta_0 L} \in(0,1)$. Then given $\alpha_t\in (0,1)$, we select $\alpha_{t+1}$ to be the unique solution in $(0,1)$ of the following equation,\footnote{{\color{black}Without causing any confusion with the $\alpha_t$ in (\ref{eq:nsc:cngd}), in the rest of the paper we abuse the notation of $\alpha_t$.} }  
\begin{equation}
\alpha_{t+1}^2 =  \frac{\eta_{t+1}}{\eta_t} (1-\alpha_{t+1})\alpha_t^2.\label{eq:nsc:alpha_update}
\end{equation}
We will consider two variants of the algorithm with the following two step size rules.
\begin{itemize}
	\item \textbf{Vanishing step size:} $\eta_t = \eta \frac{1}{(t+t_0)^\beta}$ for some $\eta\in(0,\frac{1}{L})$, $\beta\in (0,2)$ and $t_0\geq 1$.
	\item \textbf{Fixed step size:} $\eta_t = \eta>0$.
\end{itemize}

In both versions (Acc-DNGD-SC (\ref{eq:alg:update_elementwise}) and Acc-DNGD-NSC (\ref{eq:nsc:update_elementwise})), because $w_{ij}=0$ when $(i,j) \notin \mathcal{E}$,  each node $i$ only needs to send $x_i(t)$, $v_i(t)$, $y_i(t)$ and $s_i(t)$ to its neighbors. Therefore, the algorithm can be operated in a fully distributed fashion with only local communication. The additional term $s_i(t)$ allows each agent to obtain an estimate on the average gradient $\frac{1}{n}\sum_{i=1}^n \nabla f_i(y_i(t))$. Compared with distributed algorithms without this estimation term, it improves the convergence rate. We will provide intuition behind the gradient estimator $s_i(t)$ in Sec. \ref{subsec:innerloop}. Because of the use of the gradient estimation term, we call this method as \textit{Accelerated} Distributed Nesterov Gradient Descent (Acc-DNGD) method.  

\subsection{Convergence of the Algorithm}
To state the convergence results, we need to define the following average sequence $\bar{x}(t) = \frac{1}{n} \sum_{i=1}^n x_i(t) \in \R^{1\times N}$. We first provide the convergence result for Acc-DNGD-SC in Theorem~\ref{thm:str_cvx}. 
\begin{theorem}\label{thm:str_cvx}
Consider algorithm Acc-DNGD-SC (\ref{eq:alg:update_elementwise}). Under the strongly convex assumption (Assumption \ref{assump:str_cvx}), when $0< \eta< \frac{\sigma^{3}(1-\sigma)^3 }{250^2  L} (\frac{\mu}{L})^{3/7}$, we have (a) $f(\bar{x}(t)) - f^* = O((1 - \sqrt{\mu\eta})^t)$; {\color{black}(b) For any $i$, $ f( y_i(t)) -  f^*  = O((1- \sqrt{\mu\eta})^{t})$.}
\end{theorem}
 
{\color{black}The upper bound on the step size $\eta$ given in theorem \ref{thm:str_cvx}} results in a convergence rate of $O([1 - C (\frac{\mu}{L})^{5/7} ]^t)$ with $C = \frac{\sigma^{1.5}(1-\sigma)^{1.5}}{250}$. In this convergence rate, the dependence on the condition number\footnote{The quantity $\frac{L}{\mu}$ is called the condition number of the objective function. For a convergence rate of \textcolor{black}{$ O((1 - (\frac{\mu}{L})^\delta )^t)$, it takes $O ( (\frac{L}{\mu})^\delta \log \frac{1}{\epsilon})$} iterations to reach an accuracy level $\epsilon$. Therefore, the smaller $\delta $ is, the better the convergence rate is.} $\frac{L}{\mu}$ is strictly better than that of CGD ($O((1-\frac{\mu}{L})^t)$), and hence CGD based distributed algorithms \cite{shi2015extra,qu2016harnessing}. 
  \textcolor{black}{This means that when the condition number $\frac{L}{\mu}$ is sufficiently large, our method can outperform CGD and CGD-based distributed methods.} This is particularly appealing because in many machine learning applications, the condition number $\frac{L}{\mu}$ can be as large as the sample size \cite[Section 3.6]{bubeck2014convex}.  We also highlight that in the conference version of this paper \cite{qu2016accelerated}, the dependence on condition number was $(\frac{\mu}{L})^{1.5}$, worse than that of CGD. Compared to \cite{qu2016accelerated}, this paper conducts a sharper analysis on our algorithm. 
  

 We next provide the convergence result for Acc-DNGD-NSC in Theorem~\ref{thm:nsc:vanishing}.
\begin{theorem}\label{thm:nsc:vanishing}
	Consider algorithm Acc-DNGD-NSC (\ref{eq:nsc:update_elementwise}). Suppose Assumption \ref{assump:compact} is true and without loss of generality we assume $\bar{v}(0) \neq x^*$ {\color{black} where  $\bar{v}(0) = \frac{1}{n}\sum_{j=1}^n v_j(0)$}.  Let the step size be $\eta_t = \frac{\eta}{(t+t_0)^\beta}$ with $\beta= 0.6 + \epsilon$ where $\epsilon\in (0,1.4)$. Suppose the following conditions are met.
	\begin{itemize}
		\item[(i)] 
		$$t_0 > \frac{1}{ \min( ( \frac{\sigma+3}{\sigma+2} \frac{3}{4})^{\sigma/(28\beta)} , (\frac{16}{15+\sigma})^{\frac{1}{\beta}})  - 1  }.$$
		\vspace*{-\baselineskip}
		\item[(ii)] $$\eta < \min(\frac{\sigma^2}{9^3 L},   \frac{(1-\sigma)^3}{6144 L}).$$
		\vspace*{-\baselineskip}
		\item[(iii)] $$\eta< \Bigg( \frac{ D(\beta,t_0) (\beta - 0.6) (1-\sigma)^2}{9216 (t_0+1)^{2-\beta} L^{2/3} [4 + R^2/\Vert\bar{v}(0) - x^*\Vert^2 ]}\Bigg)^{3/2}$$ where $D(\beta,t_0) = \frac{1}{(t_0+3)^2 e^{16+ \frac{6}{2-\beta} }}$ and $R$ is the diameter of the $(2 f(\bar{x}(0)) - f^* + 2L\Vert\bar{v}(0) - x^*\Vert^2)$-level set of $f$.\footnote{Here we have used the fact that by Assumption \ref{assump:compact}, all level sets of $f$ are bounded. See \cite[Proposition B.9]{bertsekas1999nonlinear}.}
	\end{itemize}
	
	{\color{black} Then, we have (a) $f(\bar{x}(t)) - f^* =O( \frac{1}{t^{2-\beta}})= O(\frac{1}{t^{1.4 - \epsilon }})$; (b) $\forall i$, $f(y_i(t)) - f^* = O(\frac{1}{t^{1.4 - \epsilon }})$. }
\end{theorem}
\begin{remark}\label{rem:vanishing_stepsize}
{\color{black} The step size condition in Theorem \ref{thm:nsc:vanishing} may be difficult to implement since constant $R$ may be unknown to the individual agents. However, we believe the exact value of $\eta>0$ and $t_0>0$ do not actually matter, and in simulations we simply set $\eta = \frac{1}{2L}$ and $t_0=1$. 
The reason is as follows. The most important condition we need on $\eta_t$ is that $\eta_t$ is monotonically decaying to $0$ in the order of $O(\frac{1}{t^{\beta}})$ with $\beta\in(0.6,2)$. The condition $\beta>0.6$ is important in the proof as it makes $\eta_t$ converge to $0$ sufficiently fast to control the error.
Other than that, the condition (i)(ii)(iii) on $\eta$ and $t_0$ in Theorem~\ref{thm:nsc:vanishing} is to ensure $\eta_0$ is small, and $\frac{\eta_t}{\eta_{t+1}}$ is close to $1$, which are needed only for some technical reasons in the proof.
In fact, in the literature that uses vanishing step sizes, it is not uncommon that only the order of decaying matters, while other constants are not as important \cite[Theorem 5 (b)]{jakovetic2014fast}. Lastly, in Figure \ref{fig:stability_rem1} in the simulations, we verify the stability of step size rule $\frac{1}{2L(t+1)^\beta}$ (i.e. $\eta = \frac{1}{2L}$ and $t_0=1$) under randomly generated problem instances.
}
\end{remark}

While in Theorem \ref{thm:nsc:vanishing} we require $\beta > 0.6$, we conjecture that the algorithm will converge with rate $O(\frac{1}{t^{2-\beta}}) $ even if we choose $\beta\in[0,0.6]$, with $\beta = 0$ corresponding to the case of fixed step size. In Section \ref{sec:numerical} we will use numerical methods to support this conjecture. 

In the next theorem, we provide a $O(\frac{1}{t^2})$ convergence result  when a fixed step size is used and the objective functions belong to a special class. 
\begin{theorem}\label{thm:nsc:fixed}
Consider algorithm Acc-DNGD-NSC (\ref{eq:nsc:update_elementwise}). Assume each $f_i(x)$ can be written as $f_i(x) = h_i(x A_i ) $, where $A_i$ is a non-zero $N\times M_i$ matrix, and $h_i(x):\R^{1\times M_i}\rightarrow \R$ is a $\mu_0$-strongly convex and $L_0$-smooth function. Suppose we use the fixed step size rule $\eta_t = \eta$, with 
		$$0<\eta < \min(\frac{\sigma^2}{9^3 L} , \frac{\mu^{1.5} (1-\sigma)^3}{L^{2.5} 6912^{1.5}})$$
		where $L = L_0\nu$ with $\nu = \max_i \Vert A_i\Vert_*^2$; and $\mu = \mu_0 \gamma$ with $\gamma$ being the smallest non-zero eigenvalue of matrix $A = \frac{1}{n} \sum_{i=1}^n A_i A_i^T$. Then, we have {\color{black}(a) $f(\bar{x}(t)) - f^* = O(\frac{1}{t^2})$; (b) $\forall i$, $f(y_i(t)) - f^* = O(\frac{1}{t^2})$. }
\end{theorem}
{\color{black}\begin{remark}
	The reason we can prove a faster convergence rate under the assumption in Theorem~\ref{thm:nsc:fixed} is that, when $f_i = h_i(x A_i)$, we have $\nabla f_i = \nabla h_i (x A_i) A_i^T$ which lies in the row space of $A = \frac{1}{n}\sum_{i=1}^n A_iA_i^T$.\footnote{\color{black}To see this, let (column) vector $y$ be such that $Ay =0$. Then we have $\langle y, Ay\rangle =\frac{1}{n}\sum_{j=1}^n \Vert A_j^T y\Vert^2 = 0 $, implying $A_i^T y=0$. Hence, $\langle  \nabla f_i(x), y\rangle =  \nabla h_i(xA_i) A_i^T y = 0$, implying $\nabla f_i(x)$ is orthogonal to any vector in $A$'s (right) kernel space.} As a result, $x_i(t) - y_i(t)$ also lies within the row space of $A$ as $x_i(t) - y_i(t)$ is a linear combination of such gradients.  With this property, and the fact that $h_i$ is strongly convex, we can show $f$ behaves like a strongly convex function around $x_i(t)$ and $y_i(t)$. It is this ``local strong convexity'' property that helps us obtain a faster convergence result.
\end{remark}}
An important example of functions $f_i(x)$ in Theorem~\ref{thm:nsc:fixed} is the square loss for linear regression (cf. Case I in Section~\ref{sec:numerical}) when the sample size is less than the parameter dimension.

{\color{black}\begin{remark}
		Here we provide a few remarks on the strongly-convex parameter $\mu$, smooth parameter $L$, graph parameter $\sigma$, and the numerical constants in the step size bounds in Theorem~\ref{thm:str_cvx},\ref{thm:nsc:vanishing}, and \ref{thm:nsc:fixed}.
		\begin{itemize}
			\item Parameter $\mu$, $L$: though we have assumed each function $f_i$ has the same strongly convex and smooth parameter $\mu, L$, in reality each function $f_i$ may have heterogeneous $\mu_i$ and $L_i$ parameters. In this case, our results still hold if we take $L$ to be an upper bound of the $L_i$'s, and $\mu$ to be a lower bound of the $\mu_i$'s. This may require additional coordination among the agents before running the distributed algorithms. Recently, techniques have been proposed that allow each $f_i$ to have heterogeneous $\mu_i$ and $L_i$, and each agent only needs to know their local $\mu_i, L_i$ as opposed to the global $L$, $\mu$ \cite{nedic2016geometrically}. How to incorporate such techniques into our method remains our future work.
			\item Parameter $\sigma$: though the graph parameter $\sigma$ is not directly available to each agent, there have been techniques in the literature to replace the $\sigma$ in the step size bounds with some locally accessible surrogates when the agents have an estimate on the total number of agents in the network \cite[Corollary 5]{qu2016harnessing}.
			\item Numerical Constants:  all numerical constants in the step size bounds (like the $250^2$ in Theorem~\ref{thm:str_cvx})  provided in this section are conservative. In the proof, we have been intentionally loose on the numerical constants to keep the analysis simple and only focus on getting the tightest asymptotic dependence on critical parameters (e.g. the $(\frac{\mu}{L})^{3/7}$ in the step size bound in Theorem~\ref{thm:str_cvx}). 
			In the numerical simulations, we use much larger step sizes than the bounds in this section. How to find tighter bounds remains our future work.
		\end{itemize}
\end{remark}}

\section{Algorithm Development} \label{sec:algo_dev}
{\color{black}In this section, we will discuss the proof idea in Section \ref{subsec:proofidea} and then compare our algorithm to a class of methods that use an inner-loop of consensus iterations \cite{jakovetic2014fast,chen2012fast} in Section \ref{subsec:innerloop}. Specifically, we will explain why we avoid the use of the inner loop idea and instead use the gradient estimator $s_i(t)$ in our algorithm \eqref{eq:alg:update_elementwise} \eqref{eq:nsc:update_elementwise}.}
\subsection{Inexact Nesterov gradient methods}\label{subsec:proofidea}
{\color{black} We now explain the intuition behind Acc-DNGD-SC \eqref{eq:alg:update_elementwise} whereas the intuition for Acc-DNGD-NSC \eqref{eq:nsc:update_elementwise} is similar. The main idea is that the distributed algorithm imitates the centralized one. Since we conduct averaging steps (``$\sum_{j=1}^n w_{ij} y_j(t)$'', ``$\sum_{j=1}^n w_{ij} v_j(t)$'' and ``$\sum_{j=1}^n w_{ij} s_j(t)$'') in \eqref{eq:alg:update_elementwise_1}, \eqref{eq:alg:update_elementwise_2}, \eqref{eq:alg:update_elementwise_3}, we expect that $x_i(t) \approx \bar{x}(t):= \frac{1}{n} \sum_{j=1}^n x_j(t) $, $y_i(t) \approx \bar{y}(t):= \frac{1}{n} \sum_{j=1}^n y_j(t) $, $v_i(t) \approx \bar{v}(t):= \frac{1}{n} \sum_{j=1}^n v_j(t) $, and also $s_i(t) \approx \bar{s}(t):= \frac{1}{n} \sum_{j=1}^n s_j(t)$ which can easily shown to be equal to  $\frac{1}{n} \sum_{j=1}^n \nabla f_j(y_j(t)) \approx \nabla f(\bar{y}(t))$. Therefore, we expect \eqref{eq:alg:update_elementwise} to imitate the following update,\footnote{This ``imitation'' argument will be made more precise using inexact Nesterov Gradient Descent framework in Lemma~\ref{lem:inexact_grad}.}
		{\small	\begin{subeqnarray}\label{eq:nes:imitate}
				\bar{x}(t+1) &= &  \bar{y}(t) - \eta \nabla f(\bar{y}(t))  \\
				\bar{v}(t+1) &= &(1-\alpha) \bar{v}(t) + \alpha \bar{y}(t) - \frac{\eta}{\alpha} \nabla f(\bar{y}(t)) \\
				\bar{y}(t+1) &= &  \frac{ \bar{x}(t+1) + \alpha \bar{v}(t+1)}{1+\alpha}  
		\end{subeqnarray}}which is precisely CNGD-SC \eqref{eq:cngd_s}. 
	However, we emphasize that the distributed algorithm does not exactly follow \eqref{eq:nes:imitate} due to consensus errors, i.e. how close each individual's $x_i(t)$, $y_i(t)$, $v_i(t)$, $s_i(t)$ are to their average, $\bar{x}(t)$, $\bar{y}(t)$, $\bar{v}(t)$, $\bar{s}(t)$. With these observations, the proof can be mainly divided into two steps. 
	The first is to bound the consensus errors. It turns out that bounding $\Vert y_i(t) -  \bar{y}(t) \Vert $ is sufficient for our purpose (cf. Lemma \ref{lem:nes:rel_con_err} and \ref{lem:nsc:rel_con_err}). Secondly, with the bounded consensus error, we treat the distributed algorithm as an inexact version of \eqref{eq:nes:imitate} and then use the framework of inexact Nesterov gradient descent \cite{devolder2013first,devolder2014first} to prove convergence.
	
\begin{remark} It is well known Nesterov Gradient method is prone to accumulating noise and divergence \cite{devolder2013first,devolder2014first}. Therefore, in the literature, when adapting Nesterov Gradient Descent to inexact settings, it usually requires nontrivial error bounds and sophisticated proof techniques to show convergence (e.g. in stochastic setting \cite{allen2017katyusha}). From this perspective, our contribution can be viewed as adapting Nesterov Gradient Descent to a distributed setting and show its convergence, which is nontrivial for the reasons above. Further, our proof technique may be of independent interest in general inexact Nesterov methods. See Remark \ref{rem:inexact_contribution_sc} and \ref{rem:inexact_contribution_nsc} for more details.  
	\end{remark}
  }
\subsection{Comparison with Inner-Loop Methods and Intuition on Gradient Estimator $s_i(t)$} \label{subsec:innerloop}
{\color{black}This section will compare this work and the methods that use an inner loop of consensus steps \cite{jakovetic2014fast,chen2012fast}, in particular \cite{jakovetic2014fast}. Ref. \cite{jakovetic2014fast} proposes two algorithms, the first of which is named ``D-NG'', which is described as follows \cite[eq. (2) (3)]{jakovetic2014fast} using the notations of this paper (cf. \eqref{eq:nes:stack_notation} Section~\ref{sec:proofoverview}), 
{\small	\begin{subeqnarray}\label{eq:jakovetic_dng}
	x (t+1)&=& W y(t) - \frac{c}{t+1} \nabla (t) \slabel{eq:jakovetic_dng_x} \\
	y (t+1) &= & x(t+1) + \frac{t}{t+3} (x(t+1) - x(t))
	\end{subeqnarray}	}where $x(t)$, $y(t)$, $\nabla(t)$ are the local quantities $x_i(t)$, $y_i(t)$, $\nabla f_i(y_i(t))$ stacked together as defined in \eqref{eq:nes:stack_notation}. 

 The problem with \eqref{eq:jakovetic_dng} is that the descent direction for $x_i(t)$ is local gradient $\nabla f_i(y_i(t))$, which is incorrect. To see this, note when $x_i(t)$ and $y_i(t)$ reach an optimizer $x^*$ of $f$, $\nabla f_i(y_i(t))$ is in general non-zero. As a result, $x_i(t) = y_i(t) = x^*$ may not even be a fixed point of $\eqref{eq:jakovetic_dng}$. On the contrary, the correct direction is the average gradient, i.e. $\frac{1}{n}\sum_{j=1}^n \nabla f_j(y_j(t))$, which will be zero when $y_i(t) = x^*$. As a result of the incorrect descent direction, a $O(\frac{1}{t})$ step size is used in \eqref{eq:jakovetic_dng} to force $x^*$ be a limit point of algorithm (\ref{eq:jakovetic_dng}). However, the diminishing step size also causes the convergence rate of \eqref{eq:jakovetic_dng} to be  $O(\frac{\log t}{ t})$, slower than $O(\frac{1}{t^2})$ of CNGD.

The inner-loop approach is proposed to fix the incorrect descent direction issue. Specifically, \cite{jakovetic2014fast} proposes a second algorithm named ``D-NC'' \cite[eq. (9)(10)]{jakovetic2014fast} as follows,
{\small	\begin{subeqnarray}\label{eq:jakovetic_dnc}
		x (t+1) &=& W^{\tau_x(t)} [ y(t) - \eta  \nabla (t) ] \slabel{eq:jakovetic_dnc_x} \\
		y (t+1) &= &W^{\tau_y(t)} [x(t+1) + \frac{t}{t+3}   (x(t+1) - x(t)) ] 
\end{subeqnarray}	}where $\tau_x(t)$ and $\tau_y(t)$ are both integer functions of $t$ that grow in the order of $\Omega(\log t)$. 
Briefly speaking, at each iteration $t$, D-NC first moves along the (incorrect) local gradient direction $\nabla f_i(y_i(t))$ in \eqref{eq:jakovetic_dnc_x}, and then conducts an inner loop of $\tau_x(t)$ consensus steps (the multiplication of matrix $W^{\tau_x(t)}$ in \eqref{eq:jakovetic_dnc_x}). This effectively make the descent direction \emph{approximately} the average gradient $\frac{1}{n}\sum_{j=1}^n\nabla f_j(y_j(t))$, and hence the descent direction is approximately correct, and a constant step size can be used, achieving a $O(\frac{1}{T^{2-\epsilon}})$ (for arbitrarily small $\epsilon>0$) rate where $T$ is the total number of consensus steps.

{\color{black}
 Despite the improved convergence rate, there are a few disadvantages of the inner-loop approach. Firstly, the implementation of the inner loop places extra coordination burden on the agents, e.g. they need to agree on when to terminate the inner-loop. Secondly, the work in \cite{jakovetic2014fast} (also in \cite{chen2012fast}) requires to assume the gradients $\nabla f_i(\cdot)$ are uniformly bounded.
 This is a strong assumption since even the simple quadratic functions do not meet this assumption. Further, as \cite[Sec VII-B]{jakovetic2014fast} shows, without the bounded gradient assumption the results in \cite{jakovetic2014fast} no longer holds. Thirdly, the inner loop approach does not work well when the cost functions are strongly convex. Though neither \cite{jakovetic2014fast} nor \cite{chen2012fast} studies the strongly convex case, it can be shown that if applying the inner loop approach to the strongly convex case, at outer iteration $t$, to obtain a sufficiently accurate estimate of average gradient, the number of inner loop iterations needs to grow at least in the order of $\Omega(t)$. The large number of inner loop iterations slows down the convergence, resulting in a sublinear convergence rate. 
 
For the above reasons, in this paper we do not use an inner loop. Instead, we use a gradient estimation sequence recently proposed in \cite{xu2015augmented,qu2016harnessing,nedich2016achieving}, which is a more efficient use of consensus steps to get the correct descent direction. Briefly speaking, our algorithm uses descent direction $s_i(t)$, which has its own update formula \eqref{eq:alg:update_elementwise_4} \eqref{eq:nsc:update_elementwise_4}, and can be shown to asymptotically converge to the average gradient, the correct descent direction. 	
The way it works is that, at each time $t$, $s_i(t+1)$ averages over neighbors' $s_j(t)$, which itself is an estimate of the average gradient at time $t$, and then $s_i(t+1)$ adds ``incremental information'', i.e. the difference between the new gradient and the previous gradient ($\nabla f_i(y_i(t+1)) - \nabla f_i(y_i(t))$). The incremental information makes sure the estimator tracks the latest gradient, and the incremental information is expected to be ``small'', since $\Vert \nabla f_i(y_i(t+1)) - \nabla f_i(y_i(t))\Vert\leq L \Vert y_i(t+1) - y_i(t)\Vert $ and we can make $\Vert y_i(t+1)-y_i(t)\Vert $ small by carefully controlling the step size. Given the small incremental information, $s_i(t)$ asymptotically gives a good estimate of the average gradient, which allows us to use a larger step size than \eqref{eq:jakovetic_dng} and improve the convergence rate. The use of such gradient estimator has gained popularity recently and our earlier work \cite{qu2016harnessing} provides more details explaining how it works. Using this gradient estimator $s_i(t)$, we have avoided the disadvantages associated with the inner loop approach. Our algorithms are easier to implement, do not need the bounded gradient assumption, and achieve linear convergence rates for strongly convex and smooth functions. {\color{black} This being said, we acknowledge that for convex and smooth functions, our algorithm's rate $O(\frac{1}{t^{1.4-\epsilon}})$ (Theorem~\ref{thm:nsc:vanishing}) is slower than that of the inner-loop approach D-NC \cite{jakovetic2014fast}'s rate $O(\frac{1}{t^{2-\epsilon}})$, as our method requires smaller step size (Theorem~\ref{thm:nsc:vanishing}). With the additional assumption on the functional form, our method can also achieve a $O(\frac{1}{t^{2}})$ rate (Theorem~\ref{thm:nsc:fixed}) matching that of \cite[D-NC]{jakovetic2014fast}. Though we need the additional assumption, we comment that \cite[D-NC]{jakovetic2014fast} relies on the bounded gradient assumption which we do not need. To sum up, for the class of convex and smooth functions, the comparison between our proposed method and D-NC is complicated - one should take into account the type of functions and whether they meet the assumptions in our paper or the ones in \cite[D-NC]{jakovetic2014fast}, and one should also consider the complexity in implementation - e.g. the inner-loop approach would require additional coordination including agreeing on when to switch between the inner-loop and the outer-loop. 
  }}

}

	\section{Convergence Analysis of Acc-DNGD-SC}\label{sec:convergence_sc}
In this section, we will provide the proof of Theorem~\ref{thm:str_cvx}. We will first provide a proof overview in Section~\ref{sec:proofoverview} and then defer a few steps of the proof to the rest of the section. 

\subsection{Proof Overview} \label{sec:proofoverview}
Without causing any confusion with notations in (\ref{eq:cngd_s}), in this section we abuse the use of notation $x(t), v(t), y(t)$. We introduce notations $x(t)$, $v(t)$, $y(t)$, $s(t)$, $\nabla(t)\in\R^{n\times N}$ as follows, 
{\small\begin{subeqnarray}\label{eq:nes:stack_notation}
	  x(t) &=& [x_1(t)^T,x_2(t)^T,\ldots,x_n(t)^T]^T \\
	  v(t) &=& [v_1(t)^T,v_2(t)^T,\ldots,v_n(t)^T]^T\\
	y(t) &= &[y_1(t)^T,y_2(t)^T,\ldots,y_n(t)^T]^T  \\
	   s(t) &=& [s_1(t)^T,s_2(t)^T,\ldots,s_n(t)^T]^T\\
	\nabla(t) &= &[\nabla f_1(y_1(t))^T, \nabla f_2(y_2(t))^T, \ldots, \nabla f_n(y_n(t))^T]^T.
\end{subeqnarray}}

Now the algorithm Acc-DNGD-SC (\ref{eq:alg:update_elementwise}) can be written as
{\small\begin{subeqnarray}\label{eq:nes:update_vector}
	x(t+1) &= & W y(t) - \eta s(t)\slabel{eq:nes:update_vector_a}\\
	v(t+1) &=& (1-\alpha)W v(t) + \alpha W y(t) - \frac{\eta}{\alpha} s(t) \slabel{eq:nes:update_vector_b}\\
	y(t+1) &= &  \frac{x(t+1) +  \alpha v(t+1)}{1+\alpha}  \slabel{eq:nes:update_vector_c}\\
	s(t+1) &= & W s(t) + \nabla(t+1) - \nabla (t).\slabel{eq:nes:update_vector_d}
\end{subeqnarray}}Apart from the average sequence $\bar{x}(t) =\frac{1}{n} \sum_{i=1}^n x_i(t)$ that we have defined, we also define several other average sequences, $\bar{v}(t) =\frac{1}{n} \sum_{i=1}^n v_i(t)$, $\bar{y}(t) = \frac{1}{n} \sum_{i=1}^n y_i(t)$, $\bar{s}(t) = \frac{1}{n}  \sum_{i=1}^n s_i(t)$, and $g(t) =\frac{1}{n} \sum_{i=1}^n \nabla f_i(y_i(t))$. {\color{black} We would like to remind here that all quantities associated with each agent like $x_i(t)$, $v_i(t)$, $y_i(t)$, $s_i(t)$, $\nabla f_i(y_i(t))$, as well as their average $\bar{x}(t)$, $\bar{v}(t)$, $\bar{y}(t)$, $\bar{s}(t)$, $g(t)$, are \emph{row} vectors of dimension $N$ where $N$ is the dimension of the domain of the $f_i$'s. As a result, quantity $x(t), v(t), y(t), s(t)$, $\nabla(t)$, are matrices of dimension $n$-by-$N$. }

\vspace{6pt}
\noindent
\textbf{Overview of the Proof. }In our proof, we firstly derive the update formula for the average sequences (Lemma \ref{lem:nes:ave}). Then, it turns out that the update rule for the average sequences is in fact CNGD-SC (\ref{eq:alg:update_elementwise}) with inexact gradients \cite{devolder2014first}, and the inexactness is characterized by ``consensus error'' $\Vert y(t) - \one\bar{y}(t)\Vert$ (Lemma \ref{lem:inexact_grad}). The consensus error is then bounded in Lemma \ref{lem:nes:rel_con_err}. With the bound on consensus error, we can roughly apply the same proof steps of CNGD-SC (see e.g. \cite{nesterov2013introductory}) to the average sequence and finish the proof of Theorem \ref{thm:str_cvx}.

\begin{lemma}\label{lem:nes:ave} The following equalities hold.
{\small	\begin{subeqnarray}\label{eq:nes:update_ave}
		\bar{x}(t+1) &= &  \bar{y}(t) - \eta g(t) \slabel{eq:nes:update_ave_a}\\
		\bar{v}(t+1) &= &(1-\alpha) \bar{v}(t) + \alpha \bar{y}(t) - \frac{\eta}{\alpha} g(t) \slabel{eq:nes:update_ave_b}\\
		\bar{y}(t+1) &= &  \frac{ \bar{x}(t+1) + \alpha \bar{v}(t+1)}{1+\alpha} \slabel{eq:nes:update_ave_c}\\
		\bar{s}(t+1) &= & \bar{s}(t) + g(t+1) -  g(t) = g(t+1)\slabel{eq:nes:update_ave_d}
	\end{subeqnarray}}
\end{lemma}
\noindent\textit{Proof: } We omit the proof since these can be easily derived using the fact that $W$ is doubly stochastic. For (\ref{eq:nes:update_ave_d}) we also need to use the fact that $\bar{s}(0)=g(0)$. \qedd 

From (\ref{eq:nes:update_ave_a})-(\ref{eq:nes:update_ave_c}) we see that the sequences $\bar{x}(t)$, $\bar{v}(t)$ and $\bar{y}(t)$ follow a update rule  similar to the CNGD-SC in (\ref{eq:cngd_s}). The only difference is that $g(t)$ in (\ref{eq:nes:update_ave_a})-(\ref{eq:nes:update_ave_c})  is not the exact gradient $\nabla f(\bar{y}(t))$ in CNGD-SC. In the following Lemma, we show that $g(t)$ is an inexact gradient with error $O(\Vert y(t) - \one\bar{y}(t)\Vert ^2)$.
\begin{lemma}\label{lem:inexact_grad} Under Assumption \ref{assump:str_cvx}, $\forall t$, $g(t)$ is an inexact gradient of $f$ at $\bar{y}(t)$ with error $O(\Vert y(t) - \one\bar{y}(t)\Vert ^2)$ in the sense that if we let $ \hat{f}(t) = \frac{1}{n}\sum_{i=1}^n \big[f_i(y_i(t))+ \langle \nabla f_i(y_i(t)), \bar{y}(t) - y_i(t)\rangle \big]$, then $\forall \omega\in\R^N$,
\footnote{The exact gradient $\nabla f(\bar{y}(t))$, satisfies $f(\bar{y}(t)) + \langle \nabla f(\bar{y}(t)), \omega - \bar{y}(t)\rangle + \frac{\mu }{2}\Vert \omega - \bar{y}(t)\Vert^2 \leq f(\omega) \leq f(\bar{y}(t)) + \langle \nabla f(\bar{y}(t)), \omega - \bar{y}(t)\rangle + \frac{L }{2}\Vert \omega - \bar{y}(t)\Vert^2$. This is why $g(t)$ is called inexact gradient.}
{\small\begin{align} 
f(\omega) &\geq \hat{f} (t) + \langle g(t), \omega - \bar{y}(t)\rangle + \frac{\mu }{2}\Vert \omega - \bar{y}(t)\Vert^2 \label{eq:nes:inexact:noisy_oracle_1}\\
f(\omega) &\leq \hat{f}(t) + \langle g(t), \omega - \bar{y}(t)\rangle + L \Vert \omega - \bar{y}(t)\Vert^2  \nonumber\\
&\qquad   + L  \frac{1}{n} \Vert y(t) - \one\bar{y}(t)\Vert^2. \label{eq:nes:inexact:noisy_oracle_2}
\end{align}}
\end{lemma}

\noindent\textit{Proof: } 
For any $\omega\in\R^N$, we have
{\small\begin{align*}
&f(\omega) = \frac{1}{n} \sum_{i=1}^n f_i(\omega) \\
& \geq \frac{1}{n} \sum_{i=1}^n \big [ f_i(y_i(t)) + \langle \nabla f_i(y_i(t)), \omega-y_i(t) \rangle + \frac{\mu}{2} \Vert \omega - y_i(t)\Vert^2 \big] \\
&=  \frac{1}{n} \sum_{i=1}^n \big [ f_i(y_i(t)) + \langle \nabla f_i(y_i(t)), \bar{y}(t)-y_i(t) \rangle \big]  \\
&\qquad  +  \frac{1}{n} \sum_{i=1}^n  \langle \nabla f_i(y_i(t)), \omega - \bar{y}(t)\rangle   + \frac{1}{n} \sum_{i=1}^n\frac{\mu}{2} \Vert \omega - y_i(t)\Vert^2  \\
&\geq \hat{f}(t) + \langle g(t), \omega - \bar{y}(t)\rangle + \frac{\mu}{2} \Vert \omega - \bar{y}(t)\Vert^2
\end{align*}}which shows (\ref{eq:nes:inexact:noisy_oracle_1}). For (\ref{eq:nes:inexact:noisy_oracle_2}), similarly,
{\footnotesize \begin{align*}
&f(\omega)\\
&\leq \frac{1}{n} \sum_{i=1}^n \big [ f_i(y_i(t)) + \langle \nabla f_i(y_i(t)), \omega-y_i(t)  + \frac{L}{2} \Vert \omega - y_i(t)\Vert^2\big] \\
&= \hat{f}(t) + \langle g(t), \omega - \bar{y}(t)\rangle    +\frac{L}{2}  \frac{1}{n} \sum_{i=1}^n \Vert \omega - y_i(t) \Vert^2\\
&\leq  \hat{f}(t) + \langle g(t), \omega - \bar{y}(t)\rangle + L \Vert \omega - \bar{y}(t)\Vert^2 +L \frac{1}{n} \sum_{i=1}^n \Vert \bar{y}(t) - y_i(t)\Vert^2
\end{align*}}where in the last inequality we have used the elementary fact that $\Vert u + v\Vert^2\leq 2\Vert u\Vert^2 + 2\Vert v\Vert^2$ for all $u,v\in\R^N$. Then, (\ref{eq:nes:inexact:noisy_oracle_2}) follows by $\sum_{i=1}^n \Vert \bar{y}(t) - y_i(t)\Vert^2=\Vert y(t) - \one\bar{y}(t)\Vert^2 $.\qedd \\

The consensus error $\Vert y(t) - \one\bar{y}(t)\Vert$ in the previous lemma is bounded by the following lemma whose proof is deferred to Section~\ref{subsec:con_err}.
\begin{lemma}\label{lem:nes:rel_con_err}
	When $0<\eta < \min(\frac{1}{L} \frac{(1-\sigma)^3}{512},\frac{1}{L} \frac{\sigma^3}{64}   ) $, we have {\small $$\Vert y(k) - \one\bar{y}(k)\Vert\leq A_1(\eta) \theta^k + A_2(\eta) \sum_{\ell=0}^{k-1} \theta^{k-1-\ell} a(\ell)$$}
where $a(\ell) \triangleq  \Vert  \bar{y}(\ell)-\bar{x}(\ell)\Vert  + 2 \eta  \Vert g(\ell)   \Vert$, $\theta \triangleq \frac{1 + \sigma}{2}$, and
{\small
	\begin{align*}
	A_1(\eta) &\triangleq \frac{39[2\Vert y(0) - \one\bar{y}(0)\Vert + \frac{1}{L} \Vert s(0) - \one g(0)\Vert ]}{(\sigma\eta L)^{2/3}}\\
	A_2(\eta) &\triangleq  \frac{39 \sqrt{n} (\eta L)^{1/3}}{(\sigma )^{2/3}}.
	\end{align*} }
\end{lemma}


{\color{black}With the above preparations, we are ready to provide the proof of Theorem \ref{thm:str_cvx}.  Following similar techniques in \cite{nesterov2013introductory}, we recursively define a series of functions $\Phi_t(\omega)$.  We first define  $\Phi_0(\omega) = f(\bar{x}(0)) + \frac{\mu}{2} \Vert \omega - \bar{v}(0)\Vert^2$
and given $\Phi_t(\cdot)$, we define $\Phi_{t+1}(\cdot)$ by,
{\small\begin{align}		\label{eq:nes:inexact:phi_recursive}
	\Phi_{t+1}(\omega) 
	&= (1 - \alpha)\Phi_t(\omega) \nonumber\\
	&\quad + \alpha(\hat{f}(t) + \langle g(t), \omega - \bar{y}(t)\rangle + \frac{\mu}{2} \Vert \omega - \bar{y}(t)\Vert^2). 
	\end{align}}The following lemma gives a few properties on $\Phi_t(\cdot)$. Its proof is almost identical to the techniques used in the proof of Nesterov Gradient Descent in textbook \cite[Lemma 2.2.3]{nesterov2013introductory}.
\begin{lemma}\label{lem:nes:phi_def}
(a)	$\Phi_t(\cdot)$ can be written as, 
	\begin{equation}
	\Phi_t(\omega) = \phi_t^* + \frac{\mu}{2}\Vert \omega - \bar{v}(t)\Vert^2 \label{eq:nes:inexact:phi_quadratic}
	\end{equation}
	where $\phi_0^* = f(\bar{x}(0))$, and given $\phi_t^*$, $\phi_{t+1}^*$ is defined as,
	{\small 	\begin{align}
		\phi_{t+1}^* &= (1 - \alpha )\phi_t^* + \frac{\mu}{2}  (1 - \alpha ) \alpha \Vert \bar{v}(t) - \bar{y}(t)\Vert^2 + \alpha \hat{f}(t)\nonumber\\
		&\qquad  - \frac{1}{2}\eta  \Vert g(t)\Vert^2 + (1 - \alpha ) \alpha  \langle g(t), \bar{v}(t) - \bar{y}(t)\rangle. \label{eq:nes:inexact:phi}
		\end{align} }
(b)	For any $t$, we have 
	\begin{equation}
	\Phi_t(\omega) \leq f(\omega) + (1-\alpha)^t(\Phi_0(\omega) - f(\omega)) . \label{eq:nes:inexact:claim2}
	\end{equation}
\end{lemma}
\noindent\textit{Proof of Lemma~\ref{lem:nes:phi_def}:} We first show part (a). By (\ref{eq:nes:inexact:phi_recursive}), $\Phi_t$ is always a quadratic function. Since $ \nabla^2 \Phi_{0}(\omega) = \mu I$ and $\nabla^2 \Phi_{t+1}(\omega) = (1 - \alpha)\nabla^2 \Phi_t(\omega) + \alpha\mu I $,
we get $\nabla^2 \Phi_t(\omega) = \mu I$ for all $t$. We next show by induction that $\Phi_t(\cdot)$ achieves its minimum at $\bar{v}(t)$. Firstly, $\Phi_0$ achieves its minimum at $\bar{v}(0)$. Assume $\Phi_t(\cdot) $ achieves minimum at $\bar{v}(t)$. Since $\Phi_t$ is a quadratic function with Hessian $\mu I$, we have $\nabla \Phi_t(\bar{v}(t+1)) = \mu (\bar{v}(t+1) - \bar{v}(t))$. Then by (\ref{eq:nes:inexact:phi_recursive}), we have
{\small 	\begin{align*}
	&	\nabla \Phi_{t+1}(\bar{v}(t+1)) \\ 
	&= (1 - \alpha) \mu( \bar{v}(t+1)- \bar{v}(t) ) + \alpha( g(t) + \mu ( \bar{v}(t+1) - \bar{y}(t)) )\\
	&= \mu \big[\bar{v}(t+1) - (1-\alpha)\bar{v}(t) - \alpha \bar{y}(t) + \frac{\alpha}{\mu} g(t) \big] =0.
	\end{align*}}where the last equality follows from (\ref{eq:nes:update_ave_b}) and the fact $\frac{\eta}{\alpha} = \frac{\alpha}{\mu}$. 
Hence $\Phi_{t+1}$ achieves its optimum at $\bar{v}(t+1)$. Now we have shown that $\Phi_t(\cdot)$ is a quadratic function that achieves minimum at $\bar{v}(t)$ with Hessian $\mu I$. This implies (\ref{eq:nes:inexact:phi_quadratic}) is true. It remains to calculate $\phi_t^*$.
Clearly, $\phi_0^* = f(\bar{x}(0))$. Setting $\omega = \bar{y}(t)$ in (\ref{eq:nes:inexact:phi_recursive}), we get $\Phi_{t+1}(\bar{y}(t)) = (1 - \alpha)\Phi_t(\bar{y}(t)) + \alpha\hat{f}(t)$.
Combining this with (\ref{eq:nes:inexact:phi_quadratic}), we can get (\ref{eq:nes:inexact:phi}).

%
We next show part (b). Clearly (\ref{eq:nes:inexact:claim2}) is true for $t=0$. Assuming it's true for $t$, then for $t+1$, we have by (\ref{eq:nes:inexact:phi_recursive}) and (\ref{eq:nes:inexact:noisy_oracle_1}),
{\small\begin{align*}
	& \Phi_{t+1}(\omega) \\
	&= (1 - \alpha)\Phi_t(\omega) + \alpha(\hat{f}(t) + \langle g(t), \omega - \bar{y}(t)\rangle + \frac{\mu}{2} \Vert \omega - \bar{y}(t)\Vert^2)\\
	&\leq(1 - \alpha)\Phi_t(\omega) + \alpha f(\omega) \leq f(\omega) + (1-\alpha)^{t+1}(\Phi_0(\omega) - f(\omega)).
	\end{align*}}\qedd

The major step of proving Theorem~\ref{thm:str_cvx} is to show the following inequality,
\begin{equation}
f(\bar{x}(t)) \leq \phi_t^* + O((1-\alpha)^t). \label{eq:nes:inexact:phi_f_bound}
\end{equation} 
If (\ref{eq:nes:inexact:phi_f_bound}) is true, then combining (\ref{eq:nes:inexact:phi_f_bound}) with (\ref{eq:nes:inexact:claim2}), we have 
{\small\begin{align*}
	f(\bar{x}(t)) &\leq \phi_t^* + O((1-\alpha)^t)  \leq \Phi_t(x^*) +O((1-\alpha)^t) \\
	& \leq f^* +(1-\alpha)^t(\Phi_0(x^*) - f(x^*))+O((1-\alpha)^t),
	\end{align*}}which implies $f(\bar{x}(t)) - f^* = O((1-\alpha)^t)$, that is part (a) of the theorem. 
We put the proof of \eqref{eq:nes:inexact:phi_f_bound} in Section \ref{subsec:convergence}. We will also derive part (b) of the theorem (i.e. $ f( y_i(t)) -  f^*  = O((1- \sqrt{\mu\eta})^{t})$) in Section \ref{subsec:convergence}, which will be an easy corollary of \eqref{eq:nes:inexact:phi_f_bound}. 
\begin{remark}\label{rem:inexact_contribution_sc}
	In the context of inexact Nesterov method \cite{devolder2013first} for strongly convex and smooth functions, the final part of our proof essentially shows that for any inexact Nesterov method for sequence $\bar{x}(t)$, $\bar{y}(t)$, $\bar{v}(t)$ with inexact gradient $g(t)$ and step size $\eta$, when the error at time $t$ depends on past iterates of the form $A_1(\eta)\theta^t + A_2(\eta) \sum_{k=0}^t \theta^{t-k}(\Vert\bar{x}(k) - \bar{y}(t)\Vert + \eta \Vert g(k)\Vert )$ for some $\theta\in[0,1)$, where $A_2(\eta)$ is a positive constant that depends on $\eta$ and satisfies $\lim_{\eta\rightarrow 0} A_2(\eta)<\infty $, then the algorithm can have exact convergence if $\eta$ is small enough. This stands in contrast to the results in \cite[Sec. 5.3]{devolder2013first} where only approximate convergence can be obtained when the error is a constant throughout the iterations. 
\end{remark}
}

\subsection{Proof of the Bounded Consensus Error (Lemma \ref{lem:nes:rel_con_err})}\label{subsec:con_err} We now give the proof of Lemma \ref{lem:nes:rel_con_err}. \ifthenelse{\boolean{appendixin}}{We will frequently use the following two lemmas, whose proofs are deferred to Appendix-\ref{subsec:useful_fact}.}{We will frequently use the following two straightforward lemmas, whose proofs are omitted due to space limit. The proofs can be found in the full version of this paper \cite[Appendix-A]{fullversion}.}
\begin{lemma}\label{prop:simp_eq}
	The following inequality is true.
	\begin{equation}
	\Vert \bar{y}(t+1) - \bar{y}(t)\Vert  \leq \Vert \bar{y}(t) - \bar{x}(t) \Vert+2\eta \Vert g(t) \Vert \label{eq:nes:y_y}
	\end{equation}
\end{lemma}
\begin{lemma}\label{prob:nes:ineq_useful}
	The following inequalities are true. 
	\begin{eqnarray}
	\Vert \nabla(t+1) - \nabla(t)\Vert &\leq& L\Vert y(t+1) - y(t)\Vert\label{eq:nes:ineq_useful_a}\\
	\Vert g(t) - \nabla f(\bar{y}(t))\Vert &\leq &\frac{L}{\sqrt{n}} \Vert y(t) - \one \bar{y}(t)\Vert\label{eq:nes:ineq_useful_b}
	\end{eqnarray}
	
\end{lemma}


The proof of Lemma \ref{lem:nes:rel_con_err} is separated into two steps. In step 1, we treat the algorithm (\ref{eq:nes:update_vector}) as a linear system and derive a linear system inequality (\ref{eq:nes:rel_con:z_recursive}). In step 2, we analyze the state transition matrix in (\ref{eq:nes:rel_con:z_recursive}) and prove its spectral properties, from which we conclude the lemma. 
	
	\textbf{Step 1: A Linear System Inequality.} Define state $z(t) = [	\Vert v(t) - \one \bar{v}(t) \Vert,
	\Vert y(t) - \one \bar{y}(t)\Vert,
	\frac{1}{L}\Vert s(t) - \one g(t)\Vert ]^T\in\R^3$ and we will derive a linear system inequality that bounds the update of $z(t)$, given as follows
		\begin{align}
	z(t+1)\leq G(\eta) z(t) + b(t). \label{eq:nes:rel_con:z_recursive}
	\end{align}
	Here $b(t) = [0,0,\sqrt{n}a(t)]^T\in\R^3 $ is the input to the system with
	$$ a(t) \triangleq    \Vert  \bar{y}(t)-\bar{x}(t)\Vert  + 2 \eta  \Vert g(t)   \Vert.$$
The state transition matrix $G(\eta)\in\R^{3\times 3}$ is given by 
{\small	$$G(\eta) = \left[	\begin{array}{ccc}
	(1-\alpha)\sigma & \alpha \sigma&  \frac{\eta L}{\alpha } \\
	\frac{1-\alpha}{1+\alpha}\alpha\sigma  & \frac{1+\alpha^2}{1+\alpha} \sigma & 2\eta L  \\
	 \alpha\sigma  & 2  &\sigma+ 2\eta L
	\end{array}\right].$$  }
	
	We now prove (\ref{eq:nes:rel_con:z_recursive}). By (\ref{eq:nes:update_vector_a}) and (\ref{eq:nes:update_ave_a}), we have
	\begin{align}
		\Vert x(t+1) - \one \bar{x}(t+1)\Vert & = \Vert [Wy(t) - \one \bar{y}(t) ]- \eta [s(t) - \one g(t)] \Vert \nonumber \\
		&\leq \sigma \Vert y(t) - \one \bar{y}(t) \Vert + \eta \Vert s(t) - \one g(t)\Vert . \label{eq:nes:rel_con:ineq_x}
	\end{align}
	By (\ref{eq:nes:update_vector_b}) and (\ref{eq:nes:update_ave_b}), we have
	{\small
	\begin{align}
	&\Vert v(t+1) - \one \bar{v}(t+1) \Vert \nonumber \\
	& \leq \Vert (1-\alpha)[Wv(t) - \one\bar{v}(t) ] + \alpha [Wy(t)-\one\bar{y}(t)]  - \frac{\eta}{\alpha} [s(t) - \one g(t) ] \Vert\nonumber\\
		&\leq (1-\alpha)\sigma \Vert v(t) - \one \bar{v}(t)\Vert + \alpha\sigma\Vert y(t) - \one\bar{y}(t)\Vert +  \frac{\eta}{\alpha} \Vert s(t) - \one g(t)\Vert.\label{eq:nes:rel_con:ineq_v}
	\end{align}}
	
	By (\ref{eq:nes:update_vector_c}) and (\ref{eq:nes:update_ave_c}), we have
{\small	\begin{align}
		&\Vert y(t+1) - \one \bar{y}(t+1)\Vert\nonumber\\
		&\leq  \frac{1}{1+\alpha}\Vert x(t+1) - \one\bar{x}(t+1)\Vert  + \frac{\alpha}{1+\alpha} \Vert v(t+1) - \one \bar{v}(t+1)\Vert\nonumber\\
	&\leq \frac{1}{1+\alpha}\bigg[ \sigma \Vert y(t) - \one \bar{y}(t) \Vert + \eta \Vert s(t) - \one g(t)\Vert  \bigg]\nonumber\\
		&\qquad +  \frac{\alpha}{1+\alpha}  \bigg[ (1-\alpha)\sigma \Vert v(t) - \one \bar{v}(t)\Vert + \alpha\sigma\Vert y(t) - \one\bar{y}(t)\Vert \nonumber\\
		&\qquad  +  \frac{\eta}{\alpha} \Vert s(t) - \one g(t)\Vert \bigg] \nonumber\\
		&\leq \frac{1-\alpha}{1+\alpha} \alpha \sigma \Vert v(t) -\one \bar{v}(t)\Vert  + \frac{1+\alpha^2}{1+\alpha} \sigma \Vert y(t) - \one\bar{y}(t)\Vert \nonumber\\
	&	\qquad + 2\eta  \Vert s(t) - \one g(t)\Vert\label{eq:nes:rel_con:ineq_y}
	\end{align}}where we have used (\ref{eq:nes:rel_con:ineq_x}) and (\ref{eq:nes:rel_con:ineq_v}) in the second inequality\textcolor{black}{, and the fact that $\frac{1}{1+\alpha}<1$ in the third inequality. }

	By (\ref{eq:nes:update_vector_d}) and (\ref{eq:nes:update_ave_d}), we have
{\small	\begin{align}
	&	\Vert s(t+1) - \one g(t+1)\Vert \nonumber\\
		&= \Vert W s(t) - \one g(t)   + [\nabla(t+1) - \nabla(t) - \one (g(t+1) - g(t))]\Vert \nonumber \\
		&\stackrel{(a)} {\leq}\sigma \Vert s(t) - \one g(t)\Vert + \Vert \nabla (t+1) - \nabla(t)\Vert \nonumber\\
		&\stackrel{(b)}{\leq} \sigma \Vert s(t) - \one g(t)\Vert + L \Vert y (t+1) - y(t)\Vert  \label{eq:nes:rel_con:ineq_s_1}
	\end{align}}where in (a) we have used the fact that \textcolor{black}{by $g(t) = \frac{1}{n}\one^T\nabla(t)$, }
{\small	\begin{align*}
		&\big\Vert[\nabla(t+1)-\nabla(t)]  -   [\one g(t+1) - \one g(t)]\big\Vert^2 \\
		&= \Vert \nabla(t+1)-\nabla(t) \Vert^2 - n \Vert g(t+1) - g(t)\Vert ^2
	\end{align*}
}and in (b) we have used (\ref{eq:nes:ineq_useful_a}).
	
	Now we expand $y(t+1) - y(t)$.
{\small	\begin{align}
		&\Vert y(t+1) - y(t)\Vert \nonumber\\
		&\leq \Vert y(t+1) - \one \bar{y}(t+1)\Vert + \Vert y(t) - \one \bar{y}(t)\Vert + \Vert \one \bar{y}(t+1) - \one \bar{y}(t)\Vert\nonumber\\
		&\leq \frac{1-\alpha}{1+\alpha} \alpha \sigma \Vert v(t) -\one \bar{v}(t)\Vert  + \big[ \frac{1+\alpha^2}{1+\alpha} \sigma +1\big] \Vert y(t) - \one\bar{y}(t)\Vert  \nonumber\\
		&\qquad +  2\eta\Vert s(t) - \one g(t)\Vert + \Vert \one \bar{y}(t+1) - \one \bar{y}(t)\Vert \nonumber\\
		&\leq \alpha \sigma \Vert v(t) -\one \bar{v}(t)\Vert  +2 \Vert y(t) - \one\bar{y}(t)\Vert  \nonumber\\
		&\qquad + 2\eta \Vert s(t) - \one g(t)\Vert + \sqrt{n}   [\Vert  \bar{y}(t)-\bar{x}(t)\Vert + 2\eta \Vert g(t)   \Vert ]\nonumber
	\end{align}}where in the last inequality we have used (\ref{eq:nes:y_y}) (Lemma~\ref{prop:simp_eq}).
{\small	Combining the above with (\ref{eq:nes:rel_con:ineq_s_1}), we get
	\begin{align}
		&\frac{1}{L}\Vert s(t+1) - \one g(t+1)\Vert \nonumber\\
		&\leq \sigma \frac{1}{L}  \Vert s(t) - \one g(t)\Vert +  \Vert y (t+1) - y(t)\Vert  \nonumber\\
		&\leq  \alpha \sigma   \Vert v(t) -\one \bar{v}(t)\Vert+2\Vert y(t) - \one\bar{y}(t)\Vert \nonumber \\
		&\quad +(\sigma+ 2\eta L) \frac{1}{L}\Vert s(t) - \one g(t)\Vert   + \sqrt{n} \big[ \Vert  \bar{y}(t)-\bar{x}(t)\Vert  +2\eta   \Vert g(t)   \Vert\big]. \label{eq:nes:rel_con:ineq_s} 
	\end{align} }
	
	Combining (\ref{eq:nes:rel_con:ineq_v}) (\ref{eq:nes:rel_con:ineq_y})  (\ref{eq:nes:rel_con:ineq_s}) gives the linear system inequality (\ref{eq:nes:rel_con:z_recursive}). 

\textbf{Step 2: Spectral Properties of $G(\eta)$.} 
We give the following lemma regarding $G(\eta)$. We provide a proof-sketch here while the complete proof can be found in \ifthenelse{\boolean{appendixin}}{Appendix-\ref{subsec:G}.}{\cite[Appendix-B]{fullversion}.}
\begin{lemma}\label{lem:G_eig}
	When $0<\eta<  \min(\frac{1}{L} (\frac{1-\sigma}{8})^3, \frac{\sigma^3}{L64}   )$, the following holds.
\begin{enumerate}
	\item[(a)] We have $ \sigma +    (\sigma\eta L)^{1/3} < \rho(G(\eta))< \sigma+ 4 (\eta L)^{1/3} < \frac{1+\sigma}{2} = \theta.$
\end{enumerate}
\begin{enumerate}
 
	\item[(b)] The $(2,3)$th entry of $G(\eta)^t$ is upper bounded by $ [G(\eta)^t]_{2,3} \leq     \frac{39 (\eta L)^{1/3}}{(\sigma )^{2/3}} \rho(G(\eta))^t.$
	 \item[(c)] The entries in the 2nd row of $G(\eta)^t$ are all upper bounded by $\frac{39}{(\sigma\eta L)^{2/3}} \rho(G(\eta))^t.$
\end{enumerate}	
\end{lemma}

Notice in Lemma~\ref{lem:G_eig} (b), the constant $\frac{39 (\eta L)^{1/3}}{(\sigma )^{2/3}}$ converges to $0$ as $\eta \rightarrow 0$. This fact is crucial in the proof of eq. \eqref{eq:nes:inexact:phi_f_bound} in Sec.~\ref{subsec:convergence} (cf. (\ref{eq:nes:phi_f_geometric}) and the argument following it).

The step size condition in Lemma~\ref{lem:nes:rel_con_err} ensures the condition of Lemma~\ref{lem:G_eig}. 
By (\ref{eq:nes:rel_con:z_recursive}),
{\small$$z(t) \leq G(\eta)^t z(0) + \sum_{k=0}^{t-1} G(\eta)^{t-1-k} b(k).$$}Recall the second entry of $z(t)$ is $\Vert y(t) - \one \bar{y}(t)\Vert$, and $b(t) = [0,0, \sqrt{n}a(t)]^T$. Hence the above inequality implies,
{\small 
\begin{align*}
\Vert y(t) - \one \bar{y}(t)\Vert 
&\leq \max( [G(\eta)^t]_{2,1},[G(\eta)^t]_{2,2} , [G(\eta)^t]_{2,3}) \Vert z(0)\Vert_1  \\
&\quad +\sum_{k=0}^{t-1} [G(\eta)^{t-1-k}]_{2,3} \sqrt{n} a(k)\\
&\leq \underbrace{\frac{39[2\Vert y(0) - \one\bar{y}(0)\Vert + \frac{1}{L} \Vert s(0) - \one g(0)\Vert ]}{(\sigma\eta L)^{2/3}} }_{\triangleq A_1(\eta)}\theta^t \\
&\quad  + \underbrace{ \frac{39 \sqrt{n} (\eta L)^{1/3}}{(\sigma )^{2/3}}}_{\triangleq A_2(\eta)}\sum_{k=0}^{t-1} \theta^{t-1-k} a(k). 
\end{align*}}

\subsection{Proof of \eqref{eq:nes:inexact:phi_f_bound}}\label{subsec:convergence}

 By (\ref{eq:nes:inexact:phi}),
	{\footnotesize
	\begin{align}
	&	\phi_{t+1}^*- f(\bar{x}(t+1)) \nonumber\\
		&\stackrel{(a)}{=} (1-\alpha)(\phi_t^* - f(\bar{x}(t)))+ (1 - \alpha ) f(\bar{x}(t))+ \frac{\mu(1 - \alpha)  }{2 \alpha}   \Vert \bar{x}(t) - \bar{y}(t)\Vert^2 \nonumber\\
		&\qquad + \alpha \hat{f}(t) - \frac{1}{2}\eta  \Vert g(t)\Vert^2 + (1 - \alpha) \alpha  \langle g(t), \bar{v}(t) - \bar{y}(t)\rangle-f(\bar{x}(t+1)) \nonumber \\
		&\stackrel{(b)}{\geq} (1-\alpha)(\phi_t^* - f(\bar{x}(t))) + \frac{\mu (1-\alpha)}{2\alpha} \Vert \bar{x}(t) - \bar{y}(t)\Vert^2   \nonumber\\
	&	\qquad +(1 - \alpha )[\hat{f}(t) + \langle g(t),\bar{x}(t) - \bar{y}(t) \rangle  ]   + \alpha \hat{f}(t) - \frac{1}{2}\eta  \Vert g(t)\Vert^2 \nonumber\\
	&\qquad + (1 - \alpha) \alpha  \langle g(t), \bar{v}(t) - \bar{y}(t)\rangle-f(\bar{x}(t+1))\nonumber \\
		&\stackrel{(c)}{=} (1-\alpha)(\phi_t^* - f(\bar{x}(t)))+ \hat{f}(t)  - \frac{1}{2}\eta \Vert g(t)\Vert^2 \nonumber\\
		&\quad + \frac{\mu(1 - \alpha)  }{2 \alpha} \Vert \bar{x}(t) - \bar{y}(t)\Vert^2 -f(\bar{x}(t+1))\nonumber\\
		&\stackrel{(d)}{\geq}
		(1-\alpha)(\phi_t^* - f(\bar{x}(t))) + (\frac{1}{2}\eta - \eta^2L)\Vert g(t)\Vert^2 \nonumber\\
		&\quad +  \frac{\mu(1 - \alpha)  }{2 \alpha}  \Vert \bar{x}(t) - \bar{y}(t)\Vert^2   -\frac{L}{n}\Vert y(t) - \one\bar{y}(t)\Vert^2 \label{eq:nes:phi_f_onestep}
	\end{align}}where (a) is due to $\bar{v}(t) - \bar{y}(t) = \frac{1}{\alpha}(\bar{y}(t) - \bar{x}(t)) $ and (b) is due to (\ref{eq:nes:inexact:noisy_oracle_1}) (for $\omega = \bar{x}(t)$). In (c), we have used $(1-\alpha)\alpha (\bar{v}(t) - \bar{y}(t)) + (1-\alpha) (\bar{x}(t) - \bar{y}(t)) 
	= (1-\alpha)[\alpha\bar{v}(t) + \bar{x}(t) - (1+\alpha)\bar{y}(t)]=0$.
	In (d), we have used $f(\bar{x}(t+1)) \leq \hat{f}(t) + (\eta^2L-\eta)\Vert g(t)\Vert^2 + \frac{L}{n}\Vert y(t) - \one\bar{y}(t)\Vert^2$, which follows from (\ref{eq:nes:inexact:noisy_oracle_2}) (for $\omega=\bar{x}(t+1)$). We expand (\ref{eq:nes:phi_f_onestep}) recursively,  
	{\small	\begin{align}
			&\phi_{t+1}^*- f(\bar{x}(t+1)) \nonumber	\\
			 &\geq (1-\alpha)^{t+1}(\phi_0^* - f(\bar{x}(0))) + (\frac{1}{2}\eta - \eta^2L) \sum_{k=0}^t (1-\alpha)^{t-k} \Vert g(k)\Vert^2\nonumber \\
			&\quad   + \frac{\mu}{2\alpha } \sum_{k=0}^t(1-\alpha)^{t+1-k}\Vert \bar{x}(k) - \bar{y}(k)\Vert^2 \nonumber\\
			&\quad - \frac{L}{n}\sum_{k=0}^t(1-\alpha)^{t-k}\Vert y(k) - \one\bar{y}(k)\Vert^2.\label{eq:nes:phi_f_recursive}
		\end{align}}

%

Now we bound $\sum_{k=0}^t(1-\alpha)^{t-k}\Vert y(k) - \one\bar{y}(k)\Vert^2$. Fixing $t$, define vector $\nu,\pi_k\in \R^{t+1}$ (for $0\leq k\leq t$),
{\small\begin{align*}
\nu = &[A_1(\eta) (1-\alpha)^{\frac{t}{2}}, A_2(\eta)a(0)(1-\alpha)^{\frac{t-1}{2}},\\
&\quad  A_2(\eta)a(1)(1-\alpha)^{\frac{t-2}{2}}, \ldots, A_{2}(\eta)a(t-1)]^T \\
\pi_k  =& [\theta^k (1-\alpha)^{-\frac{t}{2}},\theta^{k-1}(1-\alpha)^{-\frac{t-1}{2}},\\
&\quad \ldots,\theta^0 (1-\alpha)^{-\frac{t-k}{2}},0,\ldots,0]^T .
\end{align*}}By Lemma \ref{lem:nes:rel_con_err} (the step size in Theorem~\ref{thm:str_cvx} implies the condition of Lemma \ref{lem:nes:rel_con_err}  holds), we have $\Vert y(k) - \one\bar{y}(k)\Vert\leq \nu^T\pi_k$, and hence $\Vert y(k) - \one\bar{y}(k)\Vert^2\leq \nu^T\pi_k\pi_k^T\nu$. Therefore,
{\small\begin{align}
\sum_{k=0}^t(1-\alpha)^{t-k}\Vert y(k) - \one\bar{y}(k)\Vert^2\leq \nu^T \Pi \nu \label{eq:nes:nupinu}
\end{align}}where $\Pi = \sum_{k=0}^t (1-\alpha)^{t-k}\pi_k\pi_k^T\in\R^{(t+1)\times(t+1)}$ is a symmetric matrix. Let $\Pi$'s $(p,q)$'th element be $\Pi_{pq}$. 
When $q\geq p$, the $(p,q)$th element of $\pi_k \pi_k^T$ is given by $[\pi_k \pi_k^T]_{p,q} = [\pi_k]_p [\pi_k]_q$, which equals $\theta^{2k+2-p-q}(1-\alpha)^{-t-1+\frac{p+q}{2}} $if $k\geq q-1$, and $0$ if $k< q-1$. 
Then, when $q\geq p$, 
{\small
\begin{align*}
\Pi_{pq} &= \sum_{k=q-1}^t (1-\alpha)^{t-k} \theta^{2k+2-p-q} (1-\alpha)^{-t-1 + \frac{p+q}{2}}\\
&= (\frac{\theta}{\sqrt{1-\alpha}})^{q-p} \frac{1- (\frac{\theta^2}{1-\alpha})^{t-q+2}}{1 - \frac{\theta^2}{1-\alpha}}.
\end{align*}}

By the step size in Theorem~\ref{thm:str_cvx} we have $\eta< \frac{(1-\sigma)^2}{4\mu}$, and hence $\alpha< \frac{1-\sigma}{2} = 1-\theta$. Therefore, $\frac{\theta}{\sqrt{1-\alpha}} < \sqrt{\theta} <1$, and $1 - \frac{\theta}{\sqrt{1-\alpha}} > 1- \sqrt{\theta}>  \frac{1-\theta}{2}= \frac{1-\sigma}{4}  $. Therefore,
{\small
\begin{align*}
\sum_{q=p+1}^{t+1} \Pi_{pq} &< \frac{1}{1 - \frac{\theta^2}{1-\alpha}} \sum_{q=p+1}^{t+1} (\frac{\theta}{\sqrt{1-\alpha}})^{q-p}\\
& = \frac{1}{1 - \frac{\theta^2}{1-\alpha}}\frac{\theta}{\sqrt{1-\alpha}} \frac{1 - (\frac{\theta}{\sqrt{1-\alpha}})^{t-p+1}}{1 - \frac{\theta}{\sqrt{1-\alpha}}}
< \frac{16}{(1-\sigma)^2}.
\end{align*}}
And similarly,
{\small
\begin{align*}
\sum_{q=1}^{p} \Pi_{pq} & = \sum_{q=1}^{p} \Pi_{qp} 
< \frac{1}{1 - \frac{\theta^2}{1-\alpha}} \sum_{q=1}^{p} (\frac{\theta}{\sqrt{1-\alpha}})^{p-q} \\
 & = \frac{1}{1 - \frac{\theta^2}{1-\alpha}}  \frac{1 - (\frac{\theta}{\sqrt{1-\alpha}})^{p}}{1 - \frac{\theta}{\sqrt{1-\alpha}}}
< \frac{16}{(1-\sigma)^2}.
\end{align*}}

Hence, by Gershgorin Disk Theorem\cite{gershgorin1931uber}, $\rho( \Pi)\leq \max_{p} (\sum_{q=1}^{t+1}\Pi_{pq}) \leq 32/(1-\sigma)^2. $
Combining the above with (\ref{eq:nes:nupinu}), 
{\small
\begin{align*}
&\sum_{k=0}^t(1-\alpha)^{t-k} \Vert y(k) - \one \bar{y}(k)\Vert^2 \leq \rho(\Pi) \Vert \nu\Vert^2 \\
&\leq   \frac{32}{(1-\sigma)^2} \Big[A_1(\eta)^2 (1-\alpha)^t + A_2(\eta)^2\sum_{k=0}^{t-1} (1-\alpha)^{t-k-1} a(k)^2\Big]\\
&\leq   \frac{32}{(1-\sigma)^2} \Big[A_1(\eta)^2 (1-\alpha)^t \\
&\quad+ 2A_2(\eta)^2\sum_{k=0}^{t-1} (1-\alpha)^{t-k-1} \Vert \bar{y}(k) -  \bar{x}(k) \Vert^2  \\
&\quad + 8\eta^2 A_2(\eta)^2\sum_{k=0}^{t-1} (1-\alpha)^{t-k-1} \Vert g(k) \Vert^2\Big]
\end{align*}}where in the last step, we have used by definition, $a(k)^2 =  (  \Vert  \bar{y}(k)-\bar{x}(k)\Vert  + 2 \eta  \Vert g(k)   \Vert)^2 \leq 2 \Vert  \bar{y}(k)-\bar{x}(k)\Vert^2 + 8\eta^2 \Vert g(k)\Vert^2$. Now returning to (\ref{eq:nes:phi_f_recursive}), we get,
		{\small	\begin{align}
	&\phi_{t+1}^*- f(\bar{x}(t+1))	\nonumber \\
	&\geq -\frac{L}{n} \frac{32 A_1(\eta)^2}{(1-\sigma)^2} (1-\alpha)^t \nonumber\\
&	\quad +    \overbrace{(\frac{1}{2}\eta - \eta^2L- \frac{256 L \eta^2 A_2(\eta)^2}{n(1-\sigma)^2(1-\alpha)}     )}^{\triangleq A_3(\eta) } \sum_{k=0}^t (1-\alpha)^{t-k} \Vert g(k)\Vert^2\nonumber \\
	&\quad + \frac{1}{\alpha} \underbrace{\big[\frac{\mu}{2}  - \frac{64 L A_2(\eta)^2 \alpha }{n (1-\sigma)^2 (1-\alpha)^2}   \big]}_{\triangleq A_4(\eta) }\sum_{k=0}^t(1-\alpha)^{t+1-k}\Vert  \bar{y}(k) - \bar{x}(k) \Vert^2. \label{eq:nes:phi_f_geometric}
	\end{align}}

	To prove (\ref{eq:nes:inexact:phi_f_bound}), it remains to check $A_3(\eta)$ and $A_4(\eta)$ are positive.
Plugging in $A_2(\eta) =  \frac{39 \sqrt{n} (\eta L)^{1/3}}{(\sigma )^{2/3}} $ into $A_3(\eta)$ and using $(1-\alpha) > \frac{1}{2}$ (equivalent to $ \eta < \frac{1}{4\mu}$ using $\alpha = \sqrt{\mu\eta}$, and implied by the step size condition in Theorem~\ref{thm:str_cvx}) , we have 
{\small	\begin{align*}
	A_3(\eta)   
	&\geq  \eta (\frac{1}{2} - \eta L - \frac{512\times 39^2 (\eta L)^{5/3} }{(1-\sigma)^2(\sigma )^{4/3}  } )> 0
	\end{align*}}where in the second inequality we have used by the step size condition in Theorem~\ref{thm:str_cvx}, $\eta L<\frac{1}{4}$, and $\frac{512\times 39^2 (\eta L)^{5/3} }{(1-\sigma)^2(\sigma )^{4/3}  }  < \frac{1}{4}$ ($\Leftarrow \eta < \frac{(1-\sigma)^{1.2} \sigma^{0.8}}{7872 L }$). For $A_4(\eta)$, similarly plugging in $A_2(\eta) $ and using $1-\alpha>\frac{1}{\sqrt{2}}$ ($\Leftarrow \eta< \frac{1}{25 \mu}$) and $\alpha =\sqrt{\mu \eta}$, we have
{\small	\begin{align*}
	 A_4(\eta)  \geq \frac{\mu}{2}  - \frac{128\times 39^2   }{ (1-\sigma)^2\sigma^{4/3}} \sqrt{L \mu} (\eta L)^{7/6}  > 0
	\end{align*} }where in the last inequality we have used by the step size condition in Theorem~\ref{thm:str_cvx}, $\frac{\mu}{2}  \geq \frac{128\times 39^2   }{ (1-\sigma)^2\sigma^{4/3}} \sqrt{L \mu} (\eta L)^{7/6}  $ ($\Leftarrow \eta< \frac{\sigma^{8/7}(1-\sigma)^{12/7} }{61909  L} (\frac{\mu}{L})^{3/7}$). So we have proven (\ref{eq:nes:inexact:phi_f_bound}).

{\color{black}
At last we will prove part (b) of Theorem \ref{thm:str_cvx}. Using (\ref{eq:nes:inexact:claim2}), $\phi_{t}^*\leq \Phi_t(x^*) \leq f^* + O((1-\alpha)^t)$. Hence $\phi_{t+1}^* - f(\bar{x}(t+1)) \leq \phi_{t+1}^* - f^* =O((1-\alpha)^t) $. Using (\ref{eq:nes:phi_f_geometric}), we have $ A_3(\eta)\Vert g(t)\Vert^2 + \frac{1-\alpha}{\alpha} A_4(\eta) \Vert \bar{y}(t) - \bar{x}(t)\Vert^2= O((1-\alpha)^t)$.
Therefore, both $\Vert g(t)\Vert^2$ and $\Vert \bar{y}(t) - \bar{x}(t)\Vert^2$ are $O((1-\alpha)^t)$. Then, the $a(t)$ defined in Lemma~\ref{lem:nes:rel_con_err} is $O((1-\alpha)^{t/2})$. By Lemma~\ref{lem:nes:rel_con_err}, we also have $\Vert y(t) - \one\bar{y}(t)\Vert = O((1-\alpha)^{t/2})$ (where we have used an easy-to-check fact: $\sqrt{1-\alpha}>\theta$). Since $f$ is $\mu$ strongly convex, we have $f(\bar{x}(t)) -f^* \geq \frac{\mu}{2} \Vert \bar{x}(t) - x^*\Vert^2$ which implies $\Vert \bar{x}(t) - x^*\Vert = O((1-\alpha)^{t/2})$. Since $\Vert \bar{x}(t) - x^*\Vert$, $\Vert \bar{y}(t) - \bar{x}(t)\Vert$ and $\Vert y(t) - \one\bar{y}(t)\Vert $ are all $O((1-\alpha)^{t/2})$, by triangle inequality we have $\Vert y(t) - \one x^*\Vert = O((1-\alpha)^{t/2})$. Hence, for each $i$, $\Vert y_i(t) -  x^*\Vert = O((1-\alpha)^{t/2})$, and therefore by $L$-smoothness of $f$ we have $f( y_i(t)) -  f^*  = O((1-\alpha)^{t})$. }

%
%
%
	\section{Convergence Analysis of Acc-DNGD-NSC}\label{sec:convergence_nsc}
In this section, we will prove Theorem~\ref{thm:nsc:vanishing} and Theorem~\ref{thm:nsc:fixed}. We will first provide a proof overview in Section~\ref{sec:nsc:proofoverview} and then defer the detailed proof to the rest of the section. \ifthenelse{\boolean{appendixin}}{}{Due to space limit, some of the technical calculations are omitted and can be found in the full version of this paper \cite{fullversion}.}

\subsection{Proof Overview} \label{sec:nsc:proofoverview}
Same as (\ref{eq:nes:stack_notation}), we introduce matrix notations $x(t)$, $v(t)$, $y(t)$, $s(t)$, $\nabla(t)\in\R^{n\times N}$. We also define $\bar{x}(t)$, $\bar{v}(t)$, $\bar{y}(t)$, $\bar{s}(t)$ and $g(t)$ analogously. 
Then our algorithm in (\ref{eq:nsc:update_elementwise}) can be written as
\begin{subeqnarray}\label{eq:nsc:update_vector}
	x(t+1) &= & W y(t) - \eta_t s(t)\slabel{eq:nsc:update_vector_a}\\
	v(t+1) &=&  W v(t)  - \frac{\eta_t}{\alpha_t} s(t) \slabel{eq:nsc:update_vector_b}\\
	y(t+1) &= &  (1-\alpha_{t+1})x(t+1) +  \alpha_{t+1} v(t+1)  \slabel{eq:nsc:update_vector_c}\\
	s(t+1) &= & W s(t) + \nabla(t+1) - \nabla (t).\slabel{eq:nsc:update_vector_d}
\end{subeqnarray}

\vspace{6pt}
\noindent
\textbf{Overview of the Proof.} We derive a series of lemmas (Lemma \ref{lem:nsc:ave}, \ref{lem:nsc:inexact_grad}, \ref{lem:nsc:rel_con_err} and \ref{lem:nsc:convergence}) that will work for both the vanishing and the fixed step size case. We firstly derive the update formula for the average sequences (Lemma \ref{lem:nsc:ave}). Then, we show that the update rule for the average sequences is in fact CNGD-NSC (\ref{eq:nsc:cngd}) with inexact gradients \cite{devolder2014first}, and the inexactness is characterized by ``consensus error'' $\Vert y(t) - \one\bar{y}(t)\Vert$ (Lemma \ref{lem:nsc:inexact_grad}). The consensus error is bounded in Lemma \ref{lem:nsc:rel_con_err}. Then, we apply the proof of CNGD (see e.g. \cite{nesterov2013introductory}) to the average sequences in spite of the consensus error, and derive an intermediate result in Lemma \ref{lem:nsc:convergence}. Lastly, we finish the proof of  Theorem \ref{thm:nsc:vanishing} and Theorem \ref{thm:nsc:fixed} in Section \ref{subsec:nsc:vanishing} and \ifthenelse{\boolean{appendixin}}{Appendix-\ref{subsec:nsc:fixed_full}}{\cite[Appendix-I]{fullversion}} respectively.

As shown above, the proof is similar to that of Acc-DNGD-SC in Section~\ref{sec:convergence_sc}. The main difference lies in how we bound the consensus error (Lemma~\ref{lem:nsc:rel_con_err}) and how we apply the CNGD proof in Section~\ref{subsec:nsc:vanishing}. In what follows, we will mainly focus on the different parts while putting details for the parts that are similar to Acc-DNGD-SC into the Appendix.

{\color{black} We first derive Lemma~\ref{lem:nsc:ave} that characterizes the update rule for the average sequences. }

\begin{lemma}\label{lem:nsc:ave} The following equalities hold.
{\small	\begin{subeqnarray}\label{eq:nsc:update_ave}
		\bar{x}(t+1) &= &  \bar{y}(t) - \eta_t g(t) \slabel{eq:nsc:update_ave_a}\\
		\bar{v}(t+1) &= & \bar{v}(t) - \frac{\eta_t}{\alpha_t} g(t) \slabel{eq:nsc:update_ave_b}\\
		\bar{y}(t+1) &= & (1-\alpha_{t+1}) \bar{x}(t+1) + \alpha_{t+1} \bar{v}(t+1)  \slabel{eq:nsc:update_ave_c}\\
		\bar{s}(t+1) &= & \bar{s}(t) + g(t+1) -  g(t) = g(t+1)\slabel{eq:nsc:update_ave_d}
	\end{subeqnarray}}
\end{lemma}
\noindent\textit{Proof: } We omit the proof since these equalities can be easily derived using the fact that $W$ is doubly stochastic and the fact that $\bar{s}(0)=g(0)$. \qedd 

From (\ref{eq:nsc:update_ave_a})-(\ref{eq:nsc:update_ave_c}) we see that the sequences $\bar{x}(t)$, $\bar{v}(t)$ and $\bar{y}(t)$ follow a update rule  similar to the CNGD-NSC in (\ref{eq:nsc:cngd}). The only difference is that the $g(t)$ in (\ref{eq:nsc:update_ave_a})-(\ref{eq:nsc:update_ave_c})  is not the exact gradient $\nabla f(\bar{y}(t))$ in CNGD-NSC. In the following Lemma, we show that $g(t)$ is an inexact gradient.
\begin{lemma}\label{lem:nsc:inexact_grad} $\forall t$, $g(t)$ is an inexact gradient of $f$ at $\bar{y}(t)$ with error $O(\Vert y(t) - \one\bar{y}(t)\Vert ^2)$ in the sense that, $\forall \omega\in\R^N$,\footnote{The exact gradient $\nabla f(\bar{y}(t))$, satisfies $f(\bar{y}(t)) + \langle \nabla f(\bar{y}(t)), \omega - \bar{y}(t)\rangle  \leq f(\omega) \leq f(\bar{y}(t)) + \langle \nabla f(\bar{y}(t)), \omega - \bar{y}(t)\rangle + \frac{L }{2}\Vert \omega - \bar{y}(t)\Vert^2$. This is why we call $g(t)$ inexact gradient.}
{\small\begin{align}
f(\omega) &\geq \hat{f} (t) + \langle g(t), \omega - \bar{y}(t)\rangle   \label{eq:nsc:inexact:noisy_oracle_1}\\
f(\omega) &\leq \hat{f}(t) + \langle g(t), \omega - \bar{y}(t)\rangle + L \Vert \omega - \bar{y}(t)\Vert^2 \nonumber\\
&\quad + L  \frac{1}{n} \Vert y(t) - \one\bar{y}(t)\Vert^2, \label{eq:nsc:inexact:noisy_oracle_2}
\end{align}}where $\hat{f}(t) = \frac{1}{n} \sum_{i=1}^n  [ f_i(y_i(t)) + \langle \nabla f_i(y_i(t)), \bar{y}(t)-y_i(t) \rangle ].$
\end{lemma}

\noindent\textit{Proof: } We omit the proof since it's almost identical as the proof of Lemma \ref{lem:inexact_grad}.  \qedd 

We then bound in Lemma~\ref{lem:nsc:rel_con_err} the consensus error $\Vert y(t) - \one\bar{y}(t)\Vert $. The proof of Lemma~\ref{lem:nsc:rel_con_err} is given in Section~\ref{subsec:nsc:con_err}.
\begin{lemma}\label{lem:nsc:rel_con_err}
	Suppose the step sizes satisfy 
	\begin{itemize}
		\item[(i)] $\eta_t\geq \eta_{t+1}>0$,
		\item[(ii)] $\eta_0< \min(  \frac{\sigma^2}{9^3 L}, \frac{(1-\sigma)^3}{6144 L})$,
		\item[(iii)] $\sup_{t\geq 0} \frac{\eta_t}{\eta_{t+1}} \leq \min( ( \frac{\sigma+3}{\sigma+2} \frac{3}{4})^{\sigma/28} , \frac{16}{15+\sigma})$.
	\end{itemize}  
Then, we have,
	\begin{align*}
	&\Vert y(t) - \one\bar{y}(t) \Vert \\
	&  \leq \kappa \sqrt{n} \chi_2(\eta_t) \bigg[L\Vert \bar{y}(t) - \bar{x}(t)\Vert + \frac{8}{1-\sigma} L \eta_t \Vert g(t)\Vert \bigg]
	\end{align*}
	where $\chi_2:\R\rightarrow \R$ is a function satisfying $0<\chi_2(\eta_t) \leq \frac{2}{L^{2/3}} \eta_t^{1/3}$, and $\kappa =  \frac{6}{(1-\sigma)}$.
\end{lemma}

We next provide the following intermediate result, which essentially uses the same construction and derivations in the standard CNGD proof in textbook \cite[Lemma 2.2.3]{nesterov2013introductory}. 
\begin{lemma}\label{lem:nsc:convergence}
	Define $\gamma_0 = \frac{\alpha_0^2}{\eta_0(1-\alpha_0)} = \frac{L}{1-\alpha_0}$ (\textcolor{black}{the second equality is due to the definition $\alpha_0=\sqrt{\eta_0 L}$}). We define a series of functions $\Phi_t:\R^N\rightarrow \R$ where $t\geq 0$,  with $\Phi_0(\omega) = f(\bar{x}(0)) + \frac{\gamma_0}{2} \Vert \omega - \bar{v}(0)\Vert^2$ and 
	\begin{equation}
	\Phi_{t+1}(\omega)= (1 - \alpha_t)\Phi_t(\omega) + \alpha_t[\hat{f}(t) + \langle g(t), \omega - \bar{y}(t)\rangle]. \label{eq:nsc:inexact:phi_recursive}
	\end{equation}
	Then, we have,
		\begin{equation}
		\Phi_t(\omega) \leq f(\omega) + \lambda_t (\Phi_0(\omega) - f(\omega)) \label{eq:nsc:phi_upperbound}
		\end{equation}
		where $\lambda_t$ is defined through $\lambda_0=1$, and $\lambda_{t+1}=(1-\alpha_t)\lambda_t$. Further, we have function $\Phi_t(\omega)$ can be written as
			\begin{equation}
		\Phi_t(\omega) = \phi_t^* + \frac{\gamma_t}{2}\Vert \omega - \bar{v}(t)\Vert^2 \label{eq:nsc:inexact:phi_quadratic}
		\end{equation}
 where $\gamma_t $ is defined through $\gamma_{t+1} = \gamma_t (1-\alpha_t)$, and $\phi_t^*$ is some real number that satisfies $\phi_0^* = f(\bar{x}(0))$, and 
		\begin{equation}
		\phi_{t+1}^* = (1 - \alpha_t )\phi_t^*   + \alpha_t \hat{f}(t) - \frac{1}{2}\eta_t \Vert g(t)\Vert^2   +  \alpha_t \langle g(t), \bar{v}(t) - \bar{y}(t)\rangle.   \label{eq:nsc:inexact:phi}
		\end{equation}
\end{lemma}
\noindent \textit{Proof: } The proof is almost the same as that of Lemma~\ref{lem:nes:phi_def} in the proof of Theorem~\ref{thm:str_cvx}. For completeness we include a proof in \ifthenelse{\boolean{appendixin}}{Appendix-\ref{subsec:nsc:convergence}.}{\cite[Appendix-C]{fullversion}.}  \qedd 

{\color{black} Lemma~\ref{lem:nsc:convergence} has laid the ground work for proving the convergence of the inexact Nesterov Gradient descent \eqref{eq:nsc:update_ave} for sequence $\bar{x}(t), \bar{y}(t), \bar{v}(t)$. What remains to be done is to follow the proof strategies of Nesterov Gradient descent, and in the meanwhile carefully control the error caused by the inexactness in Lemma~\ref{lem:nsc:rel_con_err}. With this guideline, we will finish the proof of Theorem \ref{thm:nsc:vanishing} (vanishing step size case) in Section \ref{subsec:nsc:vanishing} and Theorem \ref{thm:nsc:fixed} (fixed step size case) in \ifthenelse{\boolean{appendixin}}{Appendix-\ref{subsec:nsc:fixed_full}.}{\cite[Appendix-I]{fullversion}.} }

{\color{black}
\begin{remark}\label{rem:inexact_contribution_nsc}
	In the context of inexact Nesterov method \cite{devolder2014first} for convex and smooth functions, the final part of our proof for Theorem~\ref{thm:nsc:vanishing} essentially shows that for any inexact Nesterov method for sequence $\bar{x}(t)$, $\bar{y}(t)$, $\bar{v}(t)$ with inexact gradient $g(t)$ and step size $\eta_t$, when the error at time $t$ depends on past iterates of the form $ A(\eta_t) (\Vert\bar{x}(t) - \bar{y}(t)\Vert +\eta_t \Vert g(t)\Vert )$, where $A(\eta_t)$ is a constant that depends on $\eta_t$ and satisfies $\lim_{\eta\rightarrow 0}\frac{A_2(\eta)}{\eta^{2/3}} <\infty$, then the algorithm can have exact convergence if $\eta_t = \frac{\eta}{(t+1)^\beta} $ for small enough $\eta$ and $\beta\in(0.6, 2)$. This stands in contrast to the results in \cite[Sec. 6.2]{devolder2014first} where inexact Nesterov methods diverge when the error is fixed throughout the iterations. 
\end{remark}
}
\subsection{Proof of the Bounded Consensus Error (Lemma \ref{lem:nsc:rel_con_err})}\label{subsec:nsc:con_err}
We will frequently use the following straightforward lemma, \ifthenelse{\boolean{appendixin}}{whose proof can be found in Appendix-\ref{subsec:nsc:useful_eq}.}{whose proof is omitted and can be found in \cite[Appendix-D]{fullversion}.} We will also use Lemma~\ref{prob:nes:ineq_useful} in Section \ref{subsec:con_err}, which still holds under the setting of this section.
\begin{lemma}\label{prop:nsc:simp_eq}
	The following equalities are true.
{\small	\begin{align}
&\bar{y}(t+1) - \bar{y}(t) \nonumber \\
&=\alpha_{t+1}(\bar{v}(t) - \bar{y}(t)) - \eta_t \big[\frac{\alpha_{t+1}}{\alpha_t} + 1 - \alpha_{t+1} \big] g(t),  \label{eq:nsc:y_y}\\
&\bar{v}(t+1) - \bar{y}(t+1) \nonumber \\
 &= (1-\alpha_{t+1}) (\bar{v}(t) - \bar{y}(t))  + \eta_t (1 - \alpha_{t+1})(1 - \frac{1}{\alpha_t}) g(t) \label{eq:nsc:v_y}
\end{align}}
\end{lemma}

%
\noindent\textit{Proof of Lemma \ref{lem:nsc:rel_con_err}: }

	\textbf{Overview of the proof.} The proof is divided into three steps. In step 1, we treat the algorithm (\ref{eq:nsc:update_vector}) as a linear system and derive a linear system inequality (\ref{eq:nsc:rel_con:z_recursive}). In step 2, we analyze the state transition matrix in (\ref{eq:nsc:rel_con:z_recursive}) and prove a few spectral properties. In step 3, we further analyze the linear system (\ref{eq:nsc:rel_con:z_recursive}) and bound the state by the input, from which the conclusion of the lemma follows. Throughout the proof, we will frequently use an easy-to-check fact: $\alpha_t$ is a decreasing sequence.
		
\textbf{Step 1: A Linear System Inequality.} 	
	Define $z(t) = [\alpha_{t}\Vert v(t) - \one \bar{v}(t) \Vert,
	\Vert y(t) - \one \bar{y}(t)\Vert,
	\Vert s(t) - \one g(t)\Vert ]^T\in\R^3$, $b(t) =[0,0,\sqrt{n}a(t)]^T\in\R^3$ where 
{\small	$$ a(t) \triangleq \alpha_t L\Vert  \bar{v}(t)-\bar{y}(t)\Vert  + 2\lambda L \eta_t  \Vert g(t)   \Vert $$}in which $\lambda \triangleq \frac{4}{1-\sigma}>1$. Then, we have the following linear system inequality holds.
%
{\small	\begin{align}
	z(t+1)& \leq 
	\overbrace{\left[	\begin{array}{ccc}
		\sigma &  0 & \eta_t  \\
		\sigma  &   \sigma &  2\eta_t  \\
		L &   2L & \sigma + 2\eta_t L 
		\end{array}\right]}^{\triangleq G(\eta_t) } z(t) + b(t)\label{eq:nsc:rel_con:z_recursive}
	\end{align}}The derivation of (\ref{eq:nsc:rel_con:z_recursive}) is almost the same as that of (\ref{eq:nes:rel_con:z_recursive}) (in the proof of Lemma~\ref{lem:nes:rel_con_err}) and is deferred to \ifthenelse{\boolean{appendixin}}{Appendix-\ref{subsec:nsc:linear_sys_ineq}.}{\cite[Appendix-E]{fullversion}}
	
	\textbf{Step 2: Spectral Properties of $G(\cdot)$.} 
When $\eta$ is positive, $G(\eta)$ is a nonnegative matrix and $G(\eta)^2$ is a positive matrix. By Perron-Frobenius Theorem (\cite[Theorem 8.5.1]{horn2012matrix}) $G(\eta)$ has a unique largest (in magnitude) eigenvalue that is a positive real with multiplicity $1$, and the eigenvalue is associated with an eigenvector with positive entries. We let the unique largest eigenvalue be $\theta(\eta) = \rho(G(\eta))$ and let its eigenvector be $\chi(\eta)=[\chi_1(\eta),\chi_2(\eta),\chi_3(\eta)]^T$, normalized by $\chi_3(\eta) = 1$. We give bounds on the eigenvalue and the eigenvector in the following lemmas. We defer the proof to \ifthenelse{\boolean{appendixin}}{Appendix-\ref{subsec:nsc:G}.}{\cite[Appendix-F]{fullversion}} 
	\begin{lemma} \label{lem:nsc:G_eig} When $0<\eta L<1$,
		we have $\sigma< \theta(\eta)< \sigma+4(\eta L)^{1/3}$, and $\chi_2(\eta) \leq  \frac{2}{L^{2/3}} \eta^{1/3}$.
	\end{lemma}
	\begin{lemma}\label{lem:nsc:G_eig_lowerboud}
	When $\eta\in (0, \frac{\sqrt{\sigma}}{ L 2\sqrt{2}} )$, $\theta(\eta) \geq \sigma+  (\sigma \eta L)^{1/3}$ and $\chi_1(\eta) < \frac{\eta}{(\sigma\eta L)^{1/3}}$.
\end{lemma}
\begin{lemma} \label{lem:nsc:G_chi}
	When $0< \zeta_2< \zeta_1 <\frac{\sigma^2}{9^3 L}$, then $\frac{\chi_1(\zeta_1)}{\chi_1(\zeta_2)}\leq  (\frac{\zeta_1}{\zeta_2})^{6/\sigma} $ and $\frac{\chi_2(\zeta_1)}{\chi_2(\zeta_2)}\leq  (\frac{\zeta_1}{\zeta_2})^{28/\sigma}$.
\end{lemma}

It is easy to check that, under our step size condition (ii) in Lemma~\ref{lem:nsc:rel_con_err}, all the conditions of Lemma \ref{lem:nsc:G_eig}, \ref{lem:nsc:G_eig_lowerboud}, \ref{lem:nsc:G_chi} are satisfied.

	\textbf{Step 3: Bound the state by the input.} With the above preparations, now we prove, by induction, the following statement, 
	\begin{equation}\label{eq:nsc:rel_con:induction_assumption}
		z(t) \leq     \sqrt{n} a(t)  \kappa \chi(\eta_t) 
	\end{equation}
where $\kappa = \frac{6}{1-\sigma}$. Eq. (\ref{eq:nsc:rel_con:induction_assumption}) is true for $t=0$, since the left hand side is zero when $t=0$. 

	Assume (\ref{eq:nsc:rel_con:induction_assumption}) holds for $t$. We now show (\ref{eq:nsc:rel_con:induction_assumption}) is true for $t+1$. We divide the rest of the proof into two sub-steps. Briefly speaking, step 3.1 proves that the input to the system (\ref{eq:nsc:rel_con:z_recursive}), $a(t+1)$ does not decrease too much compared to $a(t)$ ($a(t+1) \geq \frac{\sigma+3}{4} a(t)$); while step 3.2 shows that the state $z(t+1)$, compared to $z(t)$, decreases enough for (\ref{eq:nsc:rel_con:induction_assumption}) to hold for $t+1$.
	
	\textbf{Step 3.1: We prove that $a(t+1) \geq \frac{\sigma+3}{4} a(t)$. }  
	By (\ref{eq:nsc:v_y}),
{\small	\begin{align*}
	a(t+1)
	&=  \alpha_{t+1}L \Vert  \bar{v}(t+1)-\bar{y}(t+1)\Vert  + 2\lambda \eta_{t+1} L  \Vert g(t+1)  \Vert\\
	&=   \alpha_{t+1} L \Vert  (1-\alpha_{t+1}) (\bar{v}(t) - \bar{y}(t)) \\
	&\qquad  + (1-\alpha_{t+1})(1-\frac{1}{\alpha_{t}})\eta_{t} g(t)\Vert  +  2\lambda \eta_{t+1} L \Vert g(t+1)  \Vert\\
	&\geq \alpha_{t+1}(1-\alpha_{t+1}) L \Vert \bar{v}(t) - \bar{y}(t)\Vert \\
	&\qquad   - \frac{\alpha_{t+1}}{\alpha_{t}}(1-\alpha_{t+1})(1-\alpha_{t})\eta_{t} L \Vert g(t) \Vert \\
	&\qquad + 2\lambda \eta_{t+1} L\Vert g(t)  \Vert - 2\lambda \eta_{t+1} L \Vert g(t+1) - g(t)  \Vert.
	\end{align*}}Therefore, recalling $ a(t) = \alpha_t L\Vert  \bar{v}(t)-\bar{y}(t)\Vert  + 2\lambda L \eta_t  \Vert g(t)   \Vert $, we have
{\small	\begin{align}
&	a(t)-a(t+1) \nonumber  \\
	&\leq \Big[\alpha_{t}-\alpha_{t+1}(1-\alpha_{t+1})  \Big]  L \Vert \bar{v}(t) - \bar{y}(t)\Vert  \nonumber\\
	& \quad  +   \Big[\frac{\alpha_{t+1}}{\alpha_{t}}(1-\alpha_{t+1})(1-\alpha_{t}) \eta_{t} L   + 2\lambda\eta_{t}L - 2\lambda \eta_{t+1} L \Big] \Vert g(t) \Vert \nonumber\\
	&\quad  + 2\lambda \eta_{t+1} L \Vert g(t+1) - g(t)  \Vert  \nonumber\\
	&\leq \Big[\alpha_{t}-\alpha_{t+1}(1-\alpha_{t+1})  \Big]  L \Vert \bar{v}(t) - \bar{y}(t)\Vert  \nonumber\\
	&\quad  +   ( \eta_{t} + 2\lambda (\eta_{t}-\eta_{t+1})) L \Vert g(t) \Vert  + 2\lambda \eta_{t+1} L \Vert g(t+1) - g(t)  \Vert  \nonumber\\
	&\leq \max(1-\frac{\alpha_{t+1}}{\alpha_{t}}+\frac{\alpha_{t+1}^2}{\alpha_{t}}, \frac{1}{2\lambda }+ \frac{\eta_{t} - \eta_{t+1}}{\eta_{t}}) a(t) \nonumber\\
	&\quad  +2\lambda\eta_{t+1} L \Vert g(t+1) - g(t)  \Vert \label{eq:nsc:rel_con:a_a}
	\end{align}}where in the last inequality, we have used the elementary fact that for \textcolor{black}{four positive real numbers} $a_1,a_2,a_3,a_4$ and $x,y\geq 0$, we have $a_1 x+a_2 y = \frac{a_1}{a_3} a_3 x+ \frac{a_2}{a_4} a_4 y \leq \max(\frac{a_1}{a_3}, \frac{a_2}{a_4})(a_3x +a_4y)$.

	Next, we expand $\Vert g(t+1) - g(t)\Vert$,
{\small	\begin{align}
	&	\Vert g(t+1) - g(t)  \Vert  \nonumber\\
		& \leq \Vert g(t+1) - \nabla f(\bar{y}(t+1))  \Vert + \Vert g(t) - \nabla f(\bar{y}(t))\Vert  \nonumber\\
		&\quad  + \Vert \nabla f(\bar{y}(t+1))-\nabla f(\bar{y}(t)) \Vert\nonumber\\
		&\stackrel{(a)}{\leq} \frac{L}{\sqrt{n}} \Vert y(t+1) - \one \bar{y}(t+1)\Vert + \frac{L}{\sqrt{n}} \Vert y(t) - \one \bar{y}(t)\Vert \nonumber\\
		&\quad  + L\Vert \bar{y}(t+1) -\bar{y}(t)\Vert \nonumber\\
		&\stackrel{(b)}{\leq}  \frac{L}{\sqrt{n}} \sigma \alpha_{t} \Vert v(t) -\one \bar{v}(t)\Vert  +  \frac{L}{\sqrt{n}} 2 \Vert y(t) - \one\bar{y}(t)\Vert  \nonumber\\
		&\quad  +\frac{L}{\sqrt{n}}2\eta_{t} \Vert s(t) - \one g(t)\Vert + a(t)\nonumber\\
		&\stackrel{(c)}{\leq} L \sigma \kappa \chi_1(\eta_{t}) a(t)  +  2L  \kappa \chi_2(\eta_{t}) a(t)  +  2 L\eta_{t} \kappa \chi_3(\eta_{t}) a(t)  + a(t)\nonumber\\
		&\stackrel{(d)}{\leq} a(t)\Big\{ L \sigma \kappa \frac{\eta_{t}}{(\sigma \eta_{t} L)^{1/3}} +  2L   \kappa \frac{2 \eta_{t}^{1/3}}{L^{2/3} } +2 L\eta_{t}\kappa + 1 \Big\} \stackrel{(e)}{\leq}  8\kappa  a(t). \label{eq:nsc:rel_con:g_difference}
		\end{align}}Here (a) is due to (\ref{eq:nes:ineq_useful_b}); (b) is due to the second row of (\ref{eq:nsc:rel_con:z_recursive}) and the fact that $a(t)\geq L\Vert \bar{y}(t+1) - \bar{y}(t)\Vert$ (cf. (\ref{eq:nsc:y_y}));
	(c) is due to the induction assumption (\ref{eq:nsc:rel_con:induction_assumption}). In (d), we have used the bound on $\chi_1(\cdot)$ (Lemma \ref{lem:nsc:G_eig_lowerboud}), $\chi_2(\cdot)$ (Lemma \ref{lem:nsc:G_eig}), and $\chi_3(\eta_t) = 1$. In (e), we have used $\eta_tL<1$, $\sigma<1$ and $\kappa>1$.
	
	Combining (\ref{eq:nsc:rel_con:g_difference}) with (\ref{eq:nsc:rel_con:a_a}) and recalling $\kappa = \frac{6}{1-\sigma}, \lambda = \frac{4}{1-\sigma}$, we have
{\footnotesize	\begin{align*}
	&a(t)-a(t+1) \\
	  &\leq \max(1-\frac{\alpha_{t+1}}{\alpha_{t}}+\frac{\alpha_{t+1}^2}{\alpha_{t}}, \frac{1}{2\lambda }+ \frac{\eta_{t} - \eta_{t+1}}{\eta_{t}}) a(t)  +16\kappa \lambda\eta_{t+1} L a(t) \\
	&\leq \bigg[\max(1 - \frac{\eta_{t+1}}{\eta_{t}} +2\alpha_{t+1} , \frac{1-\sigma}{8}+  \frac{\eta_{t} - \eta_{t+1}}{\eta_{t}})   + \frac{384\eta_0 L}{(1-\sigma)^2}   \bigg] a(t) 
	\end{align*}}where in the last inequality, we have used the fact that
{\small	\begin{align*}
	&1-\frac{\alpha_{t+1}}{\alpha_{t}}+\frac{\alpha_{t+1}^2}{\alpha_{t} }< 1 - \frac{\alpha_{t+1}^2}{\alpha_{t}^2} + \alpha_{t+1} \\
	& = 1 - \frac{\eta_{t+1}}{\eta_{t}}(1-\alpha_{t+1}) + \alpha_{t+1}< 1 - \frac{\eta_{t+1}}{\eta_{t}} + 2\alpha_{t+1}
	\end{align*} }where the equality follows from the update rule for $\alpha_t$ (\ref{eq:nsc:alpha_update}). By the step size condition (iii) in Lemma~\ref{lem:nsc:rel_con_err}, $\frac{\eta_{t}}{\eta_{t+1}} \leq \frac{16}{15 + \sigma} $, and hence $ 1 - \frac{\eta_{t+1}}{\eta_{t}} \leq \frac{1-\sigma}{16}$. By the step size condition (ii), $2\alpha_{t+1}\leq 2\alpha_0 = 2\sqrt{\eta_0 L} \leq \frac{1-\sigma}{16} $, and $ \eta_0 L \frac{384}{(1-\sigma)^2} < \frac{1-\sigma}{16}$. Combining the above, we have $a(t)-a(t+1) \leq \frac{1-\sigma}{4} a(t)$. 
	Hence $a(t+1) \geq \frac{3+\sigma}{4} a(t)$.
	
	\textbf{Step 3.2: Finishing the induction. }  We have,
{\small	\begin{align}
	z(t+1) 
	&\stackrel{(a)}{\leq} G(\eta_t) z(t) + b(t) \nonumber\\
	&\stackrel{(b)}{\leq} G(\eta_t)  \sqrt{n} a(t)  \kappa\chi(\eta_t) + \sqrt{n}a(t)  \chi(\eta_t) \nonumber \\
	&\stackrel{(c)}{=}  \theta(\eta_t)   \sqrt{n} a(t) \kappa\chi(\eta_t) + \sqrt{n}a(t) \chi(\eta_t) \nonumber \\
	&= \sqrt{n} a(t)\chi(\eta_t) (\kappa \theta(\eta_t)+ 1 ) \nonumber \\
	&\stackrel{(d)}{\leq}  \sqrt{n} a(t+1)\chi(\eta_{t+1}) (\kappa \frac{\sigma+1}{2} +1)\frac{4}{3+\sigma}  \nonumber\\
	&\qquad \times \max( \frac{\chi_1(\eta_t)}{\chi_1(\eta_{t+1})},\frac{\chi_2(\eta_t)}{\chi_2(\eta_{t+1})} ,1)\nonumber\\
	&\stackrel{(e)}{=}    \sqrt{n} a(t+1)\chi(\eta_{t+1}) \frac{\sigma+2}{3}\kappa \frac{4}{\sigma+3} \nonumber\\
	&\qquad \times\max( \frac{\chi_1(\eta_t)}{\chi_1(\eta_{t+1})},\frac{\chi_2(\eta_t)}{\chi_2(\eta_{t+1})} ,1) \nonumber\\
	&\stackrel{(f)}{\leq }  \sqrt{n} a(t+1)\kappa \chi(\eta_{t+1}) 
	\end{align}}where (a) is due to (\ref{eq:nsc:rel_con:z_recursive}), and (b) is due to induction assumption (\ref{eq:nsc:rel_con:induction_assumption}), and (c) is because $\theta(\eta_t)$ is an eigenvalue of $G(\eta_t)$ with eigenvector $\chi(\eta_t)$, and (d) is due to step 3.1, and $\theta(\eta_t)< \sigma + 4 (\eta_0 L)^{1/3}< \frac{1+\sigma}{2}$ (by Lemma~\ref{lem:nsc:G_eig} and step size condition (ii) in Lemma~\ref{lem:nsc:rel_con_err}), and in (e), we have used the fact $\kappa \frac{\sigma+1}{2} + 1= \frac{\sigma+2}{3}\kappa$ (since $\kappa = \frac{6}{1-\sigma}$).
For (f),  we have used that by Lemma~\ref{lem:nsc:G_chi} and step size condition (iii) (in Lemma~\ref{lem:nsc:rel_con_err}),
	{\small	\begin{align*}
		\max( \frac{\chi_1(\eta_t)}{\chi_1(\eta_{t+1})},\frac{\chi_2(\eta_t)}{\chi_2(\eta_{t+1})},1 ) \leq (\frac{\eta_t}{\eta_{t+1}})^{28/\sigma}\leq \frac{\sigma+3}{\sigma+2} \frac{3}{4}.
		\end{align*}}
	
	
	Now, (\ref{eq:nsc:rel_con:induction_assumption}) is proven for $t+1$, and hence is true for all $t$. 
	Therefore, we have
	$$\Vert y(t) - \one \bar{y}(t)\Vert  \leq \kappa  \sqrt{n} \chi_2(\eta_t) a(t). $$
	Notice that $a(t) = \alpha_t L\Vert \bar{v}(t) - \bar{y}(t)\Vert + 2\lambda L \eta_t\Vert g(t)\Vert \leq L\Vert  \bar{y}(t) - \bar{x}(t)\Vert + \frac{8}{1-\sigma} L \eta_t  \Vert g(t)\Vert$ (by $\alpha_t(\bar{v}(t) - \bar{y}(t)) = (1-\alpha_t)(\bar{y}(t) - \bar{x}(t))$). The statement of the lemma follows.
\qedd

\subsection{Proof of Theorem \ref{thm:nsc:vanishing}} \label{subsec:nsc:vanishing}
We first introduce Lemma \ref{lem:nsc:alpha} regarding the asymptotic behavior of $\alpha_t$ and $\lambda_t$, the proof of which can be found in \ifthenelse{\boolean{appendixin}}{Appendix-\ref{subsec:nsc:alpha}.}{\cite[Appendix-G]{fullversion}.}

\begin{lemma}\label{lem:nsc:alpha}
	When the vanishing step size is used ($\eta_t = \frac{\eta}{ (t+t_0)^\beta}$, $t_0\geq 1$, $\beta\in(0,2)$), and $\eta_0 < \frac{1}{4L}$ (equivalently $\alpha_0<\frac{1}{2}$), we have
	\begin{itemize}
		\item [(i)] $\alpha_t\leq \frac{2}{t+1}$.
		\item [(ii)] $\lambda_t = O(\frac{1}{t^{2-\beta}}) $.
		\item [(iii)] $\lambda_t \geq \frac{D(\beta,t_0)}{(t+t_0)^{2-\beta}  } $ where $D(\beta,t_0)$ is some constant that only depends on $\beta$ and $t_0$, given by $D(\beta,t_0) = \frac{1}{(t_0+3)^2 e^{16+ \frac{6}{2-\beta} }}$.
	\end{itemize}
\end{lemma}


To prove Theorem \ref{thm:nsc:vanishing}, we first note that with the step size condition in Theorem~\ref{thm:nsc:vanishing}, all the conditions of Lemma~\ref{lem:nsc:rel_con_err} and \ref{lem:nsc:alpha} are satisfied, hence the conclusions of Lemma~\ref{lem:nsc:rel_con_err} and \ref{lem:nsc:alpha} hold. The major step of proving Theorem~\ref{thm:nsc:vanishing} is to show the following inequality,
	\begin{equation}
	\lambda_t(\Phi_0(x^*) - f^*) +\phi_t^*\geq f(\bar{x}(t) ).\label{eq:nsc:inexact:f_y_phi}
	\end{equation}
If (\ref{eq:nsc:inexact:f_y_phi}) is true, by (\ref{eq:nsc:inexact:f_y_phi}) and (\ref{eq:nsc:phi_upperbound}), we have 
{\small	\begin{align*}
	f(\bar{x}(t))&\leq \phi_t^* + \lambda_t(\Phi_0(x^*) - f^*) \leq \Phi_t(x^*)+ \lambda_t(\Phi_0(x^*) - f^*)\\
	& \leq f^* +2 \lambda_t(\Phi_0(x^*) - f^*).
	\end{align*}}Hence $f(\bar{x}(t)) - f^* = O(\lambda_t) = O(\frac{1}{t^{2-\beta}})$, i.e. the desired result of part (a) of Theorem ~\ref{thm:nsc:vanishing}. In what follows, we will prove (\ref{eq:nsc:inexact:f_y_phi}), after which we will also prove part (b) (bounding objective error $f(y_i(t)) - f^*$ for each individual agent), which will then be an easy corollary.  

	 We use induction to prove (\ref{eq:nsc:inexact:f_y_phi}). Firstly, (\ref{eq:nsc:inexact:f_y_phi}) is true for $t=0$, since $\phi_0^* = f(\bar{x}(0)) $ and $\Phi_0(x^*) \geq f(\bar{x}(0))\geq f^*$. Suppose it's true for $0,1,2,\ldots, t$. For $0\leq k\leq t$, by (\ref{eq:nsc:phi_upperbound}),
	$\Phi_k(x^*)\leq f^* + \lambda_k(\Phi_0(x^*) - f^*)$. Combining the above with (\ref{eq:nsc:inexact:phi_quadratic}),
	$$ \phi_k^* + \frac{\gamma_k}{2}\Vert x^*  - \bar{v}(k)\Vert^2 \leq f^* + \lambda_k (\Phi_0(x^*) - f^*).$$
	Using the induction assumption, we get 
	\begin{equation}
	f(\bar{x}(k)) + \frac{\gamma_k}{2}\Vert x^*  - \bar{v}(k)\Vert^2 \leq f^* + 2\lambda_k (\Phi_0(x^*) - f^*).\label{eq:nsc:x_levelset}
	\end{equation}
	Since $f(\bar{x}(k))\geq f^*$ and $\gamma_k = \lambda_k\gamma_0$, we have $\Vert x^*  - \bar{v}(k)\Vert^2 \leq  \frac{4}{\gamma_0} (\Phi_0(x^*) - f^*)$. 
	Since $\bar{v}(k) =\frac{1}{\alpha_k}( \bar{y}(k) -  \bar{x}(k)) + \bar{x}(k)$, we have $\Vert \bar{v}(k) - x^*\Vert^2 = \Vert \frac{1}{\alpha_k}( \bar{y}(k) -  \bar{x}(k)) + \bar{x}(k) - x^*\Vert^2 \geq \frac{1}{2 \alpha_k^2} \Vert\bar{y}(k) - \bar{x}(k)\Vert^2 - \Vert\bar{x}(k) - x^*\Vert^2 $. By (\ref{eq:nsc:x_levelset}), $f(\bar{x}(k))\leq 2 \Phi_0(x^*) - f^* = 2f(\bar{x}(0)) - f^* + \gamma_0\Vert \bar{v}(0) - x^*\Vert^2 $. Also since $\gamma_0 = \frac{L}{1-\alpha_0}<2L$, we have $\bar{x}(k)$ lies within the ($2 f(\bar{x}(0)) - f^* + 2L\Vert\bar{v}(0) - x^*\Vert^2$)-level set of $f$. By Assumption \ref{assump:compact} and  \cite[Proposition B.9]{bertsekas1999nonlinear}, we have the level set is compact. Hence we have $\Vert \bar{x}(k) -x^*\Vert\leq R$ where $R$ is the diameter of that level set. Combining the above arguments, we get
{\small	\begin{align}
	\Vert \bar{y}(k) - \bar{x}(k)\Vert^2  
	&\leq 2 \alpha_k^2 (\Vert \bar{v}(k) - x^* \Vert^2 + \Vert \bar{x}(k)-x^*\Vert^2 )\nonumber\\
	&\leq 2 \alpha_k^2 (\frac{4}{\gamma_0} (\Phi_0(x^*) - f^*) + R^2 ) \nonumber\\
	&=  2\alpha_k^2 [ \frac{4}{\gamma_0} (f(\bar{x}(0)) - f^*) + 2 \Vert \bar{v}(0) - x^*\Vert^2 +R^2 ] \nonumber\\ 
	&\leq 2\alpha_k^2 \underbrace{[  4\Vert \bar{v}(0) - x^*\Vert^2 + R^2]}_{\triangleq C_1}  \label{eq:nsc:x_y_bound}
	\end{align}}where $C_1$ is a constant that does \textit{not} depend on $\eta$, and in the last inequality we have used by the $L$-smoothness of $f$, $f(\bar{x}(0)) - f^*\leq \frac{L}{2} \Vert \bar{x}(0) - x^*\Vert^2 \leq  \frac{\gamma_0}{2} \Vert \bar{x}(0) - x^*\Vert^2$.

	Next, we consider (\ref{eq:nsc:inexact:phi}), 
{\small	\begin{align}
	& \phi_{t+1}^*   \nonumber \\ 
	&= (1 - \alpha_t )\phi_t^*  + \alpha_t \hat{f}(t) - \frac{1}{2}\eta_t \Vert g(t)\Vert^2  +  \alpha_t \langle g(t), \bar{v}(t) - \bar{y}(t)\rangle \nonumber \nonumber\\   
	&= (1 - \alpha_t)( \phi_t^* - f(\bar{x}(t)))+(1 - \alpha_t )f(\bar{x}(t))   + \alpha_t \hat{f}(t)\nonumber \\
	&\qquad- \frac{1}{2}\eta_t \Vert g(t)\Vert^2  +  \alpha_t \langle g(t), \bar{v}(t) - \bar{y}(t)\rangle   \nonumber\\
	&\stackrel{(a)}{\geq} (1 - \alpha_t)( \phi_t^* - f(\bar{x}(t)))  + (1 - \alpha_t )\{ \hat{f}(t) + \langle g(t), \bar{x}(t) - \bar{y}(t)\rangle   \}  \nonumber\\
	&\qquad + \alpha_t \hat{f}(t) - \frac{1}{2}\eta_t  \Vert g(t)\Vert^2+  \alpha_t  \langle g(t), \bar{v}(t) - \bar{y}(t)\rangle  \nonumber\\
	&\stackrel{(b)}{=} (1 - \alpha_t)( \phi_t^* - f(\bar{x}(t)))+ \hat{f}(t)   - \frac{1}{2}\eta_t \Vert g(t)\Vert^2  \label{eq:nsc:phi_f_recursive}
	\end{align}}where (a) is due to (\ref{eq:nsc:inexact:noisy_oracle_1}) and (b) is due to $ \alpha_t (\bar{v}(t) - \bar{y}(t)) + (1-\alpha_t) (\bar{x}(t) - \bar{y}(t)) =0$.

	By (\ref{eq:nsc:inexact:noisy_oracle_2}) (setting $\omega = \bar{x}(t+1)$) and Lemma \ref{lem:nsc:rel_con_err}, 
{\footnotesize	\begin{align}
	f(\bar{x}(t+1))  
&\leq \hat{f}(t) + \langle g(t), \bar{x}(t+1) - \bar{y}(t)\rangle  + L\Vert \bar{x}(t+1) - \bar{y}(t)\Vert^2 \nonumber\\
&\quad + \frac{L}{n}\Vert y(t)-\one\bar{y}(t)\Vert^2\nonumber\\
	&\leq  \hat{f}(t) -(\eta_t -L\eta_t^2)\Vert g(t)\Vert^2 + 2 L \kappa^2 \chi_2(\eta_t)^2   [ L^2 \Vert   \bar{y}(t) - \bar{x}(t) \Vert^2 \nonumber\\
	&\quad +  \frac{64}{(1-\sigma)^2}L^2\eta_t^2 \Vert g(t)\Vert^2 ].\label{eq:nsc:f_upperbound}
	\end{align}}
	
	Combining the above with (\ref{eq:nsc:phi_f_recursive}) and recalling $\kappa =\frac{6}{1-\sigma}$, we get, 
{\color{black}{\small	\begin{align}
	&\phi_{t+1}^*- f(\bar{x}(t+1)) \nonumber\\
	 &\geq(1 - \alpha_t)( \phi_t^* - f(\bar{x}(t)))  \nonumber\\
	 &\qquad + (\frac{1}{2}\eta_t - L\eta_t^2 - \frac{4608 L^3 \chi_2(\eta_t)^2\eta_t^2 }{(1-\sigma)^4})\Vert g(t)\Vert^2 \nonumber\\
	&\qquad  - 2\kappa^2\chi_2(\eta_t)^2  L^3\Vert   \bar{y}(t) - \bar{x}(t) \Vert^2\nonumber\\
	&\geq (1 - \alpha_t)( \phi_t^* - f(\bar{x}(t)))  - 2\kappa^2\chi_2(\eta_t)^2  L^3\Vert \bar{y}(t) - \bar{x}(t) \Vert^2 \nonumber\\
	&\quad + \frac{1}{4} \eta_t \Vert g(t)\Vert^2 \label{eq:nsc:phi_f_recursive_2}
	\end{align}}}where in the last inequality we have used that, recalling $\chi_2(\eta_t) \leq \frac{2}{L^{2/3}} \eta_t^{1/3}$,
{\small	\begin{align*}
&	\frac{1}{2}\eta_t - L\eta_t^2 - \frac{4608 L^3 \chi_2(\eta_t)^2\eta_t^2 }{(1-\sigma)^4} \\
	&\geq \frac{1}{2}\eta_t - L\eta_t^2 - \frac{4608 L^3  \eta_t^2 }{(1-\sigma)^4} \frac{4 \eta_t^{2/3}}{L^{4/3}}\\
&= \frac{1}{2}\eta_t - L\eta_t^2 - \frac{18432  \eta_t (L  \eta_t)^{5/3} }{(1-\sigma)^4}  \\
	&\geq \eta_t(1/2  - \eta L - \frac{18432 (L  \eta)^{5/3} }{(1-\sigma)^4}  ) > \frac{1}{4} \eta_t
	\end{align*}}where the last inequality follows from $\eta L < \frac{1}{8}$, and $ \frac{18432 (L  \eta)^{5/3} }{(1-\sigma)^4}  < \frac{1}{8}$ ($\Leftarrow \eta L< \frac{(1-\sigma)^{2.4}}{1263} $ cf. step size condition (ii) in Theorem~\ref{thm:nsc:vanishing}). 
	Next, expanding (\ref{eq:nsc:phi_f_recursive_2}) recursively, we get
{\footnotesize	\begin{align}
	&\phi_{t+1}^*- f(\bar{x}(t+1))  \nonumber \\
	&\geq \prod_{k=0}^t(1-\alpha_k)( \phi_0^* - f(\bar{x}(0))) + \sum_{k=0}^t \frac{1}{4} \eta_k \Vert g(k)\Vert^2 \prod_{\ell=k+1}^{t}(1-\alpha_\ell) \nonumber \\
	&\qquad    -\sum_{k=0}^t 2 \kappa^2\chi_2(\eta_k)^2  L^3 \Vert  \bar{y}(k)- \bar{x}(k)  \Vert^2\prod_{\ell=k+1}^{t}(1-\alpha_\ell) \nonumber\\
	& = \sum_{k=0}^t \frac{1}{4} \eta_k \Vert g(k)\Vert^2 \prod_{\ell=k+1}^{t}(1-\alpha_\ell) \nonumber\\
&	\quad -\sum_{k=0}^t 2 \kappa^2\chi_2(\eta_k)^2  L^3 \Vert  \bar{y}(k)- \bar{x}(k)  \Vert^2\prod_{\ell=k+1}^{t}(1-\alpha_\ell) \label{eq:nsc:phi_f_for_individual_err}\\
&\geq -\sum_{k=0}^t 2 \kappa^2\chi_2(\eta_k)^2  L^3 \Vert  \bar{y}(k)- \bar{x}(k)  \Vert^2\prod_{\ell=k+1}^{t}(1-\alpha_\ell).\nonumber
	\end{align}}
	Therefore to finish the induction, we need to show
{\small	\begin{align}
&	\sum_{k=0}^t 2 \kappa^2\chi_2(\eta_k)^2  L^3 \Vert   \bar{y}(k) - \bar{x}(k)\Vert^2\prod_{\ell=k+1}^t(1-\alpha_\ell) \nonumber\\
& \leq  (\Phi_0(x^*) - f^*)\lambda_{t+1}. \label{eq:nsc:phi_f_for_individual_err_2}
	\end{align}}Notice that 
{\footnotesize	\begin{align*}
&	\sum_{k=0}^t \frac{2\kappa^2\chi_2(\eta_k)^2  L^3 \Vert \bar{y}(k) -  \bar{x}(k) \Vert^2\prod_{\ell=k+1}^t(1-\alpha_\ell)}{(\Phi_0(x^*) - f^*)\lambda_{t+1}}  \\
	&\stackrel{(a)}{\leq} \sum_{k=0}^t  \frac{4\kappa^2  L^2 }{ \Vert\bar{v}(0)-x^*\Vert^2 }\chi_2(\eta_k)^2\Vert \bar{y}(k) -  \bar{x}(k)\Vert^2 \frac{1}{\lambda_{k+1}}\\
	&\stackrel{(b)}{\leq} \sum_{k=0}^t  \frac{4(6/(1-\sigma))^2  L^2 }{ \Vert\bar{v}(0)-x^*\Vert^2 } (\frac{2}{L^{2/3}} \eta_k^{1/3} )^2  2C_1\alpha_k^2 \frac{1}{\lambda_{k+1}}\\
	&= \sum_{k=0}^t \underbrace{\frac{1152 L^{2/3} C_1}{(1-\sigma)^2\Vert\bar{v}(0)-x^*\Vert^2}   }_{\triangleq C_2} \eta_k^{2/3} \alpha_k^2 \frac{1}{\lambda_{k+1}}
	\end{align*}}where $C_2$ is a costant that does \textit{not} depend on $\eta$, and in (a) we have used $\Phi_0(x^*) - f^* \geq \frac{L}{2}\Vert \bar{v}(0) - x^*\Vert^2>0$, and $\prod_{\ell=k+1}^t(1-\alpha_\ell) = \lambda_{t+1}/\lambda_{k+1}$; in (b), we have plugged in $\kappa = \frac{6}{1-\sigma}$, used $\chi_2(\eta_k) \leq \frac{2\eta_k^{1/3}}{L^{2/3}}$ (cf. Lemma \ref{lem:nsc:rel_con_err}) and the bound on $\Vert   \bar{y}(k) -\bar{x}(k)\Vert^2$ (equation (\ref{eq:nsc:x_y_bound})). Now by Lemma \ref{lem:nsc:alpha}, we get, 
{\small	\begin{align}
\sum_{k=0}^t\eta_k^{2/3} \alpha_k^2 \frac{1}{\lambda_{k+1}} & \leq \sum_{k=0}^t  \frac{\eta^{2/3}}{(k+t_0)^{\frac{2}{3}\beta}}  \frac{4}{(k+1)^2} \frac{(k+1+t_0)^{2-\beta}}{ D(\beta,t_0)}\nonumber\\
&\stackrel{(a)}{\leq} \eta^{2/3} \frac{4 (t_0+1)^{2-\beta}}{D(\beta,t_0)}  \sum_{k=0}^\infty \frac{1}{(k+1)^{\frac{5}{3}\beta}}\nonumber\\
&\stackrel{(b)}{\leq} \eta^{2/3} \frac{4 (t_0+1)^{2-\beta}}{D(\beta,t_0)} \times \frac{2}{\beta - 0.6}\label{eq:nsc:covthm4:sumstepsize}
	\end{align}}where in (a) we have used, $k+t_0 \geq k+1$, $k+1+t_0\leq (t_0+1)(k+1)$; in (b) we have used $\frac{5}{3}\beta > 1$. So, we have 
{\small	\begin{align*}
	&\sum_{k=0}^t \frac{2\kappa^2\chi_2(\eta_k)^2  L^3 \Vert  \bar{y}(k)-\bar{x}(k)\Vert^2\prod_{\ell=k+1}^t(1-\alpha_\ell)}{(\Phi_0(x^*) - f^*)\lambda_{t+1}}  \\
	&	\leq   C_2 \eta^{2/3} \frac{8 (t_0+1)^{2-\beta}}{D(\beta,t_0) (\beta - 0.6)}   < 1
	\end{align*}}where in the last inequality, we have simply required $\eta^{2/3} < \frac{D(\beta,t_0) (\beta - 0.6)}{8 (t_0+1)^{2-\beta} C_2   }$ (step size condition (iii) in Theorem~\ref{thm:nsc:vanishing}), which is possible since the constants $C_2$ and $D(\beta,t_0)$ do not depend on $\eta$. So the induction is complete and we have (\ref{eq:nsc:inexact:f_y_phi}) is true. 
{\color{black}Part (b) of the theorem, i.e. $f(y_i(t)) - f^* = O(\frac{1}{t^{1.4-\epsilon}})$ will be an easy corollary. Its proof can be found in \ifthenelse{\boolean{appendixin}}{Appendix-\ref{appendix:theorem4b}.}{\cite[Appendix-H]{fullversion}.}
}

\section{Numerical Experiments}\label{sec:numerical}


\textbf{Graphs and Matrix $W$:} We consider three different graphs. \textbf{Random graph:} the graph has $n=100$ agents and is generated using the Erdos-Renyi model \cite{erdds1959random} with connectivity probability $0.3$. \textbf{$k$-cycle:} the graph has $n=100$ agents and it is a $k$-cycle with $k = 20$, i.e. we arrange the agents into a cycle, and each agent is connected to $20$ agents to its left and $20$ agents to its right. \textbf{2D-grid:} the graph has $n=25$ nodes and is a $5\times 5$ 2-D grid. The weight matrix $W$ is chosen using the Laplacian method \cite[Sec 2.4. (ii)]{shi2015extra}. 


\textbf{Cost functions: } For the functions $f_i$, we consider three cases. 
\textbf{Case I:} The functions $f_i$ are square losses for linear regression, i.e. $
f_i(x) = \frac{1}{M_i}\sum_{m=1}^{M_i} ( \langle u_{im},x\rangle - v_{im} )^2 $
 where $u_{im}\in\R^N$ ($N=3$) are the features and $v_{im}\in\R$ are the observed outputs, and $\{(u_{im},v_{im})\}_{m=1}^{M_i}$ are $M_i=50$ data samples for agent $i$.
  We fix a predefined parameter $\tilde{x}\in\R^N$ with each element drawn uniformly from $[0,1]$. For each sample $(u_{im}, v_{im})$, the last element of $u_{im}$ is fixed to be $1$, and the rest elements are drawn from i.i.d. Gaussian with mean $0$ and variance $400$. Then we generate $v_{im} = \langle \tilde{x}, u_{im}\rangle + \epsilon_{im}$ where $\epsilon_{im}$ are independent Gaussian noises with mean $0$ and variance $100$. 
\textbf{Case II: } $f_i$ is the loss function for logistic regression \cite{logit_regression}, i.e. $f_i(x) = \frac{1}{M_i}\sum_{m=1}^{M_i} \big[\ln(1 + e^{\langle u_{im},x\rangle }) - v_{im}\langle u_{im},x\rangle \big] $ where $u_{im}\in\R^N$ ($N=3$) are the features and $v_{im}\in\{0,1\}$ are the observed labels, and $\{(u_{im},v_{im})\}_{m=1}^{M_i}$ are $M_i=100$ data samples for agent $i$. We first fix a predefined parameter $\tilde{x}\in\R^N$ with each element drown uniformly from $[0,1]$. For each sample $(u_{im}, v_{im})$, the last element of $u_{im}$ is fixed to be $1$, and the rest elements of $u_{im}$ are drawn from i.i.d. Gaussian with mean $0$ and variance $100$. We then generate $v_{im}$ from a Bernoulli distribution, with probability of $v_{im}=1$ being $\frac{1}{1 + e^{-\langle \tilde{x},u_{im} \rangle} }$. 
\textbf{Case III:} the objective functions are given by,
\begin{align*}
f_i(x) = \left\{ \begin{array}{ll}
\frac{1}{m} \langle a_i, x\rangle^m + \langle b_i,x\rangle & \text{ if } | \langle a_i, x\rangle|\leq 1,\\
|\langle a_i, x\rangle| - \frac{m-1}{m} + \langle b_i,x\rangle & \text{ if } | \langle a_i, x\rangle|> 1,
\end{array}\right.
\end{align*}
where $m=12$, $a_i,b_i\in\R^N$ ($N = 4$) are vectors whose entries are i.i.d. Gaussian with mean $0$ and variance $1$, with the exception that $b_n$ is set to be $b_n = - \sum_{i=1}^{n-1} b_i$ s.t. $\sum_i b_i = 0$. It is easy to check that $f_i$ is convex and smooth, but not strongly convex (around the minimizer). 

\textbf{Algorihtms to compare:} We will compare our algorithm (Acc-DNGD-SC or Acc-DNGD-NSC) with the ``D-NG'' method in \cite{jakovetic2014fast}, Distributed Gradient Descent (DGD) in \cite{nedic2009distributed} with a vanishing step size, the ``EXTRA'' algorithm in \cite{shi2015extra} (with $\tilde{W} = \frac{W+I}{2}$), the algorithm studied in \cite{xu2015augmented,di2015distributed,di2016next,qu2016harnessing,nedich2016achieving,xi2016add,nedic2016geometrically} (we name it ``Acc-DGD'').  We will also compare with two centralized methods that directly optimize $f$: Centralized Gradient Descent (CGD) and Centralized Nesterov Gradient Descent (CNGD-SC (\ref{eq:cngd_s}) or CNGD-NSC (\ref{eq:nsc:cngd})). Each element of the initial point $x_i(0)$ is drawn from i.i.d. Gaussian with mean $0$ and variance $25$. 

{\color{black}\textbf{Step size selection: } For D-NG, $\eta_t = \frac{1}{2L(t+1)}$; for DGD, $\eta_t = \frac{1}{L\sqrt{t}}$; for CGD and CNGD, $\eta= \frac{1}{L}$; for EXTRA (convex and smooth case), $\eta = \frac{1}{L}$. The above step size selections all follow the guideline in the respective paper. For EXTRA (strongly convex and smooth case), Acc-DGD and Acc-DNGD-SC, the bounds in the respective paper are known to be not tight, and therefore we do trial and error to get a step size that has the fastest convergence rate. Lastly, for Acc-DNGD-NSC, we test both vanishing step size $\eta_t = \frac{1}{2L (t+1)^{0.61}}$ and fixed step size $\eta_t = \eta$ where the fixed step size is obtained through trial and error. The exact values of the step sizes used in this section will be reported in detail in \ifthenelse{\boolean{appendixin}}{Appendix-\ref{appendix:stepsize}.}{\cite[Appendix-J]{fullversion}.}}

\begin{figure*}[t]
	\begin{center}
\begin{subfigure}[b]{0.3\textwidth}
	\centering
\includegraphics[width=\textwidth]{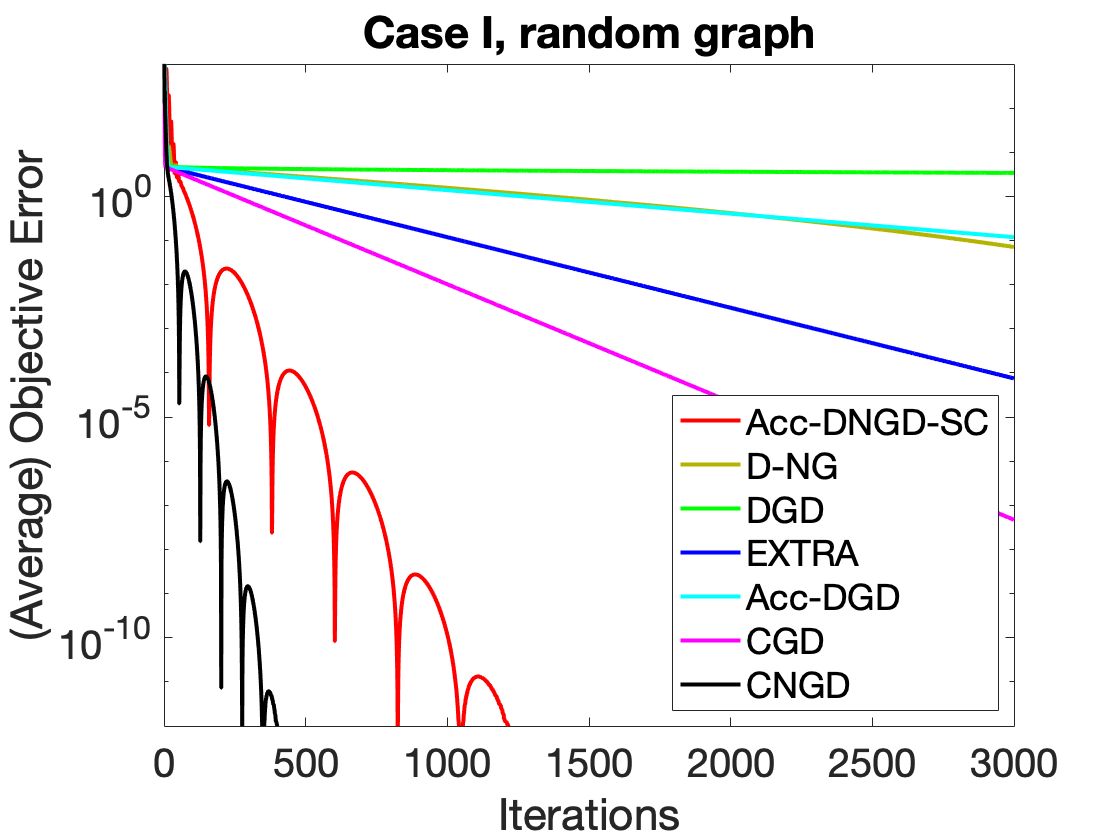} 
\end{subfigure}
\begin{subfigure}[b]{0.3\textwidth}
	\centering
	\includegraphics[width=\textwidth]{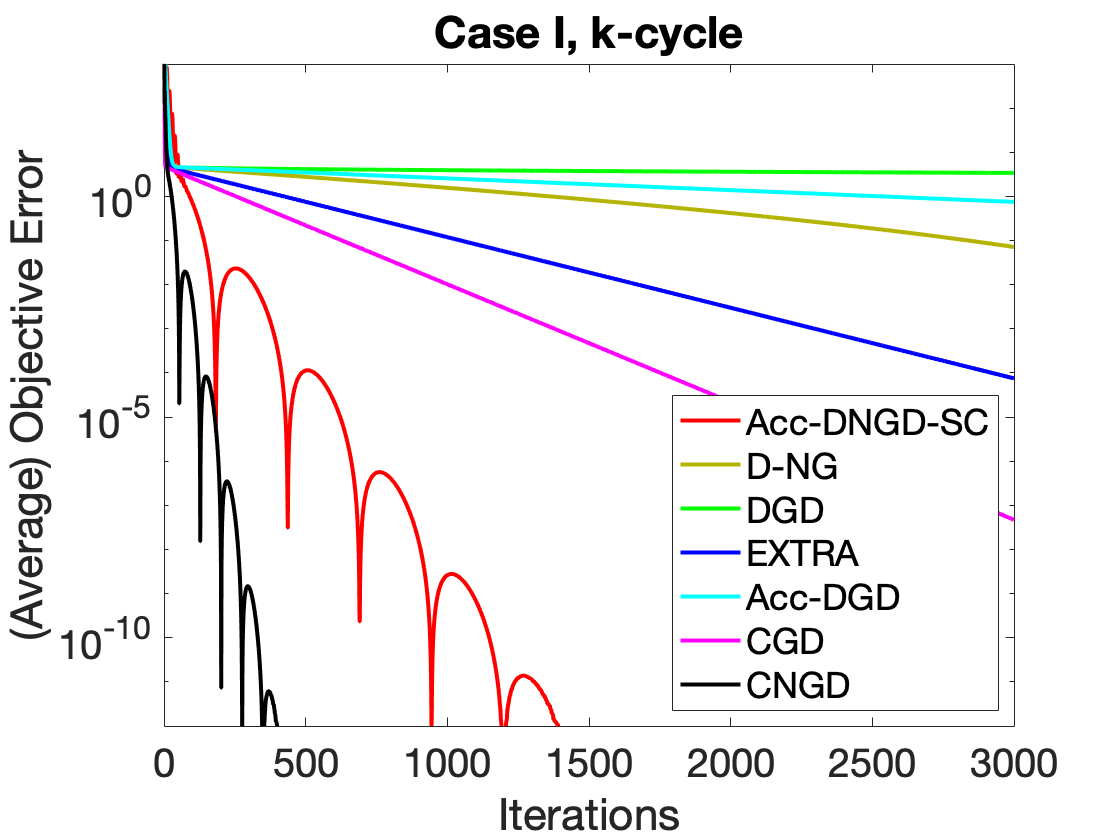} 
\end{subfigure}
\begin{subfigure}[b]{0.3\textwidth}
	\centering
	\includegraphics[width=\textwidth]{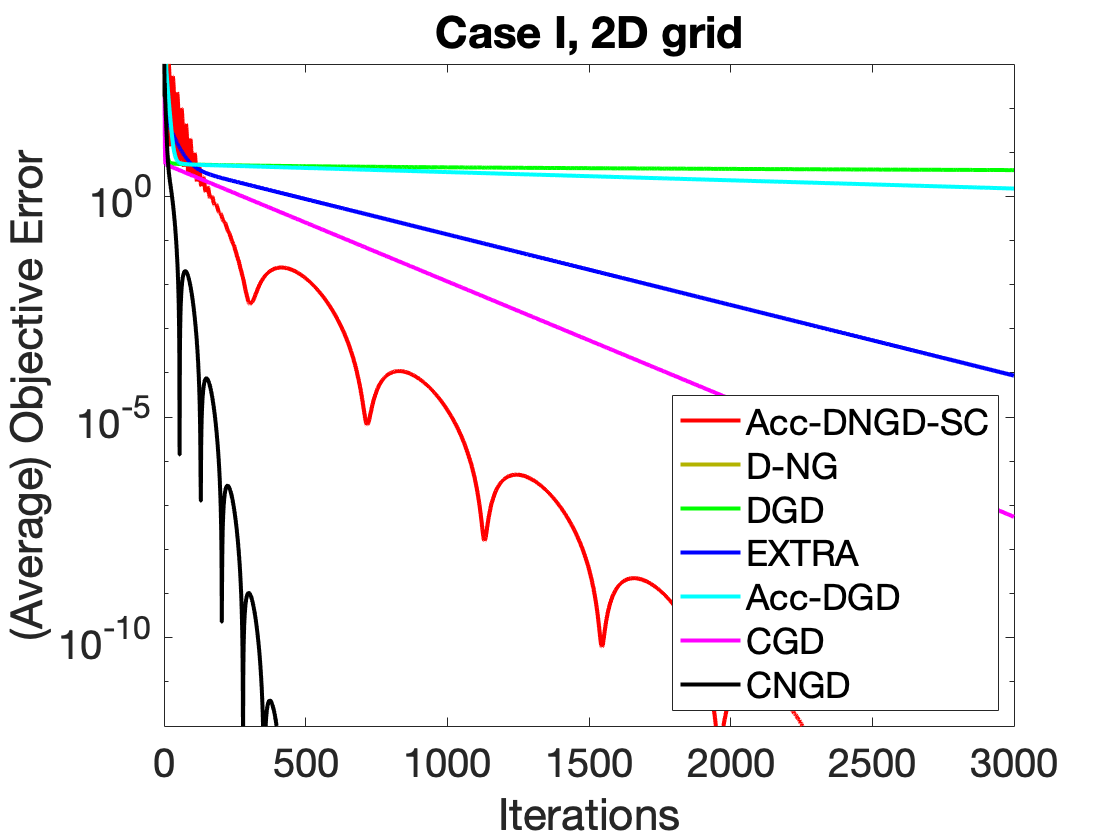} 
\end{subfigure}
\end{center}
	\caption{Case I. (a) random graph: $\frac{L}{\mu} \approx 793, \sigma \approx 0.59$; 
		 (b) $k$-cycle: $\frac{L}{\mu} \approx 793, \sigma \approx 0.75$; 
		 (c) 2D grid: $\frac{L}{\mu} \approx 773, \sigma \approx 0.92$. 
	 }\label{fig:case1}
\end{figure*}

\begin{figure*}[t]
	\begin{center}
\begin{subfigure}[b]{0.3\textwidth}
	\centering
	\includegraphics[width=\textwidth]{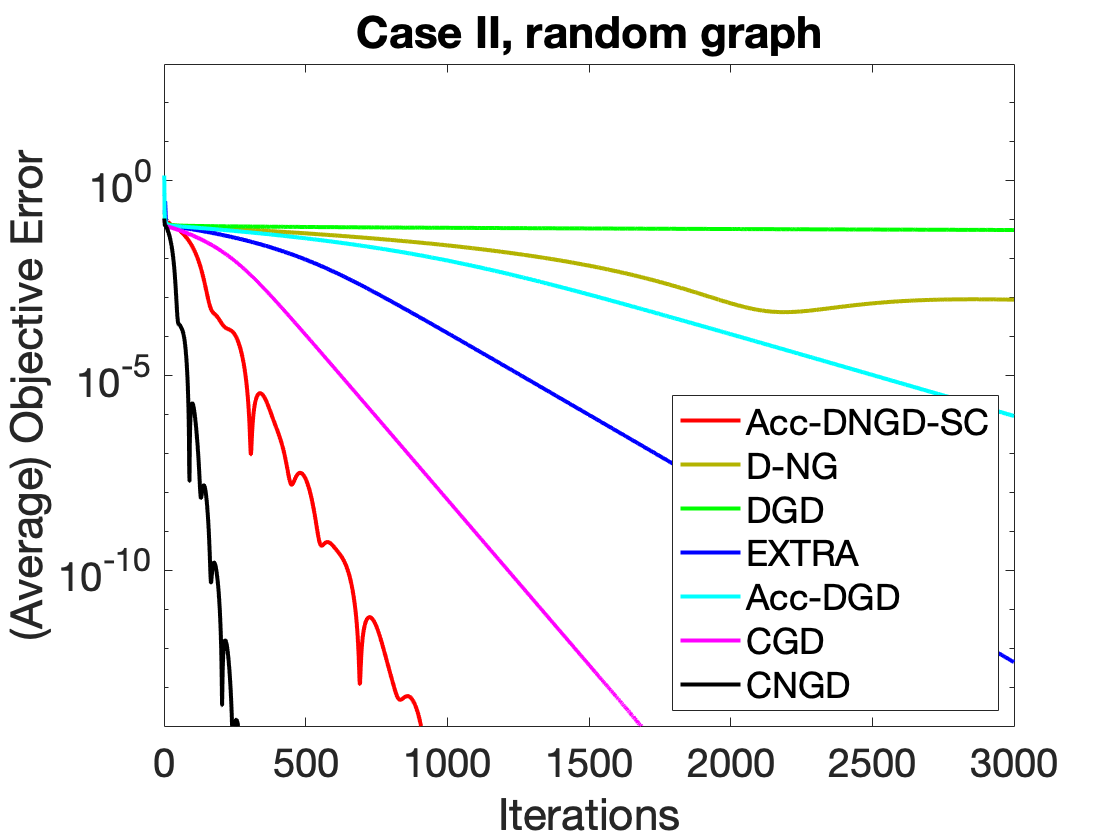} 
\end{subfigure}
\begin{subfigure}[b]{0.3\textwidth}
	\centering
	\includegraphics[width=\textwidth]{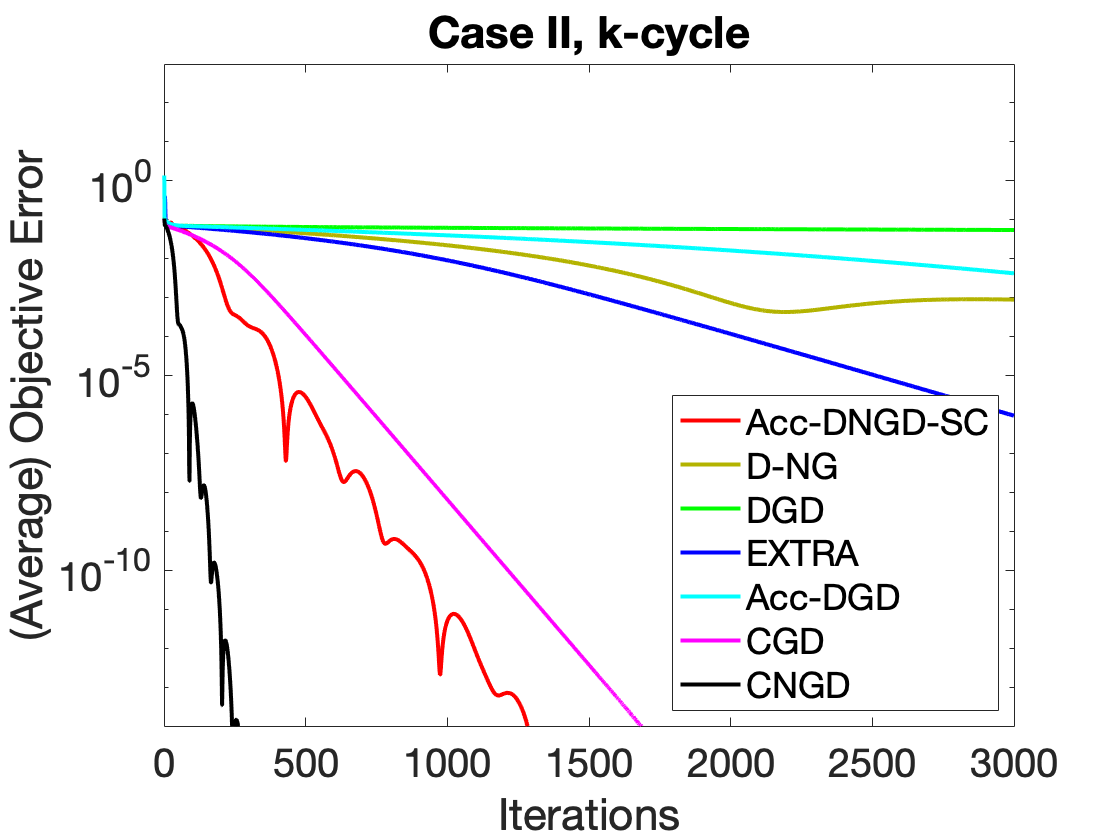} 
\end{subfigure}
\begin{subfigure}[b]{0.3\textwidth}
	\centering
	\includegraphics[width=\textwidth]{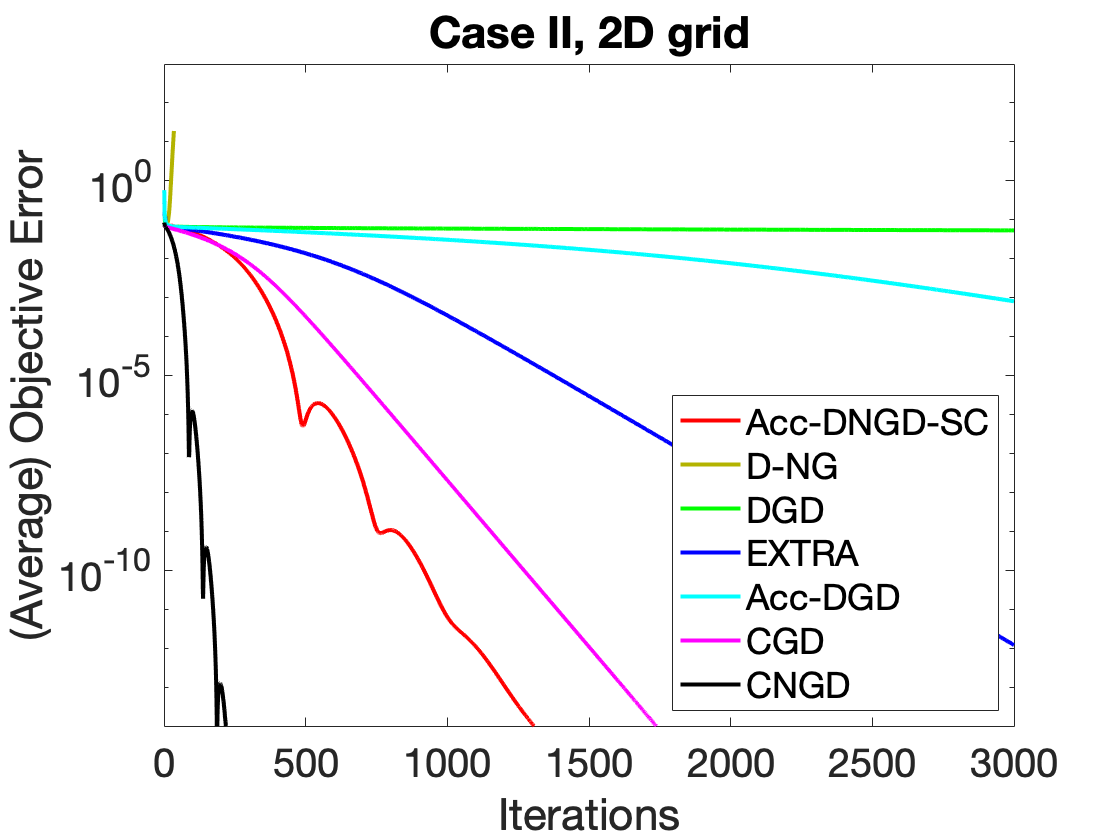} 
\end{subfigure}
\end{center}
	\caption{Case II. (a) random graph: $\frac{L}{\mu} \approx 658, \sigma = 0.59$;
		 (b) $k$-cycle: $\frac{L}{\mu} \approx 658, \sigma \approx 0.75$;
		 (c) 2D grid: $\frac{L}{\mu} \approx 345, \sigma \approx 0.92$. 
	 }\label{fig:case2}
\end{figure*}

\textbf{Simulation Results for Case I and II.} In Case I and II, the functions are strongly convex and smooth, so we test the strongly-convex variant of our algorithm Acc-DNGD-SC and the version of CNGD we compare with is CNGD-SC. 
The simulation results are shown in Figure \ref{fig:case1} and \ref{fig:case2} for all the three graphs, where the $x$-axis is iteration number $t$, and the $y$-axis is the {\color{black}average objective error for distributed methods ($\frac{1}{n} \sum f(y_i(t)) - f^*$ for the proposed method, $\frac{1}{n} \sum f(x_i(t)) - f^*$ for other distributed methods), or objective error $f(x(t)) - f^*$ for centralized methods.} It can be seen that our algorithm Acc-DNGD-SC indeed performs significantly better than CGD, CGD-based distributed methods (DGD, EXTRA, Acc-DGD) and also D-NG. \textcolor{black}{Also, recall that our theoretic results show that our methods have better \textit{asymptotic} dependence on condition number $\frac{L}{\mu}$ in the large $\frac{L}{\mu}$ regime than CGD and existing distributed methods. The simulation results show that when the condition number is between $300$ to $800$, our method can already outperform CGD and existing distributed methods.   }

\begin{figure*}[t]
	\begin{center}
\begin{subfigure}[b]{0.32\textwidth}
	\centering
	\includegraphics[width=\textwidth]{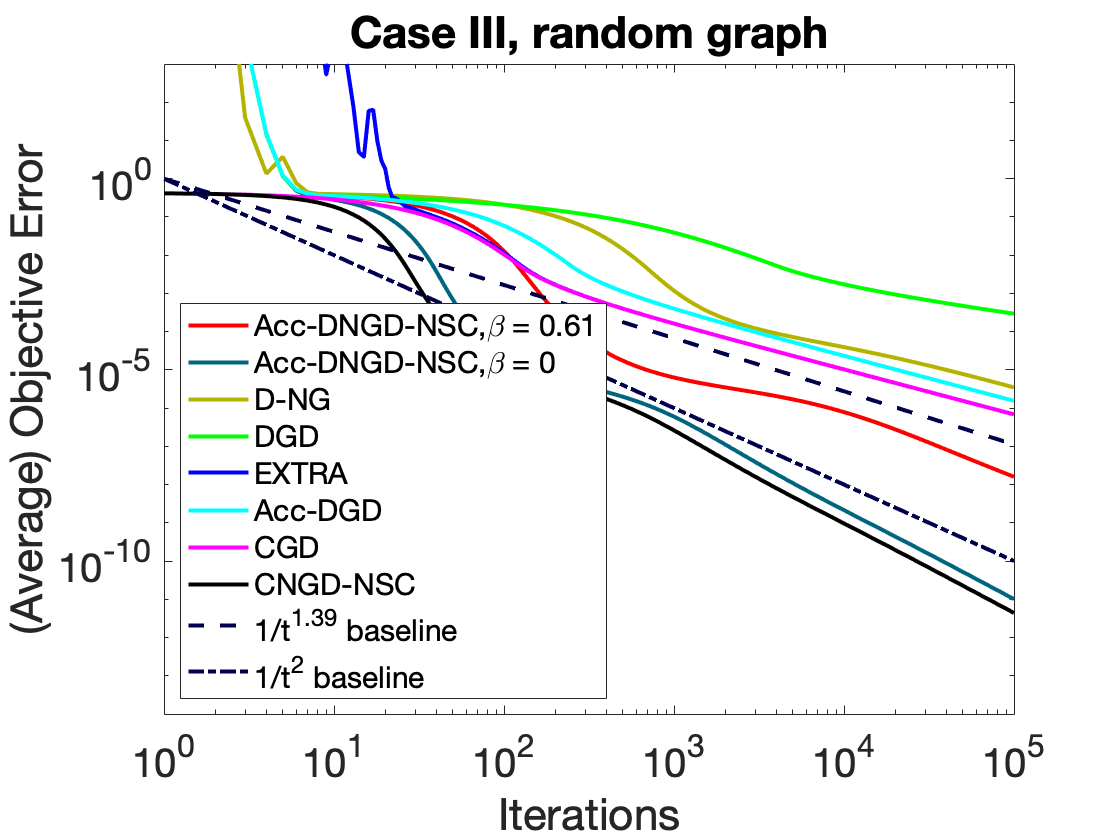} 
\end{subfigure}
\begin{subfigure}[b]{0.32\textwidth}
	\centering
	\includegraphics[width=\textwidth]{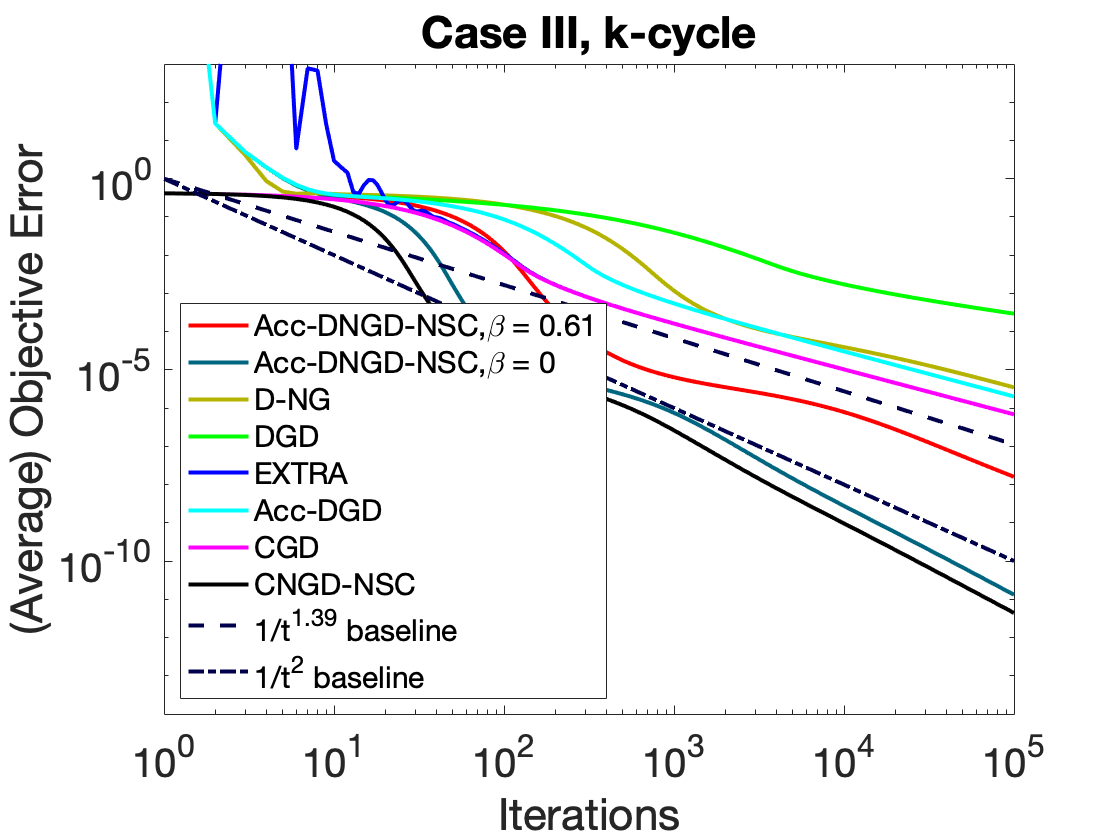} 
\end{subfigure}
\begin{subfigure}[b]{0.32\textwidth}
	\centering
	\includegraphics[width=\textwidth]{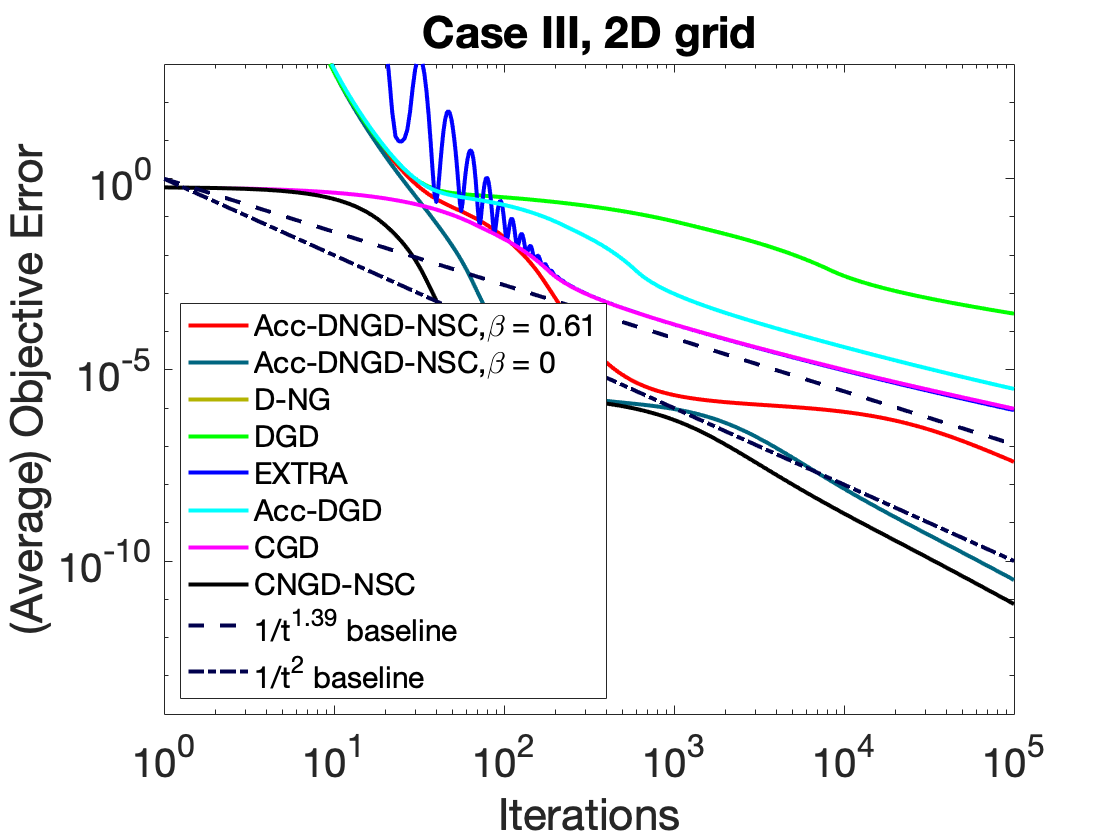} 
\end{subfigure}
\end{center}
	\caption{Case III. (a) For random graph, $\sigma \approx 0.59$; 
		(b) For $k$-cycle, $\sigma \approx 0.75$; 
		(c) For 2D grid, $\sigma \approx 0.92$. 
	  }\label{fig:case3}
\end{figure*}

\textbf{Simulation Results for Case III.} Case III is to test the sublinear convergence rate $\frac{1}{t^{2-\beta}}$ ($\beta > 0.6$) of the Acc-DNGD-NSC (\ref{eq:nsc:update_elementwise}) and the conjecture that the $\frac{1}{t^{2}}$ rate still holds even if $\beta=0$ (i.e. fixed step size, cf. Theorem~\ref{thm:nsc:vanishing} and the comments following it). Therefore, we do two runs of Acc-DNGD-NSC, one with $\beta = 0.61$ and the other with $\beta = 0$. The results are shown in Figure~\ref{fig:case3}, where the $x$-axis is the iteration $t$, and the $y$-axis is the (average) objective error. 
It shows that Acc-DNGD-NSC with $\beta = 0.61$ is faster than $\frac{1}{t^{1.39}}$, while D-NG, CGD and CGD-based distributed methods (DGD, Acc-DGD, EXTRA) are slower than $\frac{1}{t^{1.39}}$. Further, both Acc-DNGD-NSC with $\beta = 0$ and CNGD-NSC are faster than $\frac{1}{t^2}$. 

\begin{figure}[t]
	\centering
	\includegraphics[width=0.7\columnwidth]{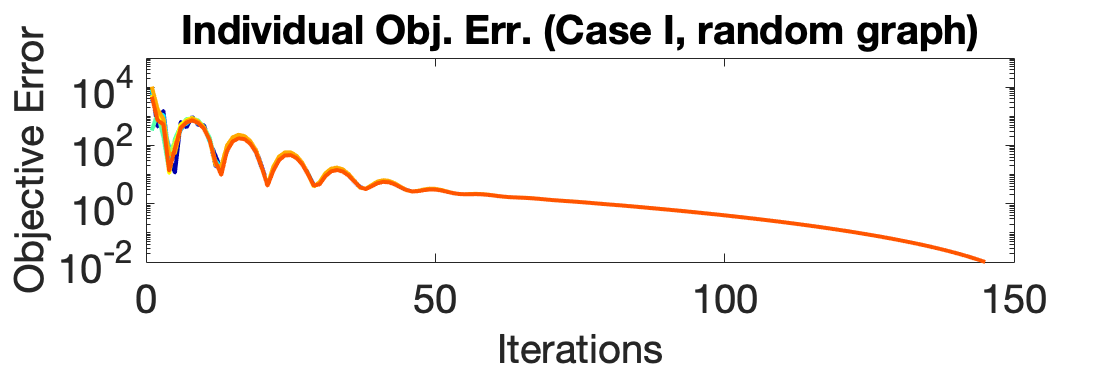}
	\includegraphics[width=0.7\columnwidth]{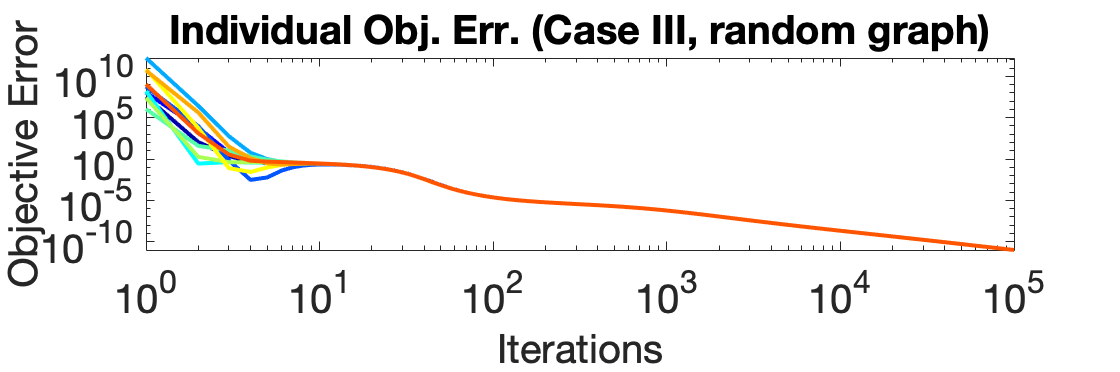}
	\caption{Individual objective errors. 
	}\label{fig:individual_err}
\end{figure}

{\color{black}\textbf{Convergence from each agent's perspective.} For case I and case III with random graph, we also plot the objective error $f(y_i(t)) - f^*$ of our algorithm (Acc-DNGD-SC for case I and Acc-DNGD-NSC with fixed step size for case III) for agent $i=10,20,\ldots,100$ in Figure~\ref{fig:individual_err}. The figure shows that the individual objective errors mix very fast and become indistinguishable after about 100 iterations. }

\begin{figure}[t]
\centering
\includegraphics[width=0.7\columnwidth]{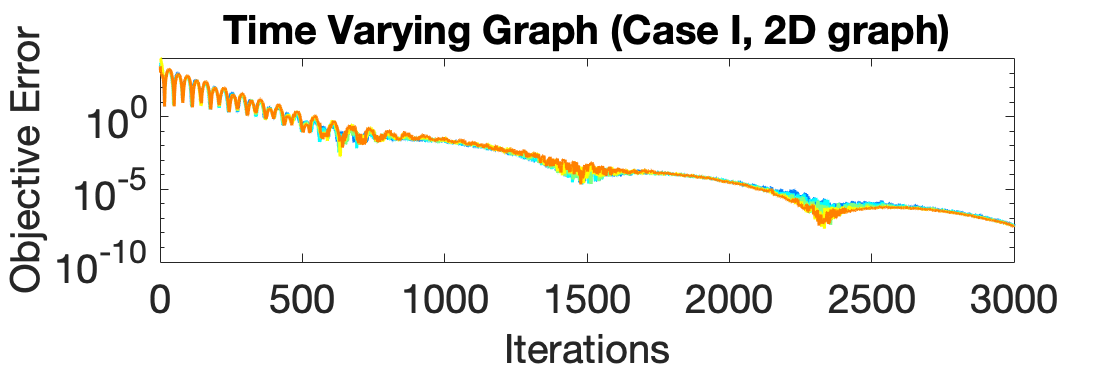}
\includegraphics[width=0.7\columnwidth]{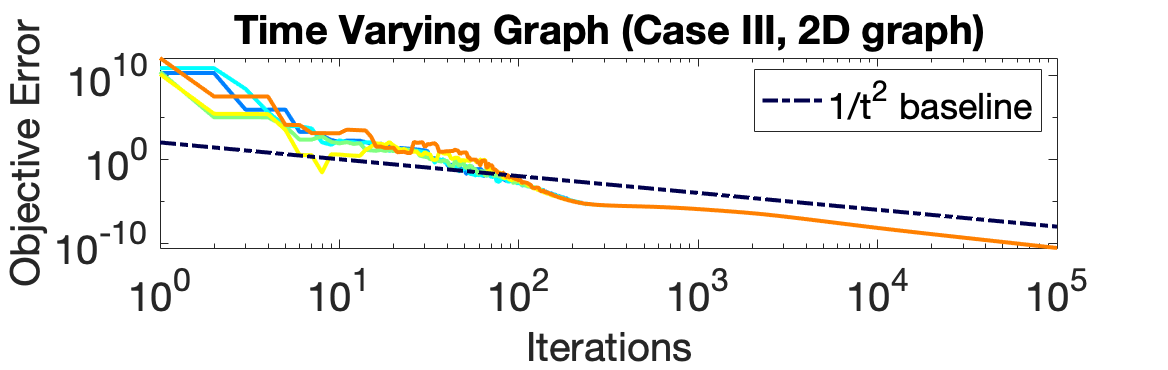}
\caption{Time varying graph. 
}\label{fig:timevarying}
\end{figure}

{\color{black}\textbf{Simulation on time-varying graphs.} We test our algorithm on a time-varying graph. We use the cost function of Case I and case III for strongly convex case (using Acc-DNGD-SC) and non-strongly convex case (using Acc-DNGD-NSC with fixed step size) respectively. We let the ground graph be the $5\times 5$ 2D grid. At each iteration $t$, we randomly select a $0.75$ portion of the edges from the ground graph and remove the selected edges.
	 Correspondingly at each iteration $t$ the agents recalculate the weights $w_{ij}$ according to the current graph using the method in \cite[page 8, eq. after Assump. 3]{nedich2016achieving}, and then implement our algorithm. The simulation results are shown in Figure~\ref{fig:timevarying}, where for both cases we plot $f(y_i(t)) - f^*$ for $5$ randomly selected $i$'s. Figure~\ref{fig:timevarying} shows our algorithm converges even when the graph is time varying. It also shows that in the time varying graph case, while overall the agents' objective errors converge to $0$, their trajectories are more volatile compared to Figure~\ref{fig:individual_err}.} 

\begin{figure}[t]
	\centering
	\includegraphics[width=0.7\columnwidth]{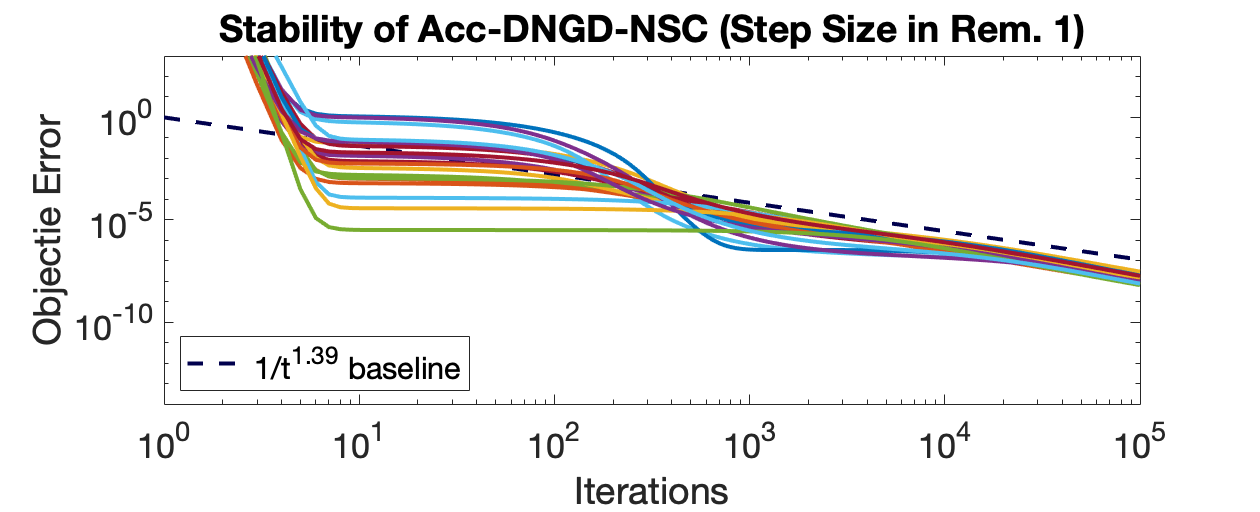}
	\caption{\color{black}Stability of the step size rule in Remark~\ref{rem:vanishing_stepsize}. }\label{fig:stability_rem1}
\end{figure}

{\color{black} \textbf{Stability of the step size rule in Remark~\ref{rem:vanishing_stepsize}.} We verify the stability of Acc-DNGD-NSC under the step size rule $\eta_t = \frac{1}{2L(t+1)^\beta}$ proposed in Remark~\ref{rem:vanishing_stepsize}. To this end, we use the same setting as Case III, random graph (Figure \ref{fig:case3} (a)) and we generate $20$ random instances (including the function parameters and the graph). We run Acc-DNGD-NSC on the instances using the step size in Remark~\ref{rem:vanishing_stepsize} with $\beta = 0.61$, and plot the average objective error trajectory of the $20$ runs in Figure~\ref{fig:stability_rem1}. Figure~\ref{fig:stability_rem1} shows that Acc-DNGD-NSC is stable under the step size rule in Remark~\ref{rem:vanishing_stepsize}. Moreover, the objective error decays faster than the baseline $O(\frac{1}{t^{1.39}})$, further confirming the results in Theorem~\ref{thm:nsc:vanishing}.}
	
\section{Conclusion} \label{sec:conclusion}
In this paper we have proposed an Accelerated Distributed Nesterov Gradient Descent (Acc-DNGD) method. The first version works for convex and $L$-smooth functions and we show that it achieves a $O(\frac{1}{t^{1.4-\epsilon}})$ convergence rate for all $\epsilon\in(0,1.4)$. We also show the convergence rate can be improved to $O(\frac{1}{t^2})$ if the objective function is a composition of a linear map and a strongly-convex and smooth function. The second version works for $\mu$-strongly convex and $L$-smooth functions, and we show that it achieves a linear convergence rate of $O([ 1 - C  (\frac{\mu}{L})^{5/7} ]^t)$ for some constant $C$ independent of $L$ or $\mu$. All the rates are better than CGD and CGD-based distributed methods. In the future, we plan to tighten our analysis to obtain better convergence rates. 
\bibliographystyle{IEEEtran}
\bibliography{con_opt_ref}

	\appendix

\section{Proofs of Auxiliary Lemmas in Section~\ref{sec:convergence_sc}} 
\subsection{Proof of Lemma \ref{prop:simp_eq} and Lemma \ref{prob:nes:ineq_useful}}\label{subsec:useful_fact} \ 

\noindent\textit{Proof of Lemma \ref{prop:simp_eq}:} 
We have, 
\begin{align*}
&\Vert \bar{y}(t+1) - \bar{y}(t)\Vert \\
&=\Vert \frac{1}{1+\alpha} [\bar{x}(t+1) - \bar{x}(t)] + \frac{\alpha}{1+\alpha} [\bar{v}(t+1) - \bar{v}(t)]\Vert \\
&=\Vert  \frac{1}{1+\alpha}[\bar{y}(t)-\bar{x}(t)] - \frac{\eta}{1+\alpha} g(t)   + \frac{\alpha}{1+\alpha} \alpha (\bar{y}(t) - \bar{v}(t) ) - \frac{\eta}{1+\alpha} g(t) \Vert \\
&= \Vert \frac{1}{1+\alpha}[\bar{y}(t)-\bar{x}(t)]  + \frac{\alpha}{1+\alpha} \alpha [\bar{y}(t) -    \frac{1+\alpha}{\alpha}\bar{y}(t) + \frac{1}{\alpha} \bar{x}(t) ]   - \frac{2\eta}{1+\alpha}  g(t)\Vert \\
&\leq  \frac{1-\alpha}{1+\alpha}\Vert\bar{y}(t)-\bar{x}(t)\Vert +  \frac{2\eta}{1+\alpha}  \Vert g(t)\Vert \\
&\leq \Vert\bar{y}(t)-\bar{x}(t)\Vert + 2\eta \Vert g(t)\Vert.
\end{align*}
\qedd

\noindent\textit{Proof of Lemma \ref{prob:nes:ineq_useful}: }
For (\ref{eq:nes:ineq_useful_a}), we have
\begin{align*}
\Vert \nabla(t+1) - \nabla(t) \Vert &= \sqrt{\sum_{i=1}^n \Vert \nabla f_i(y_i(t+1)) - \nabla f_i(y_i(t))\Vert^2}\\
&\leq  \sqrt{\sum_{i=1}^n L^2 \Vert y_i(t+1) - y_i(t)\Vert^2}\\
&= L \Vert y(t+1) - y(t)\Vert.
\end{align*}
For (\ref{eq:nes:ineq_useful_b}), we have 
\begin{align*}
\Vert g(t) - \nabla f(\bar{y}(t))\Vert & = \Vert  \frac{\sum_{i=1}^n (\nabla f_i(y_i(t)) - \nabla f_i (\bar{y}(t)))}{n}\Vert\\
&\leq \frac{1}{n} \sum_{i=1}^n \Vert \nabla f_i(y_i(t)) - \nabla f_i (\bar{y}(t))\Vert\\
&\leq \frac{L}{n} \sum_{i=1}^n \Vert y_i(t)  - \bar{y}(t)\Vert\\
&\leq L \sqrt{\frac{1}{n} \sum_{i=1}^n \Vert y_i(t)  - \bar{y}(t)\Vert^2}\\
& = \frac{L}{\sqrt{n}} \Vert y(t) - \one \bar{y}(t)\Vert.
\end{align*}
\qedd

\subsection{Proof of Lemma \ref{lem:G_eig}}\label{subsec:G}
In this section, we provide the proof of Lemma \ref{lem:G_eig}.

\noindent\textit{Proof of Lemma \ref{lem:G_eig}: }
(a) We first calculate the characteristic polynomial of $G(\eta)$ as 
{\small\begin{align*}
p(\zeta) &= \overbrace{(\zeta-\sigma)(\zeta - \frac{1-\alpha}{1+\alpha} \sigma)(\zeta - \sigma - 2\eta L) }^{\triangleq p_0(\zeta)} \\
& \quad  - \overbrace{[k_1(\zeta - \sigma - 2\eta L) + k_2 ]}^{\triangleq p_1(\zeta)}
\end{align*}}where $k_1$ and $k_2$ are positive constants given by
\begin{align*}
k_1 = 4\eta L + \sigma\eta L < 5\eta L
\end{align*}
and 
\begin{align*}
k_2 
&= 2\eta^2 L^2(4+\sigma) + \eta L\sigma (2+ \alpha\sigma) + \frac{2\eta L \alpha^2 \sigma (2+\alpha \sigma )}{1+\alpha}\\
&< 10 \eta^2 L^2 + 3\eta L + 6\eta L   < 19\eta L.
\end{align*}
where we have used $\eta L<1$.
 Let $\zeta_0 = \sigma + 4(\eta L)^{1/3}$. Then $\zeta_0 > \sigma + 2\eta L$. Since $p(\sigma + 2\eta L)<0 $, and $p(\zeta)$ is a strictly  increasing function on $[\zeta_0,+\infty)$ (since when $\zeta\geq\zeta_0$, $p'(\zeta) > (\zeta_0 - \sigma)(\zeta_0 - \frac{1-\alpha}{1+\alpha} \sigma) - k_1 > 16(\eta L)^{2/3}  - 5\eta L>0$), and also since $G(\eta)$'s largest eigenvalue in magnitude must be a positive real number (Perron-Frobenius Theorem \cite[Theorem 8.2.11]{horn2012matrix}), we have $\rho(G(\eta))$ is a real root of $p$ on $[\sigma+2\eta L,\infty)$. Then $\rho(G(\eta))<\zeta_0$ will follow from $p(\zeta_0)>0$, which is shown below. 
{\small\begin{align*}
&p(\zeta_0) \\
&= \Big[4(\eta L)^{1/3}( 4(\eta L)^{1/3} + \frac{2\alpha }{1+\alpha}\sigma ) - k_1 \Big] (4(\eta L)^{1/3} - 2\eta L) - k_2\\
&> ( 16 (\eta L)^{2/3} - 5 \eta L )(4(\eta L)^{1/3}  -2\eta L ) - k_2\\
&>  ( 16 (\eta L)^{2/3} - 5 (\eta L)^{2/3} )(4(\eta L)^{1/3}  -2 (\eta L)^{1/3} ) - k_2 \\
&= 22\eta L - k_2 >0.
\end{align*}}For the lower bound, notice that $k_2^{1/3}>2 (\eta L)^{2/3} > 4 \eta L  $ (which is equivalent to $\eta L<\frac{1}{8}$, cf. the step size condition in Lemma~\ref{lem:G_eig}). Therefore, we have
{\small\begin{align*}
&p(\sigma + 0.84 k_2^{1/3} ) \\
& =  (0.84 k_2^{1/3} )(0.84 k_2^{1/3} + \frac{2\alpha}{1+\alpha} \sigma )(0.84 k_2^{1/3} - 2\eta L )\\
&\qquad  -k_1(0.84  k_2^{1/3} - 2\eta L ) - k_2\\
&<  (0.84  k_2^{1/3} )(0.84  k_2^{1/3} + \frac{1}{2^{5/6}} k_2^{1/3}  )(0.84  k_2^{1/3}  ) - k_2  <0
\end{align*}}where we have used that, noticing $k_2> 2\eta L \sigma$ and $\alpha = \sqrt{\mu \eta} < \sqrt{\eta L}< \sqrt{1/8^3}$ (cf. the step size condition in Lemma~\ref{lem:G_eig}), $[\frac{\alpha}{1+\alpha} \sigma]^3 < \alpha^3 \sigma <  \eta L \alpha  \sigma < \frac{\alpha}{2} k_2 < \frac{k_2}{2} \sqrt{1/ 8^3} \Rightarrow \frac{\alpha}{1+\alpha} \sigma < \frac{k_2^{1/3}}{2\cdot 2^{5/6}}.$
Hence we have $\rho(G(\eta)) > \sigma + 0.84 k_2^{1/3}$. Also noticing $k_2> 2\eta L\sigma$, we have $  \rho(G(\eta)) >  \sigma+ 0.84(2\eta L\sigma   )^{1/3} > \sigma+ (\sigma \eta L)^{1/3} $.

Before we proceed to (b) and (c), we prove the following claim. 

 \textbf{Claim: } $G(\eta)$'s spectral gap is at least $\rho(G(\eta)) - \sigma$. 

To prove the claim, we consider two cases. If $G(\eta)$ has three real eigenvalues, then notice that $p_0(\zeta)$ is nonpositive on $[\sigma,\sigma + 2\eta L]$, and $p_1(\zeta)$ is affine and positive on $[\sigma,\sigma + 2\eta L]$ (since $p_1(\sigma+2\eta L) = k_2>0$, $p_1(\sigma) = k_2 - 2\eta L k_1= k_2 - 2\eta^2 L^2 (4+\sigma)>0$), therefore $p(\zeta)$ has no real roots on $[\sigma,\sigma+2\eta L]$. Also notice $p(\zeta)$ has exactly one real root on $[\sigma+2\eta L, \infty)$ (the leading eigenvalue; here we have used the fact that $p(\zeta)$ is strictly convex on $[\sigma+2\eta L,\infty)$ and $p(\sigma+2\eta L)<0$), hence $p$'s other two real roots must be less than $\sigma$. We next show that the two other real roots must be nonnegative. Let the three eigenvalues be $\gamma_1$, $\gamma_2$ and $\gamma_3$ with $\gamma_1$ being the leading eigenvalue. Then,
$$\gamma_1 + \gamma_2 + \gamma_3 = Trace[G(\eta)] = (1-\alpha)\sigma + \frac{1+\alpha^2}{1+\alpha}\sigma + \sigma + 2\eta L.$$
Then, using $\alpha = \sqrt{\mu\eta}<\frac{1}{2}$ (cf. the step size condition in Lemma~\ref{lem:G_eig}), we have $1 - \alpha> \frac{1}{2}$ and $ \frac{1+\alpha^2}{1+\alpha} > \frac{1}{2}$, and hence $\gamma_2 + \gamma_3 > 2 \sigma - \gamma_1 > \sigma - 4(\eta L)^{1/3} >0$,
where in the last inequality, we have used $\eta L < \frac{\sigma^3}{64}$ (cf. the step size condition in Lemma~\ref{lem:G_eig}).  Next, we calculate
{\small
\begin{align*}
&\gamma_1\gamma_2\gamma_3 \\
& = - p(0) = \sigma^2(\sigma+2\eta L) \frac{1-\alpha}{1+\alpha} + k_2 - k_1 (\sigma+2\eta L) \\
&=\sigma^2(\sigma+2\eta L) \frac{1-\alpha}{1+\alpha}  - 2\eta L\sigma \frac{1+\alpha - 2\alpha^2 }{1+\alpha} - \sigma^2\eta L \frac{1-\alpha^2-2\alpha^3}{1+\alpha}\\
&= \sigma^2\eta L \frac{1 - 2\alpha + \alpha^2 + 2\alpha^3 }{1+\alpha} + \sigma^3 \frac{1-\alpha}{1+\alpha} - 2\eta L \sigma  \frac{1+\alpha - 2\alpha^2 }{1+\alpha}\\
&> \frac{1}{2} \sigma^3  - 2\eta L \sigma
\end{align*}}where in the last inequality, we have used $1 - 2\alpha>0$, $\frac{1-\alpha}{1+\alpha} >\frac{1}{2}$ ($\Leftarrow \eta < \frac{1}{9\mu} $, implied by our step size selection in Lemma~\ref{lem:G_eig}). At last, we use the fact that our step size bound in Lemma~\ref{lem:G_eig} implies $2\eta L \sigma < \frac{1}{2}\sigma^3 $ ($\iff \eta L < \frac{1}{4} \sigma^2$). Then, $\gamma_1 \gamma_2 \gamma_3>0$. Since $\gamma_1>0$, we have $\gamma_2\gamma_3>0$. We have already shown $\gamma_2+ \gamma_3>0$. Hence, both $\gamma_2$ and $\gamma_3$ are positive. Now that we have already shown that, $\gamma_2$ and $\gamma_3$ are positive reals less than $\sigma$. This implies the spectral gap is at least $\rho(G(\eta))-\sigma$.  
  
  On the other hand, if $G(\eta)$ has one real eigenvalue (the leading one) and two conjugate complex eigenvalues, then let the modulus of the two conjugate eigenvalues be $\tau$ and we have 
\begin{align*}
\tau^2\rho(G(\eta))& = - p(0) = \sigma^2(\sigma+2\eta L) \frac{1-\alpha}{1+\alpha} + k_2 - k_1 (\sigma+2\eta L)
\end{align*}
We calculate
$k_2 - k_1 (\sigma+2\eta L) = - 2\eta L\sigma \frac{1+\alpha - 2\alpha^2 }{1+\alpha} - \sigma^2\eta L \frac{1-\alpha^2-2\alpha^3}{1+\alpha}<0$,
where we have used $2\alpha^2 < 1$ and $\alpha^2 + 2\alpha^3<1 $ (since $\alpha<\sqrt{\eta L} < \sqrt{1/8^3}$). Therefore,
$\tau^2\rho(G(\eta)) < \sigma^2 (\sigma+2\eta L)$. Noticing that $\rho(G(\eta))> \sigma+2\eta L$, we have $\tau <\sigma$. Therefore, we can conclude that the spectral gap is at least $\rho(G(\eta)) - \sigma$. 

Now we proceed to prove (b) and (c).

(b) Same as before, we let the three eigenvalues of $G(\eta)$ be $\gamma_1$, $\gamma_2$ and $\gamma_3$ with $\gamma_1 = \rho(G(\eta))$ being the leading eigenvalue and $\gamma_1 > |\gamma_2|\geq |\gamma_3|$. We assume $\gamma_2\neq \gamma_3$.\footnote{There will be some values of $\eta$ and $\alpha$ for which $G(\eta)$ will have only two eigenvalues, one of which has multiplicity $2$. This case can be dealt with by taking the $\gamma_2\rightarrow\gamma_3$ limit and won't affect our result.} Then by the claim on spectral gap, we have $\gamma_1 - |\gamma_3| \geq \gamma_1 - |\gamma_2| \geq (\sigma \eta L)^{1/3}$. We will use the following lemma.
\begin{lemma} \label{lem:G23_aux}
	If a series $(B_t)_{t\geq 0}$ can be written as $B_t =\alpha_1\gamma_1^t + \alpha_2\gamma_2^t + \alpha_3\gamma_3^t$ for some $\alpha_1,\alpha_2,\alpha_3\in\mathbb{C}$, then we have
	$$\vert B_t\vert  \leq \frac{  (5\vert B_0 \vert+ 6\vert B_1\vert +3  \vert B_2\vert )  \gamma_1^t}{(\sigma\eta L)^{2/3} }.  $$
\end{lemma}

It is not hard to see that the proof of Lemma~\ref{lem:G23_aux} is straightforward and relies on the fact that $\gamma_1^t$ is asymptotically larger than $\gamma_2^t$ and $\gamma_3^t$. \ifthenelse{\boolean{appendixin}}{The proof is deferred to Appendix-\ref{subsec:G23}.}{Due to space limit, the proof is omitted and can be found in \cite[Appendix-I]{fullversion}.}

We know that through diagonalization, all the entries of $G(\eta)^t$ can be written as the form described in Lemma \ref{lem:G23_aux}. We know that $[G(\eta)^0]_{23} = 0$, $[G(\eta)^1]_{23} = 2\eta L$, $[G(\eta)^2]_{23} = 2\eta L(\sigma + 2\eta L+ \sigma \frac{1+\alpha^2}{1+\alpha}) +\eta L \sigma \frac{1-\alpha}{1+\alpha} < 9\eta L$. As a result,
$ [G(\eta)^t]_{23} \leq \frac{39\eta L}{(\sigma\eta L)^{2/3}} \gamma_1^t =  \frac{39 (\eta L)^{1/3}}{(\sigma )^{2/3}} \gamma_1^t.$
%

(c) Similarly to (b), we can bound $[G(\eta)^t]_{21}$ and $[G(\eta)^t]_{22}$. We calculate, $[G(\eta)^0]_{21}=0$,  $[G(\eta)^1]_{21}<1$, $[G(\eta)^2]_{21}< 4$. Hence
$[G(\eta)^t]_{21} < \frac{18}{(\sigma\eta L)^{2/3}} \gamma_1^t. $
Similarly, $[G(\eta)^0]_{22}=1$,  $[G(\eta)^1]_{22}<1$, $[G(\eta)^2]_{22}< 6$. Hence
$[G(\eta)^t]_{22} < \frac{29}{(\sigma\eta L)^{2/3}} \gamma_1^t. $ 
As a result, $\max([G(\eta)^t]_{21}, [G(\eta)^t]_{22}, [G(\eta)^t]_{23}  )<\frac{39}{(\sigma\eta L)^{2/3}} \gamma_1^t.$
\qedd

\subsection{Proof of the Intermediate Result (Lemma \ref{lem:nsc:convergence})}\label{subsec:nsc:convergence}

First, we prove (\ref{eq:nsc:phi_upperbound}). Notice (\ref{eq:nsc:phi_upperbound}) is true for $t=0$. Then, assume it's true for $t$. For $t+1$, we have 
\begin{align*}
\Phi_{t+1}(\omega)&= (1 - \alpha_t)\Phi_t(\omega) + \alpha_t(\hat{f}(t) + \langle g(t), \omega - \bar{y}(t)\rangle )\\
&\leq(1 - \alpha_t)\Phi_t(\omega) + \alpha_t f(\omega)  \\
&\leq f(\omega) + \lambda_{t+1}(\Phi_0(\omega) - f(\omega)).  
\end{align*}
where in the first inequality we have used (\ref{eq:nsc:inexact:noisy_oracle_1}) and in the second inequality we have used the induction assumption.

Next, we prove (\ref{eq:nsc:inexact:phi_quadratic}). It's clear that $\Phi_t$ is always a quadratic function. Notice that $ \nabla^2 \Phi_{0}(\omega) = \gamma_0 I$ and $\nabla^2 \Phi_{t+1}(\omega) = (1 - \alpha_t)\nabla^2 \Phi_t(\omega). $
We get $\nabla^2 \Phi_t(\omega) = \gamma_t I$ for all $t$, by the definition of $\gamma_t$. 

We next claim, $\Phi_t$, as a quadratic function, achieves minimum at $\bar{v}(t)$. We prove the claim by induction. Firstly, $\Phi_0$ achieves minimum at $\bar{v}(0)$. Assume $\Phi_t$ achieves minimum at $\bar{v}(t)$. Then, $\nabla \Phi_t(\omega) = \gamma_t( \omega - \bar{v}(t))$. Then by (\ref{eq:nsc:inexact:phi_recursive}),
{\small\begin{align*}
	\nabla \Phi_{t+1}(\bar{v}(t+1)) &= (1 - \alpha_t ) \gamma_t ( \bar{v}(t+1)- \bar{v}(t) ) + \alpha_t g(t)  \\
	&= [-(1-\alpha_{t})\gamma_t \frac{\eta_t}{\alpha_t} + \alpha_t ] g(t)
	=0
	\end{align*}}where the last equality follows from the fact that $\alpha_{t}^2 = \eta_t (1-\alpha_t)\gamma_t$, which can be proved recursively. It's true for $t=0$ by definition of $\gamma_0$. And note $\alpha_{t+1}^2 = \frac{\eta_{t+1}}{\eta_t}(1-\alpha_{t+1})\alpha_t^2 = (1-\alpha_{t+1}  ) \eta_{t+1} (1-\alpha_t)\gamma_t = \eta_{t+1} (1-\alpha_{t+1}) \gamma_{t+1}$.  Hence $\Phi_{t+1}$ achieves optimum at $\bar{v}(t+1)$, and hence the claim. 

We have proven $\Phi_t(\omega)$ is a quadratic function that achieves minimum at $\bar{v}(t)$ and satisfies $\nabla^2\Phi_t = \gamma_t I$. This implies $\Phi_t$ can be written in the form of (\ref{eq:nsc:inexact:phi_quadratic}) for some unique $\phi_t^*\in\R$. 

We next show $\phi_t^*$ satisfies (\ref{eq:nsc:inexact:phi}). 
Clearly, $\phi_0^* = f(\bar{x}(0))$. We now derive a recursive formula for $\phi_t^*$. By (\ref{eq:nsc:inexact:phi_recursive})
$\Phi_{t+1}(\bar{y}(t)) = (1 - \alpha_t)\Phi_t(\bar{y}(t)) + \alpha_t\hat{f}(t).$ Plugging in (\ref{eq:nsc:inexact:phi_quadratic}), we get 
{\small\begin{align*}
	& \phi_{t+1}^* + \frac{\gamma_{t+1}}{2} \Vert \bar{y}(t) - \bar{v}(t+1)\Vert^2 \\
	&= (1 - \alpha_t) [\phi_t^* + \frac{\gamma_t}{2} \Vert \bar{y}(t) - \bar{v}(t)\Vert^2] +  \alpha_t\hat{f}(t).
	\end{align*}}Plugging in $  \bar{v}(t+1) - \bar{y}(t) = (\bar{v}(t) - \bar{y}(t)) - \frac{\eta_t}{\alpha_t } g(t)  $, we get the desired equality (\ref{eq:nsc:inexact:phi}).
\subsection{Proof of Lemma \ref{prop:nsc:simp_eq}}\label{subsec:nsc:useful_eq}

The derivations are shown below.\begin{align*} 
\bar{y}(t+1) - \bar{y}(t) 
& = (1-\alpha_{t+1})(\bar{y}(t) - \eta_t g(t)) + \alpha_{t+1}(\bar{v}(t) - \frac{\eta_t}{\alpha_t} g(t)) - \bar{y}(t)\\
&= \alpha_{t+1}(\bar{v}(t) - \bar{y}(t)) - \eta_t \big[\frac{\alpha_{t+1}}{\alpha_t} + 1 - \alpha_{t+1} \big] g(t)
\end{align*}
\begin{align*}
\bar{v}(t+1) - \bar{y}(t+1) 
&= \bar{v}(t) - \frac{\eta_t}{\alpha_t} g(t)    - (1-\alpha_{t+1})(\bar{y}(t) - \eta_t g(t)) - \alpha_{t+1}(\bar{v}(t) - \frac{\eta_t}{\alpha_t} g(t))\\
& = (1-\alpha_{t+1}) (\bar{v}(t) - \bar{y}(t)) + \eta_t (1 - \alpha_{t+1})(1 - \frac{1}{\alpha_t}) g(t)
\end{align*}

\subsection{Derivation of (\ref{eq:nsc:rel_con:z_recursive}) in the Proof of Lemma \ref{lem:nsc:rel_con_err}}\label{subsec:nsc:linear_sys_ineq}
By (\ref{eq:nsc:update_vector_a}) and (\ref{eq:nsc:update_ave_a}), we have
\begin{align}
&\Vert x(t+1) - \one \bar{x}(t+1)\Vert  \nonumber	\\
& = \Vert [Wy(t) - \one \bar{y}(t) ]- \eta_t [s(t) - \one g(t)] \Vert \nonumber \\
&\leq \sigma \Vert y(t) - \one \bar{y}(t) \Vert + \eta_t \Vert s(t) - \one g(t)\Vert. \label{eq:nsc:rel_con:ineq_x}
\end{align}

By (\ref{eq:nsc:update_vector_b}) and (\ref{eq:nsc:update_ave_b}), we have
\begin{align}
&\Vert v(t+1) - \one \bar{v}(t+1) \Vert\nonumber \\
\ifthenelse{\boolean{isfullversion}}{&\leq \Vert  [Wv(t) - \one\bar{v}(t) ]  - \frac{\eta_t}{\alpha_t} [s(t) - \one g(t) ]\Vert\nonumber\\ }{ }  
&\leq \sigma \Vert v(t) - \one \bar{v}(t)\Vert +  \frac{\eta_t}{\alpha_t} \Vert s(t) - \one g(t)\Vert.\label{eq:nsc:rel_con:ineq_v}
\end{align}
\ifthenelse{\boolean{isfullversion}}{	Hence,
	\begin{align}
	&\alpha_{t+1}\Vert v(t+1) - \one \bar{v}(t+1) \Vert \nonumber\\
	 &\leq \sigma \alpha_t \Vert v(t) - \one \bar{v}(t)\Vert +  \eta_t \Vert s(t) - \one g(t)\Vert.\label{eq:nsc:rel_con:ineq_v_2}
	\end{align}}{}
By (\ref{eq:nsc:update_vector_c}) and (\ref{eq:nsc:update_ave_c}), we have
{\small\begin{align}
&	\Vert y(t+1) - \one \bar{y}(t+1)\Vert \nonumber\\
&\leq  (1-\alpha_{t+1})\Vert x(t+1) - \one\bar{x}(t+1)\Vert  \nonumber\\
&\qquad  +\alpha_{t+1} \Vert v(t+1) - \one \bar{v}(t+1)\Vert\nonumber\\
\ifthenelse{\boolean{isfullversion}}{&\leq  (1-\alpha_{t+1}) \bigg[ \sigma \Vert y(t) - \one \bar{y}(t) \Vert + \eta_t \Vert s(t) - \one g(t)\Vert  \bigg]\nonumber\\
	&\qquad +  \alpha_{t+1}  \bigg[\sigma \Vert v(t) - \one \bar{v}(t)\Vert +  \frac{\eta_t}{\alpha_t} \Vert s(t) - \one g(t)\Vert \bigg] \nonumber\\}{}
&\leq \alpha_{t} \sigma \Vert v(t) -\one \bar{v}(t)\Vert  +  \sigma \Vert y(t) - \one\bar{y}(t)\Vert +  2\eta_t \Vert s(t) - \one g(t)\Vert \label{eq:nsc:rel_con:ineq_y}
\end{align}}where we have used (\ref{eq:nsc:rel_con:ineq_x}) and (\ref{eq:nsc:rel_con:ineq_v}) in the second inequality.

By (\ref{eq:nsc:update_vector_d}) and (\ref{eq:nsc:update_ave_d}), we have
{\small\begin{align}
&\Vert s(t+1) - \one g(t+1)\Vert \nonumber \\
&= \Vert W s(t) - \one g(t)   + [\nabla(t+1) - \nabla(t) - \one (g(t+1) - g(t))]\Vert \nonumber \\
&\stackrel{(a)} {\leq}\sigma \Vert s(t) - \one g(t)\Vert + \Vert \nabla (t+1) - \nabla(t)\Vert \nonumber\\
&\stackrel{(b)}{\leq} \sigma \Vert s(t) - \one g(t)\Vert + L \Vert y (t+1) - y(t)\Vert  \label{eq:nsc:rel_con:ineq_s_1} 
\end{align}}
where in (a) we have used the fact that
	\begin{align*}
	&\Big\Vert[\nabla(t+1)-\nabla(t)]  -   [\one g(t+1) - \one g(t)]\Big\Vert^2 \\
	&= \Vert \nabla(t+1)-\nabla(t) \Vert^2 - n \Vert g(t+1) - g(t)\Vert ^2\\
	&\leq \Vert \nabla(t+1)-\nabla(t) \Vert^2
	\end{align*}
	and in (b) we have used (\ref{eq:nes:ineq_useful_a}).

Now we expand $y(t+1) - y(t)$,
\begin{align}
&\Vert y(t+1) - y(t)\Vert \nonumber\\
& \leq \Vert y(t+1) - \one \bar{y}(t+1)\Vert + \Vert y(t) - \one \bar{y}(t)\Vert  \nonumber\\
&\qquad + \Vert \one \bar{y}(t+1) - \one \bar{y}(t)\Vert\nonumber\\
&\leq \alpha_t  \Vert v(t) -\one \bar{v}(t)\Vert + 2 \Vert y(t) - \one\bar{y}(t)\Vert  \nonumber\\
&\qquad + 2\eta_t \Vert s(t) - \one g(t)\Vert + \Vert \one \bar{y}(t+1) - \one \bar{y}(t)\Vert.
\end{align}

Combining the above with (\ref{eq:nsc:rel_con:ineq_s_1}), we get,
{\small\begin{align}
&\Vert s(t+1) - \one g(t+1)\Vert \nonumber\\
&\leq 	\alpha_{t}  L \Vert v(t) -\one \bar{v}(t)\Vert  + 2 L\Vert y(t) - \one\bar{y}(t)\Vert \nonumber\\
&\quad +  \bigg[\sigma + 2 \eta_t L  \bigg] \Vert s(t) - \one g(t)\Vert + L\sqrt{n} \Vert \bar{y}(t+1) -  \bar{y}(t)\Vert.    \label{eq:nsc:rel_con:ineq_s}
\end{align}}

By (\ref{eq:nsc:y_y}) we have $a(t)\geq L\Vert\bar{y}(t+1)-\bar{y}(t)\Vert$. Now we combine \ifthenelse{\boolean{isfullversion}}{(\ref{eq:nsc:rel_con:ineq_v_2})}{(\ref{eq:nsc:rel_con:ineq_v})} (\ref{eq:nsc:rel_con:ineq_y})  (\ref{eq:nsc:rel_con:ineq_s}) and get the desired linear system inequality (\ref{eq:nsc:rel_con:z_recursive}).
\qedd

\subsection{Proof of Lemma \ref{lem:nsc:G_eig}, Lemma \ref{lem:nsc:G_eig_lowerboud} and Lemma \ref{lem:nsc:G_chi}}\label{subsec:nsc:G}
	We first prove Lemma \ref{lem:nsc:G_eig}.
	
	\noindent\textit{Proof of Lemma \ref{lem:nsc:G_eig}: }
	We write down the charasteristic polynomial of $G(\eta)$ as
	\begin{align}
	p(\zeta) & = (\zeta-\sigma)^2(\zeta - \sigma - 2\eta L)  - 5\eta L (\zeta- \sigma-2\eta L) \nonumber\\
	&\qquad -2\eta L \sigma-10\eta^2 L^2.\label{eq:G_charpoly}
	\end{align}

We evaluate $p(\cdot)$ on $\sigma + 4(\eta L)^{1/3}$ and get,
	\begin{align*}
&	p(\sigma + 4(\eta L)^{1/3})\\
	  &= (16 (\eta L)^{2/3} - 5\eta L)(4(\eta L)^{1/3} - 2\eta L) -2\eta L \sigma-10\eta^2 L^2\\
	&> 22\eta L - 2\eta L\sigma -10  \eta^2 L^2>0
	\end{align*}
	where we have used $(\eta L)^2<(\eta L) < (\eta L)^{2/3} <(\eta L)^{1/3}$. 
	It's easy to check $p(\sigma)<0$. Also, $p'(\zeta) = 2(\zeta-\sigma)(\zeta - \sigma - 2\eta L) + (\zeta-\sigma)^2 - 5\eta L>0$ on $[\sigma+4(\eta L)^{1/3},\infty)$, so $p$'s largest real root must lie within $(\sigma,\sigma+4(\eta L)^{1/3})$.  Therefore, $\sigma<\theta(\eta) <\sigma + 4(\eta L)^{1/3}$.
	For the eigenvectors, notice that
	$$ L \chi_1(\eta) + 2L \chi_2(\eta) + (\sigma+2\eta L)\chi_3(\eta) =\theta(\eta)\chi_3(\eta).$$ Hence, $\theta(\eta)\chi_3(\eta) \geq \sigma \chi_3(\eta) + 2L\chi_2(\eta)$, and $\chi_2(\eta)/\chi_3(\eta) \leq \frac{\theta(\eta)-\sigma}{2 L} \leq \frac{2}{L^{2/3}} \eta^{1/3}$.	\qedd
	
For the rest of this section, we will use the following formula for $\chi_1(\cdot)$ and $\chi_2(\cdot)$.
\begin{lemma}\label{lem:nsc:G_chi_formula}
	We have
{\small	\begin{align*}
	\chi_1(\eta)& = \frac{\eta}{\theta(\eta) - \sigma},
	\chi_2(\eta)= \frac{\theta(\eta) - \sigma}{2L} - \eta - \frac{\eta}{2(\theta(\eta)- \sigma)}.
	\end{align*}}
\end{lemma}
\noindent\textit{Proof: } Since $G(\eta)\chi(\eta) = \theta(\eta)\chi(\eta)$, we have, writing down the first line,
$\sigma \chi_1(\eta) + \eta \chi_3(\eta) = \theta(\eta)\chi_1(\eta) $,
from which we can get the formula for $\chi_1(\eta)$. Writing the third line of $G(\eta)\chi(\eta) = \theta(\eta)\chi(\eta)$, we get
$$L \chi_1(\eta) + 2L\chi_2(\eta) + (\sigma+2\eta L) \chi_3(\eta) = \theta(\eta)\chi_3(\eta)$$
from which we can derive the formula for $\chi_2(\eta)$.\qedd

Now we proceed to prove Lemma \ref{lem:nsc:G_eig_lowerboud}.

\noindent\textit{Proof of Lemma~\ref{lem:nsc:G_eig_lowerboud}: } Since $\eta L < \frac{\sqrt{\sigma}}{\sqrt{8}}$, we have $(\sigma \eta L)^{1/3}  > 2\eta L$. Hence by (\ref{eq:G_charpoly}), $p(\sigma +  (\sigma \eta L)^{1/3}) \leq (\sigma\eta L)^{2/3}( (\sigma\eta L)^{1/3} - 2\eta L) - 2\sigma \eta L < 0$. Hence $\theta(\eta)> \sigma+(\sigma\eta L)^{1/3}$. As a result, $\chi_1(\eta) < \frac{\eta}{(\sigma\eta L)^{1/3}}$.\qedd

The rest of the section will be devoted to Lemma~\ref{lem:nsc:G_chi}. Before proving Lemma~\ref{lem:nsc:G_chi}, we prove an axillary lemma first.
\begin{lemma}\label{lem:nsc:G_thetad}
	We have 
	$$\theta'(\eta) =\frac{2L (\theta(\eta) - \sigma)^2 + 5 L (\theta(\eta) - \sigma) + 2\sigma L}{3(\theta(\eta) - \sigma)^2 - 4\eta L(\theta(\eta)-\sigma) - 5\eta L  }. $$
	Further, when $0<\eta < \frac{\sigma^2}{L9^3} $, we have $0<\theta'(\eta)  < 5 \frac{L^{1/3}}{\eta^{2/3}}.$
\end{lemma}
\noindent\textit{Proof: } Since $p(\theta(\eta)) = 0$ (where $p$ is the characteristic polynomial of $G(\eta)$ as defined in (\ref{eq:G_charpoly}), we take derivative w.r.t. $\eta$ on both sides of $p(\theta(\eta)) = 0$, and get
\begin{align*}
0 =& 2(\theta(\eta) - \sigma) \theta'(\eta)(\theta(\eta) - \sigma - 2\eta L) \\
& + (\theta(\eta) - \sigma)^2(\theta'(\eta) - 2 L) - 5L(\theta(\eta) - \sigma - 2\eta L)\\
& - 5\eta L (\theta'(\eta) - 2L) - 2\sigma L - 20\eta L^2.
\end{align*}
From the above equation we obtain 
\begin{equation}
\theta'(\eta) = \frac{2L (\theta(\eta) - \sigma)^2 + 5 L (\theta(\eta) - \sigma) + 2\sigma L}{3(\theta(\eta) - \sigma)^2 - 4\eta L(\theta(\eta)-\sigma) - 5\eta L  }.
\end{equation}
By Lemma \ref{lem:nsc:G_eig}, the step size $\eta < \frac{\sigma^2}{L9^3} $ implies that 
$\theta(\eta) - \sigma<4(\eta L)^{1/3} < \sigma^{2/3}$.
Hence, the nominator part of $\theta'(\eta)$ satisfy, 
\begin{align*}
2L (\theta(\eta) - \sigma)^2 + 5 L (\theta(\eta) - \sigma) + 2\sigma L< 9L\sigma^{2/3}.
\end{align*}
Also notice
$4\eta L(\theta(\eta)-\sigma) + 5\eta L< 9\eta L < (\sigma\eta L)^{2/3}$
where the last inequality comes from the step size condition $\eta < \frac{\sigma^2}{L9^3} $. The step size condition also implies $\eta < \frac{\sqrt{\sigma}}{L 2\sqrt{2}}$ and hence we can use Lemma \ref{lem:nsc:G_eig_lowerboud} to get $(\theta(\eta) - \sigma)^2 \geq (\sigma \eta L)^{2/3}$. Hence, the denominator part of $\theta'(\eta)$ satisfy
{\small\begin{align*}
&3(\theta(\eta) - \sigma)^2 - 4\eta L(\theta(\eta)-\sigma) - 5\eta L \\
& \geq 3 (\sigma \eta L)^{2/3} - (\sigma\eta L)^{2/3} = 2(\sigma\eta L)^{2/3}.
\end{align*}}Combining the bound on the nominator and the denominator, we get $0<\theta'(\eta) \leq \frac{9 L \sigma^{2/3}}{2(\sigma \eta L)^{2/3}}< 5 \frac{L^{1/3}}{\eta^{2/3}}.$
\qedd

Now we proceed to prove Lemma \ref{lem:nsc:G_chi}. 

\noindent\textit{Proof of Lemma \ref{lem:nsc:G_chi}: } Define $\xi_1,\xi_2:\R\rightarrow\R$ to be $\xi_1(y) = \log \chi_1(e^y)$ and $\xi_2(y) = \log \chi_2(e^y)$. The statement of the lemma can be rephrased as,
{\small\begin{align*}
\log \chi_1( e^{\log \zeta_1 } ) - \log \chi_1( e^{\log \zeta_2 } ) &\leq \frac{6}{\sigma} | \log \zeta_1 - \log \zeta_2|,\\
\log \chi_2( e^{\log \zeta_1 } ) - \log \chi_2( e^{\log \zeta_2 } ) &\leq \frac{28}{\sigma} | \log \zeta_1 - \log \zeta_2|.
\end{align*}}So we only need to show that $\xi_1$ is $\frac{6}{\sigma}$-Lipschitz continuous and $\xi_2$ is $\frac{28}{\sigma}$-Lipschitz continuous, on interval $y\in(-\infty, \log( \frac{\sigma^2}{9^3 L}))$. 
We calculate the derivative of $\xi_1$ and let $\eta = e^y\in (0, \sigma^2/(9^3 L))$,  
{\small\begin{align*}
\xi_1'(y) &= \frac{\chi_1' (\eta) \eta}{\chi_1(\eta)}  = \frac{ \eta}{\chi_1(\eta)} \big[ \frac{(\theta(\eta) - \sigma) - \eta \theta'(\eta)}{  (\theta(\eta) - \sigma)^2 }\big]
=1 - \frac{\eta\theta'(\eta)}{\theta(\eta) - \sigma}.
\end{align*}}
Notice that $\eta$ satisfies all the step size conditions of Lemma~\ref{lem:nsc:G_eig},\ref{lem:nsc:G_eig_lowerboud},\ref{lem:nsc:G_thetad}. We have
\begin{align*}
|\xi_1'(y)| &\leq1  + \frac{\eta\theta'(\eta)}{\theta(\eta) - \sigma}
\leq 1 + \frac{\eta  \cdot 5 L^{1/3}/\eta^{2/3} }{(\sigma\eta L)^{1/3}}\\
&\leq 6/\sigma^{1/3} < 6/\sigma.
\end{align*}
This implies $\xi_1$ is $6/\sigma$-Lipschitz continuous.
Similarly, for $\xi_2$, we have
$\xi_2'(y) = \frac{\chi_2' (\eta) \eta}{\chi_2(\eta)}.$
By Lemma \ref{lem:nsc:G_chi_formula}, we have
{\small\begin{align*}
\chi_2(\eta) & \geq \frac{(\sigma \eta L)^{1/3}}{2L} - \eta - \frac{\eta}{2 (\sigma\eta L)^{1/3}} \\
& \stackrel{(a)}{\geq} \frac{(\sigma \eta L)^{1/3}}{2L} -2 \frac{\eta}{(\sigma\eta L)^{1/3}}
 \stackrel{(b)}{\geq} \frac{(\sigma \eta L)^{1/3}}{4L}
\end{align*}}where in (a), we have used $(\sigma \eta L)^{1/3} < \sigma/9 < 1$, and in (b), we have used $\frac{(\sigma \eta L)^{1/3}}{4L} \geq 2 \frac{\eta}{(\sigma\eta L)^{1/3}}$, which is equivalent to $\eta L \leq \sigma^2/8^3$ and follows from our step size condition. Now, we calculate $\chi_2'(\eta)$, 
{\small\begin{align*}
|\chi_2'(\eta)| &= | \frac{\theta'(\eta)}{2L} - 1 - \frac{(\theta(\eta)- \sigma) - \eta\theta'(\eta)}{2 (\theta(\eta)- \sigma)^2} |\\
&\leq \frac{\theta'(\eta)}{2L} + 1 + \frac{1}{2 (\theta(\eta)-\sigma)} + \frac{\eta \theta'(\eta)}{2 (\theta(\eta)-\sigma)^2}\\
&\leq \frac{5}{2} \frac{1}{(\eta L)^{2/3}} + 1 + \frac{1}{2 (\sigma \eta L)^{1/3}} + \frac{1}{2} \frac{\eta 5 L^{1/3} / \eta^{2/3}}{(\sigma \eta L)^{2/3}}\\
&= \frac{5}{2} \frac{1}{(\eta L)^{2/3}} + 1 + \frac{1}{2 (\sigma \eta L)^{1/3}} + \frac{5}{2} \frac{1}{\sigma^{2/3} (\eta L)^{1/3}}\\
&\leq  \frac{7}{(\sigma\eta L)^{2/3} }.
\end{align*}}

Therefore, $
|\xi_2'(y)| \leq \frac{ \eta \cdot 7/ (\sigma \eta L)^{2/3}  }{ (\sigma \eta L)^{1/3} / (4L)} = \frac{28}{\sigma}.$
Hence, $\xi_2$ is $\frac{28}{\sigma}$-Lipschitz continuous.\qedd

\qedd

\qedd


\subsection{Proof of Lemma \ref{lem:nsc:alpha}}\label{subsec:nsc:alpha}
We provide a comparison lemma that will be helpful for the analysis.
\begin{lemma}\label{lem:nsc:comparison}
	Given two step size sequences $(\eta_t)_t$ and $(\eta_t')_t$ that satisfy, $\eta_0 \leq \eta_0'$, and $\forall t$, $\frac{\eta_{t+1}}{\eta_{t}} \leq \frac{\eta_{t+1}'}{\eta_{t}'}$, then $\alpha_t\leq \alpha_t'$.
\end{lemma} 
\noindent\textit{Proof: } We prove the statement by induction. First, $\alpha_0 = \sqrt{\eta_0 L}\leq \sqrt{\eta_0' L} = \alpha_0'$. Next, assume $\alpha_t \leq\alpha_t'$, then $\frac{\eta_{t+1}}{\eta_{t}} \alpha_t^2 \leq \frac{\eta_{t+1}'}{\eta_{t}'} \alpha_t'^2$. Define function $\xi: (0,1) \rightarrow (0,\infty)$ with $\xi(y) = y^2/(1-y)$. It's easy to check that $\xi$ is a strictly increasing function and is a bijiection. Notice that $\alpha_{t+1} = \xi^{-1} ( \frac{\eta_{t+1}}{\eta_{t}} \alpha_t^2  )\leq \xi^{-1} ( \frac{\eta_{t+1}'}{\eta_{t}'} \alpha_t'^2) = \alpha_{t+1}'$. So we are done. 
\qedd 

\noindent\textit{Proof of Lemma \ref{lem:nsc:alpha}:}
Since $\eta_{t+1}<\eta_t$, we have $\alpha_{t+1}^2<\alpha_t^2$ and hence $\alpha_t$ is decreasing. Now we derive the asymptotic convergence rate of $\alpha_t$. 

\noindent\textbf{Proof of (i).} Firstly, we have,
$$ \alpha_{t+1}^2 + \frac{\eta_{t+1}}{\eta_t}\alpha_t^2 \alpha_{t+1}  -  \frac{\eta_{t+1}}{\eta_t}\alpha_t^2 = 0$$
Hence
\begin{align*}
\alpha_{t+1} &= \frac{1}{2}\bigg(-  \frac{\eta_{t+1}}{\eta_t}\alpha_t^2 + \sqrt{\frac{\eta_{t+1}^2}{\eta_t^2}\alpha_t^4  + 4 \frac{\eta_{t+1}}{\eta_t}\alpha_t^2 } \bigg)\\
&=   \frac{2 \frac{\eta_{t+1}}{\eta_t} \alpha_t^2 }{ \frac{\eta_{t+1}}{\eta_t}\alpha_t^2 + \sqrt{\frac{\eta_{t+1}^2}{\eta_t^2}\alpha_t^4  + 4 \frac{\eta_{t+1}}{\eta_t}\alpha_t^2 }}  \\
&=  \frac{2   }{ 1 + \sqrt{1  + 4 \frac{\eta_{t}}{\eta_{t+1} \alpha_t^2} }}.  
\end{align*}
Hence
\begin{align}
\frac{1}{\alpha_{t+1}} - \sqrt{\frac{\eta_t}{\eta_{t+1}}} \frac{1}{\alpha_t} &= \frac{1}{2} + \frac{1}{2}\sqrt{1  + 4 \frac{\eta_{t}}{\eta_{t+1} \alpha_t^2} } - \sqrt{\frac{\eta_t}{\eta_{t+1}}} \frac{1}{\alpha_t}\nonumber\\
&=\frac{1}{2} + \frac{1}{4} \frac{1}{\frac{1}{2}\sqrt{1  + 4 \frac{\eta_{t}}{\eta_{t+1} \alpha_t^2} } + \sqrt{\frac{\eta_t}{\eta_{t+1}}} \frac{1}{\alpha_t}}\label{eq:alpha:onehalf}\\
&> \frac{1}{2}.\nonumber
\end{align}
Hence,  $ \frac{\sqrt{\eta_{t}}}{\alpha_{t}} - \frac{\sqrt{\eta_{t-1}}}{\alpha_{t-1}} > \frac{1}{2}\sqrt{\eta_t}$.
Therefore, for $t\geq 1$,
\begin{align}
\frac{\sqrt{\eta_{t}}}{\alpha_{t}} &> \frac{\sqrt{\eta_{0}}}{\alpha_{0}} + \frac{1}{2} \sum_{k=1}^t \sqrt{\eta_k}.\label{eq:alpha:alpha_sum}
\end{align}
Now we use the comparison lemma (Lemma \ref{lem:nsc:comparison}) with a fixed step size $\eta_t' = \eta_0$. Notice that $\eta_{t+1}/\eta_t < \eta_{t+1}' /\eta_t'$, we have $\alpha_t\leq \alpha_t'$. Repeating the argument for (\ref{eq:alpha:alpha_sum}), we get (\ref{eq:alpha:alpha_sum}) is also true if we replace $\eta_t$ with $\eta_t'$ and $\alpha_t$ with $\alpha_t'$. Hence, for $t\geq 1$, 
$$\frac{\sqrt{\eta_{0}}}{\alpha_{t}'} > \frac{\sqrt{\eta_{0}}}{\alpha_{0}'} + \frac{1}{2} t \sqrt{\eta_0} > \frac{1}{2} (t+1) \sqrt{\eta_0}.  $$
This implies that for $t\geq 0$, $\alpha_t' \leq \frac{2}{t+1}$. Hence,
\begin{align}
\alpha_t \leq \alpha_t' \leq \frac{2}{t +1}. \label{eq:alpha:ub}
\end{align}
This gives part (i) of the lemma. Now we derive a tighter upper bound for $\alpha_t$ which will be used later. Returning to (\ref{eq:alpha:alpha_sum}), we have
\begin{align*}
\frac{\sqrt{\eta_{t}}}{\alpha_{t}} &> \frac{\sqrt{\eta_{0}}}{\alpha_{0}} + \frac{1}{2} \sum_{k=1}^t \sqrt{\eta_k}\\
&> \frac{1}{2} \sum_{k=1}^t \frac{\sqrt{\eta}}{  (k+t_0)^{\beta/2}} \\
&> \frac{1}{2} \sqrt{\eta} \int_1^{t+1} \frac{1}{(y+t_0)^{\beta/2}} dy\\
&= \frac{1}{2} \sqrt{\eta} \frac{2}{2-\beta} [(t+1+t_0)^{1-\frac{\beta}{2}} - (t_0+1)^{1-\frac{\beta}{2}} ].
\end{align*}
Therefore,
\begin{align*}
\frac{1}{\alpha_t} &> \frac{1}{2} \sqrt{\eta} \frac{2}{2-\beta} [(t+1+t_0)^{1-\frac{\beta}{2}} - (t_0+1)^{1-\frac{\beta}{2}} ] \times \frac{ (t+t_0)^{\beta/2}}{\sqrt{\eta}}\\
&> \frac{1}{2-\beta}(t + t_0 -  (t+t_0)^{\beta/2} (t_0+1)^{1- \frac{\beta}{2}}).
\end{align*}
Notice that when $t\geq 2$, the right hand side of the above formula is positive. Hence, 
\begin{align}
\alpha_t < \frac{2-\beta}{t + t_0 -  (t+t_0)^{\beta/2} (t_0+1)^{1- \frac{\beta}{2}} }, \forall t\geq 2.\label{eq:alpha:ub_tight}
\end{align}

\noindent\textbf{Proof of (ii).} We return to (\ref{eq:alpha:onehalf}), and use (\ref{eq:alpha:ub}) to get,
\begin{align*}
\frac{1}{\alpha_{t+1}} - \sqrt{\frac{\eta_t}{\eta_{t+1}}} \frac{1}{\alpha_t} &< \frac{1}{2} + \frac{1}{8} \alpha_t \sqrt{\frac{\eta_{t+1}}{\eta_t}} \leq \frac{1}{2} + \frac{1}{4(t+1)}.
\end{align*}
Hence, for $t\geq 1$,
$$ \frac{\sqrt{\eta_{t}}}{\alpha_{t}} - \frac{\sqrt{\eta_{t-1}}}{\alpha_{t-1}} \leq \frac{1}{2}\sqrt{\eta_t}  + \frac{\sqrt{\eta_t}}{4t}.$$
Therefore, for $t\geq 1$,
\begin{align*}
\frac{\sqrt{\eta_{t}}}{\alpha_{t}}
&\leq \frac{\sqrt{\eta_{0}}}{\alpha_{0}} + \frac{1}{2} \sum_{k=1}^t\sqrt{\eta_k} + \sum_{k=1}^{t} \frac{\sqrt{\eta_k}}{4k} \\
&< \frac{\sqrt{\eta_0}}{\alpha_{0}} + \frac{1}{2} \sum_{k=1}^t\sqrt{\eta}\frac{1}{(k+t_0)^{\beta/2}} + \sum_{k=1}^{t} \sqrt{\eta}\frac{1}{4k ^{\beta/2+1}} \\
&\leq \frac{\sqrt{\eta_0}}{\alpha_{0}}  +  \frac{\sqrt{\eta}}{2}\int_0^t \frac{1}{(y+t_0)^{\beta/2}} dy  +\frac{\sqrt{\eta}}{4} [1 + \int_1^t \frac{1}{y^{\beta/2+1}} dy]\\
&=\frac{\sqrt{\eta_0}}{\alpha_{0}} +  \frac{\sqrt{\eta}}{2} \frac{2}{2-\beta}[(t+t_0)^{1-\beta/2} - t_0^{1-\beta/2} ] + \frac{\sqrt{\eta}}{4} [1 +\frac{2}{\beta} - \frac{2}{\beta}\frac{1}{t^{\beta/2}}  ] \\
&\leq   \frac{\sqrt{\eta_0}}{\alpha_{0}}  + \frac{\sqrt{\eta}}{2-\beta}(t+t_0)^{1-\beta/2} + \frac{\sqrt{\eta}}{4} [1 +\frac{2}{\beta} ] .
\end{align*}
Therefore, for $t\geq 1$,
\begin{align*}
\frac{1 }{\alpha_{t}} \leq \frac{1}{2-\beta} (t+t_0)+ (\frac{1}{\alpha_0 t_0^{\beta/2}} + \frac{1}{2\beta} + \frac{1}{4}) (t+t_0)^{\beta/2}.
\end{align*}
And then, for $t\geq 1$,
\begin{align*}
\frac{2-\beta}{t+t_0} - \alpha_t 
&\leq \frac{2-\beta}{t+t_0} - \frac{2-\beta}{ (t+t_0)+ (2-\beta)(\frac{1}{\alpha_0 t_0^{\beta/2}} + \frac{1}{2\beta} + \frac{1}{4}) (t+t_0)^{\beta/2}  }\\
&=  \frac{(2-\beta)^2 (\frac{1}{\alpha_0 t_0^{\beta/2}} + \frac{1}{2\beta} + \frac{1}{4}) }{(t+t_0) [ (t+t_0)^{1-\frac{\beta}{2}} +  (2-\beta) (\frac{1}{\alpha_0 t_0^{\beta/2}} + \frac{1}{2\beta} + \frac{1}{4})  ] }\\
& = O( \frac{1}{(t+t_0)^{2 - \frac{\beta}{2}} }).
\end{align*}
Now we consider $\lambda_t = \prod_{k=0}^{t-1} (1-\alpha_k)$. We have 
\begin{align*}
\log \lambda_t &\leq - \sum_{k=0}^{t-1} \alpha_k\\
&=   -\alpha_0 + \sum_{k=1}^{t-1} (  \frac{2-\beta}{k+t_0} - \alpha_k  ) - \sum_{k=1}^{t-1} \frac{2-\beta}{k+t_0} \\
&\leq -\alpha_0 + O(\sum_{k=1}^\infty \frac{1}{(k+t_0)^{2 - \frac{\beta}{2}} }) - (2-\beta)\log\frac{t+t_0}{1+t_0}\\
& = - (2-\beta) \log (t+t_0) + O(1)
\end{align*}
where we have used $\sum_{k=1}^\infty \frac{1}{(k+t_0)^{2 - \frac{\beta}{2}} } < \infty$ since $2-\frac{\beta}{2}>1$. 
Hence, $\lambda_t = O(\frac{1}{(t+t_0)^{2-\beta}}) = O(\frac{1}{t^{2-\beta}})$, i.e. the part (ii) of this lemma.

\noindent\textbf{Proof of (iii).} It is easy to check that $\forall y\in(-1,\infty)$, $\log (1+y) \geq \frac{y}{1+y}$. Therefore,
\begin{align}
\log \lambda_t \geq -\sum_{k=0}^{t-1} \frac{\alpha_k}{1-\alpha_k} = -\sum_{k=0}^{t-1} \alpha_k  - \sum_{k=0}^{t-1} \frac{\alpha_k^2}{1-\alpha_k}. \label{eq:lambda_lb:sum}
\end{align}
By (\ref{eq:alpha:ub}) and the fact $\alpha_k\leq \alpha_0<\frac{1}{2}$, we have
\begin{align}
\sum_{k=0}^{t-1} \frac{\alpha_k^2}{1-\alpha_k} &\leq 2 \sum_{k=0}^{t-1} \alpha_k^2   \leq 8 \sum_{k=0}^{t-1} \frac{1}{(k+1)^2}  \leq 8 \cdot \frac{\pi^2}{6} < 14\label{eq:lambda_lb:term2}.
\end{align}
We next bound $\sum_{k=0}^{t-1} \alpha_k$. Notice that when $k\geq t_0+2$, we have 
$$ (k+t_0) \geq 2^{1-\frac{\beta}{2}} (k+t_0)^{\beta/2} (1+t_0)^{1-\frac{\beta}{2}}.$$
Hence, by (\ref{eq:alpha:ub_tight}), we have when $k\geq t_0+2$,
\begin{align*}
\alpha_k &<  \frac{2-\beta}{k+ t_0 -  (k+t_0)^{\beta/2} (t_0+1)^{1- \frac{\beta}{2}} } - \frac{2-\beta}{k+t_0} + \frac{2-\beta}{k+t_0}\\
&=\frac{ (2-\beta) (k+t_0)^{\beta/2} (t_0+1)^{1- \frac{\beta}{2}} }{(k+t_0) (k+ t_0 -  (k+t_0)^{\beta/2} (t_0+1)^{1- \frac{\beta}{2}} )}+ \frac{2-\beta}{k+t_0}\\
&\leq \frac{(2-\beta)(t_0+1)^{1- \frac{\beta}{2}} }{ (1-2^{\frac{\beta}{2} -1}) (k+t_0)^{2-\frac{\beta}{2}}} + \frac{2-\beta}{k+t_0}.
\end{align*}
Hence, for $t>t_0+3$,
\begin{align*}
&\sum_{k=t_0+3}^{t-1} \alpha_k \leq \sum_{k=t_0+3}^{t-1}  \frac{(2-\beta)(t_0+1)^{1- \frac{\beta}{2}} }{ (1-2^{\frac{\beta}{2} -1}) (k+t_0)^{2-\frac{\beta}{2}}}  + \sum_{k=t_0+3}^{t-1}  \frac{2-\beta}{k+t_0}\\
&\leq \frac{2(t_0+1)^{1- \frac{\beta}{2}} }{ (1-2^{\frac{\beta}{2} -1}) (2t_0+2)^{1-\frac{\beta}{2}}} + (2-\beta)\log\frac{t-1+t_0}{2t_0+2}\\
&=\frac{2}{ 2^{1- \frac{\beta}{2} }-1 } + (2-\beta)\log\frac{t-1+t_0}{2t_0+2}\\
&<  \frac{4}{(2-\beta)\log 2} + (2-\beta)\log (t-1+t_0).
\end{align*}
Notice that, when $t\leq t_0+3$, the above inequality is still true since the left hand side is $0$ while the right hand side is positive.
Hence,
\begin{align*}
\sum_{k=0}^{t-1} \alpha_k 
&\leq \sum_{k=0}^{t_0+2} \alpha_k + \sum_{k=t_0+3}^{t-1} \alpha_k \\
&\leq  \sum_{k=0}^{t_0+2} \frac{2}{k+1} + \sum_{k=t_0+3}^{t-1} \alpha_k \\
&\leq 2+ 2\log(t_0+3) +  \frac{4}{(2-\beta)\log 2}   + (2-\beta)\log(t-1+t_0)\\
&\leq (2-\beta)\log(t-1+t_0) + 2\log(t_0+3) +2+\frac{6}{2-\beta}.  
\end{align*}
Therefore, combining the above with (\ref{eq:lambda_lb:term2}) and (\ref{eq:lambda_lb:sum}), we get
\begin{align*}
\log\lambda_t& \geq -(2-\beta)\log(t-1+t_0)  - 2\log(t_0+3) - \frac{6}{2-\beta}  -16.
\end{align*}
Hence, 
\begin{align*}
\lambda_t \geq \frac{1}{(t+t_0)^{2-\beta}  (t_0+3)^2 e^{16+ \frac{6}{2-\beta} }}.
\end{align*}
\qedd

\subsection{Proof of Theorem~\ref{thm:nsc:vanishing} (b)} \label{appendix:theorem4b}
Finally, we show part (b) of the Theorem, i.e. upper bounding individual error $f(y_i(t)) - f^*$. We continue inequality \eqref{eq:nsc:phi_f_for_individual_err} and use \eqref{eq:nsc:phi_f_for_individual_err_2}, 
{\small\begin{align}
	& \phi_{t+1}^* - f(\bar{x}(t+1)) \nonumber\\
	&\geq \sum_{k=0}^t \frac{1}{4} \eta_k \Vert g(k)\Vert^2 \prod_{\ell=k+1}^{t}(1-\alpha_\ell) \nonumber\\
	&	\quad -\sum_{k=0}^t 2 \kappa^2\chi_2(\eta_k)^2  L^3 \Vert  \bar{y}(k)- \bar{x}(k)  \Vert^2\prod_{\ell=k+1}^{t}(1-\alpha_\ell) \nonumber\\
	&\geq \frac{1}{4} \eta_t \Vert g(t)\Vert^2 - (\Phi_0(x^*) - f^*)\lambda_{t+1} \label{eq:nsc:individual:phi_f_lb}
	\end{align}}
On the other hand, we upper bound $\phi_{t+1}^* - f(\bar{x}(t+1))$
{\small \begin{align}
	\phi_{t+1}^*  - f(\bar{x}(t+1))  &\leq \Phi_{t+1}(x^*) - f^*   \leq \lambda_{t+1}(\Phi_0(x^*) - f^*)\label{eq:nsc:individual:phi_f_ub}
	\end{align}}where the last inequality follows from (\ref{eq:nsc:phi_upperbound}). Combining \eqref{eq:nsc:individual:phi_f_lb} and \eqref{eq:nsc:individual:phi_f_ub}, we have, 
$\Vert g(t)\Vert^2 \leq 8 \frac{\lambda_{t+1}}{\eta_t} (\Phi_0(x^*) - f^*) = O(\frac{1}{t^{2-2\beta}} )\Rightarrow \Vert g(t)\Vert = O(\frac{1}{t^{1-\beta}})$. 
Next, we upper bound $\Vert y(t) - \one y(t)\Vert $. Since we know that by $\eqref{eq:nsc:x_y_bound}$, $ \Vert \bar{x}(t) - \bar{y}(t) \Vert = O(\alpha_t) = O(\frac{1}{t})$, then by Lemma~\ref{lem:nsc:rel_con_err}, we have,
{\small\begin{align*}
	&\Vert y(t) - \one\bar{y}(t) \Vert \\
	&  \leq \kappa \sqrt{n} \chi_2(\eta_t) \bigg[L\Vert \bar{y}(t) - \bar{x}(t)\Vert + \frac{8}{1-\sigma} L \eta_t \Vert g(t)\Vert \bigg]\\
	& = O(\eta_t^{1/3} (\frac{1}{t} + \eta_t \frac{1}{t^{1-\beta}} ) ) = O(\frac{1}{t^{1+\frac{\beta}{3}}}).
	\end{align*}}Then, we have, $
\Vert \nabla f(\bar{y}(t)) \Vert    \leq  \Vert   \nabla f(\bar{y}(t)) - g(t)\Vert + \Vert g(t)\Vert
\leq \frac{L}{\sqrt{n}} \Vert y(t) - \one \bar{y}(t)\Vert + \Vert g(t)\Vert = O(\frac{1}{t^{1-\beta}})
$. Similarly,
$\Vert \nabla f(\bar{x}(t)) \Vert   \leq  \Vert   \nabla f(\bar{x}(t)) - \nabla f(\bar{y}(t))\Vert + \Vert \nabla f(\bar{y}(t))\Vert 
\leq  L \Vert    \bar{x}(t)) - \bar{y}(t) \Vert + \Vert \nabla f(\bar{y}(t))\Vert = O(\frac{1}{t^{1-\beta}})$.
We next upper bound $f(\bar{y}(t))$,
{\small
	\begin{align*}
	f(\bar{y}(t)) &\leq f(\bar{x}(t)) + \langle \nabla f(\bar{x}(t)), \bar{y}(t) - \bar{x}(t)\rangle + \frac{L}{2}  \Vert \bar{y}(t) - \bar{x}(t)\Vert^2\\
	&\leq f(\bar{x}(t)) + \Vert \nabla f(\bar{x}(t))\Vert \cdot \Vert \bar{y}(t) - \bar{x}(t)\Vert + \frac{L}{2}  \Vert \bar{y}(t) - \bar{x}(t)\Vert^2\\
	&= f^* + O(\frac{1}{t^{2-\beta}}  ) .
	\end{align*}}
Finally, for any $i$, we have 
{\small
	\begin{align*}
	f(y_i(t)) &\leq f(\bar{y}(t)) + \langle \nabla f(\bar{y}(t)) , y_i(t) - \bar{y}(t) \rangle + \frac{L}{2} \Vert y_i(t) - \bar{y}(t)\Vert^2\\
	&\leq f(\bar{y}(t)) + \Vert \nabla f(\bar{y}(t)) \Vert \cdot \Vert y_i(t) - \bar{y}(t) \Vert + \frac{L}{2} \Vert y_i(t) - \bar{y}(t)\Vert^2\\
	&= f(\bar{y}(t)) + O(\frac{1}{t^{2-\beta}}) = f^* + O(\frac{1}{t^{2-\beta}}).
	\end{align*}}So we have $f(y_i(t)) - f^* =  O(\frac{1}{t^{2-\beta}})$, the desired result of the Theorem. 

\subsection{Proof of Theorem~\ref{thm:nsc:fixed}}\label{subsec:nsc:fixed_full}
In this section, we provide a detailed proof for Theorem~\ref{thm:nsc:fixed}.
\begin{lemma}\label{lem:fixed:inexact_grad}
	Under the conditions of Theorem~\ref{thm:nsc:fixed}, we have inequality (\ref{eq:nsc:inexact:noisy_oracle_1}) in Lemma~\ref{lem:nsc:inexact_grad}, when evaluated at $\omega = \bar{x}(t)$, can be strengthened to,
	\begin{align*}
	f(\bar{x}(t)) &\geq \hat{f}(t) + \langle g(t), \bar{x}(t) - \bar{y}(t) \rangle + \frac{\mu}{4} \Vert \bar{x}(t) - \bar{y}(t)\Vert^2 - \frac{L}{2n}\Vert y(t) - \one\bar{y}(t)\Vert^2
	\end{align*}
	where $\mu = \mu_0 \gamma$, and $\gamma$ is the smallest non-zero eigenvalue of the positive semidefinite matrix $A = \frac{1}{n}  \sum_{i=1}^n A_i A_i^T $ (Matrix $A$ has at least one nonzero eigenvalue since otherwise, all $A_i$ would be zero.); $L = L_0\nu$ and $\nu = \max_i \Vert A\Vert_*^2$. 
\end{lemma} 
\noindent\textit{Proof: } It is easy to check that $A$ is a symmetric and positive semidefinite matrix. Let $y\in\R^{1\times N}$ be any vector in the null space of $A$. Notice,
$$ 0 =\langle  yA,y\rangle = \frac{1}{n} \sum_{i=1}^n \Vert y A_i  \Vert^2\geq 0.$$
This implies $\forall i$, $y A_i  = 0$. Then, 
\begin{align*}
\langle g(t), y\rangle &= \langle \frac{1}{n}\sum_{i=1}^n \nabla h_i(y_i(t)A_i) A_i^T, y\rangle =\frac{1}{n} \sum_{i=1}^n  \nabla h_i(y_i(t)A_i) A_i^T y^T = 0.
\end{align*}
This implies that, $g(t)$ lies in the space spanned by the rows of $A$. By (\ref{eq:nsc:v_y}) and the fact that $\bar{v}(0) - \bar{y}(0) = 0$, we have $\bar{v}(t) - \bar{y}(t)$ lies in the row space of $A$. Hence $\bar{x}(t) - \bar{y}(t) = \frac{\alpha_t}{1-\alpha_t} (\bar{y}(t) - \bar{v}(t))$ also lies in the row space of $A$. Therefore, 
$$\langle (\bar{x}(t) - \bar{y}(t))A,\bar{x}(t) - \bar{y}(t)\rangle \geq \gamma \Vert \bar{x}(t) - \bar{y}(t)\Vert^2. $$ 
Therefore,
\begin{align*}
f(\bar{x}(t))
& = \frac{1}{n} \sum_{i=1}^n h_i(\bar{x}(t)  A_i) \\
&\stackrel{(a)}{\geq} \frac{1}{n} \sum_{i=1}^n \big [ h_i(y_i(t) A_i) + \langle \nabla h_i(y_i(t)A_i), \bar{x}(t)A_i-y_i(t)A_i \rangle   + \frac{\mu_0}{2} \Vert (\bar{x}(t)-y_i(t))A_i  \Vert^2  \big]\\
&= \frac{1}{n} \sum_{i=1}^n \big [ f_i(y_i(t) ) + \langle \nabla f_i(y_i(t)), \bar{x}(t)-y_i(t) \rangle + \frac{\mu_0}{2} \Vert ( \bar{x}(t)-y_i(t))A_i  \Vert^2  \big]\\
&\geq  \frac{1}{n} \sum_{i=1}^n \big [ f_i(y_i(t)) + \langle \nabla f_i(y_i(t)), \bar{y}(t)-y_i(t) \rangle \big]  +  \frac{1}{n} \sum_{i=1}^n  \langle \nabla f_i(y_i(t)), \bar{x}(t) - \bar{y}(t)\rangle   \\
&\qquad +  \frac{\mu_0}{2n}\sum_{i=1}^n \Big[ \frac{1}{2}\Vert ( \bar{x}(t) - \bar{y}(t))A_i \Vert^2 - \Vert( \bar{y}(t) - y_i(t))A_i\Vert^2 \Big] \\
&\stackrel{(b)}{\geq} \hat{f}(t) + \langle g(t), \bar{x}(t) - \bar{y}(t)\rangle  + \frac{\mu_0}{4n} \sum_{i=1}^n \langle (\bar{x}(t) - \bar{y}(t)) A_i A_i^T, \bar{x}(t) - \bar{y}(t)\rangle\\
&\qquad  - \frac{\mu_0 \nu }{2n}  \sum_{i=1}^n \Vert y_i(t)-\bar{y}(t)\Vert^2 \\
&= \hat{f}(t) + \langle g(t), \bar{x}(t) - \bar{y}(t)\rangle + \frac{\mu_0}{4} \langle (\bar{x}(t) - \bar{y}(t)) A, \bar{x}(t) - \bar{y}(t)\rangle  - \frac{\mu_0 \nu }{2n} \Vert y(t)-\one\bar{y}(t)\Vert^2 \\
&\geq \hat{f}(t) + \langle g(t), \bar{x}(t) - \bar{y}(t)\rangle + \frac{1}{4} \mu_0\gamma \Vert \bar{x}(t) - \bar{y}(t)\Vert^2   - \frac{L }{2}\frac{1}{n}  \Vert y(t) - \one\bar{y}(t)\Vert^2
\end{align*}
where in (a) we have used the fact that $h_i(\cdot)$ is $\mu_0$-strongly convex, and in (b), we have used the definition of $\hat{f}(t)$, and the fact that $\Vert (\bar{y}(t) - y_i(t))A_i\Vert\leq \sqrt{\nu} \Vert \bar{y}(t) - y_i(t)\Vert$ (since $A_i$'s spectral norm is upper bounded by $\sqrt{\nu}$). 
\qedd 

We now finish the proof of Theorem \ref{thm:nsc:fixed}.

\noindent\textit{Proof of Theorem \ref{thm:nsc:fixed}:} It is easy to check that $f_i$ is convex and $L$-smooth. It is easy to check under the step size conditions, Lemma~\ref{lem:nsc:rel_con_err} holds. We will prove by induction that,
\begin{equation}
\phi_t^* \geq f(\bar{x}(t)). \label{eq:fixed:induction}
\end{equation}
Equation (\ref{eq:fixed:induction}) is true for $t=0$. Next, by (\ref{eq:nsc:inexact:phi}),  
\begin{align}
\phi_{t+1}^*  
&= (1 - \alpha_t )\phi_t^*  + \alpha_t \hat{f}(t) - \frac{1}{2}\eta_t \Vert g(t)\Vert^2  +  \alpha_t \langle g(t), \bar{v}(t) - \bar{y}(t)\rangle  \nonumber \nonumber\\
&\stackrel{(a)}{\geq} (1 - \alpha_t) f(\bar{x}(t))      + \alpha_t \hat{f}(t) - \frac{1}{2}\eta_t \Vert g(t)\Vert^2  +  \alpha_t \langle g(t), \bar{v}(t) - \bar{y}(t)\rangle\nonumber\\
&\stackrel{(b)}{\geq} (1 - \alpha_t )\{ \hat{f}(t) + \langle g(t), \bar{x}(t) - \bar{y}(t)\rangle  + \frac{\mu}{4} \Vert\bar{x}(t) - \bar{y}(t)\Vert^2 \nonumber\\
&\qquad - \frac{L}{2n} \Vert y(t) - \one\bar{y}(t)\Vert^2 \}  + \alpha_t \hat{f}(t) - \frac{1}{2}\eta_t  \Vert g(t)\Vert^2    +  \alpha_t  \langle g(t), \bar{v}(t) - \bar{y}(t)\rangle \nonumber\\
&\stackrel{(c)}{\geq}   \hat{f}(t)  - \frac{1}{2}\eta_t \Vert g(t)\Vert^2  +(1-\alpha_t)\frac{\mu}{4} \Vert \bar{x}(t) - \bar{y}(t)\Vert^2   -  \frac{L}{2n} \Vert y(t) - \one\bar{y}(t)\Vert^2.   \label{eq:fixed:phi_f_recursive}
\end{align}
where (a) is due to the induction assumption (\ref{eq:fixed:induction}), (b) is due to Lemma \ref{lem:fixed:inexact_grad} and (c) is due to $ \alpha_t (\bar{v}(t) - \bar{y}(t)) + (1-\alpha_t) (\bar{x}(t) - \bar{y}(t)) =0$.
By (\ref{eq:nsc:inexact:noisy_oracle_2}) (Lemma \ref{lem:nsc:inexact_grad}), 
\begin{align*}
f(\bar{x}(t+1)) 
&\leq \hat{f}(t) + \langle g(t), \bar{x}(t+1) - \bar{y}(t)\rangle  + L\Vert \bar{x}(t+1) - \bar{y}(t)\Vert^2 + \frac{L}{n}\Vert y(t)-\one\bar{y}(t)\Vert^2\\
&\leq  \hat{f}(t) -(\eta_t -L\eta_t^2)\Vert g(t)\Vert^2 +\frac{L}{n} \Vert y(t) - \one\bar{y}(t)\Vert^2. 
\end{align*}
Combining the above with (\ref{eq:fixed:phi_f_recursive}) and using Lemma \ref{lem:nsc:rel_con_err}, we have,
\begin{align}
\phi_{t+1}^* - f(\bar{x}(t+1)) \nonumber
&\geq ( \frac{1}{2}\eta_t - L\eta_t^2) \Vert g(t)\Vert^2 + (1-\alpha_t)\frac{\mu}{4}\Vert\bar{x}(t)-\bar{y}(t)\Vert^2 - \frac{3}{2}\frac{L}{n}\Vert y(t) - \one\bar{y}(t)\Vert^2\nonumber\\
&\geq ( \frac{1}{2}\eta_t - L\eta_t^2) \Vert g(t)\Vert^2 + (1-\alpha_t)\frac{\mu}{4}\Vert\bar{x}(t)-\bar{y}(t)\Vert^2\nonumber\\
&\quad- 3 L \kappa^2 \chi_2(\eta_t)^2  ( L^2 \Vert \bar{x}(t) - \bar{y}(t)\Vert^2 +  \frac{64}{(1-\sigma)^2}L^2\eta_t^2 \Vert g(t)\Vert^2 )\nonumber\\
&= ( \frac{1}{2}\eta_t - L\eta_t^2 - \frac{192 L^3\kappa^2 \chi_2(\eta_t )^2 \eta_t^2}{(1-\sigma)^2} )\Vert g(t)\Vert^2\nonumber  \\
& \quad +( (1-\alpha_t) \frac{\mu}{4} - 3 \kappa^2\chi_2(\eta_t)^2L^3)\Vert\bar{x}(t) - \bar{y}(t)\Vert^2.\label{eq:fixed:phi_f}
\end{align}
Since $\eta_t = \eta$, and $\chi_2(\eta) < \frac{2 \eta^{1/3}}{L^{2/3}}$, and recalling $\kappa = \frac{6}{1-\sigma}$, we have 
{\color{black}\begin{align*}
\frac{1}{2}\eta_t - L\eta_t^2 - \frac{192 L^3\kappa^2 \chi_2(\eta_t )^2 \eta_t^2}{(1-\sigma)^2} & \geq \frac{1}{2}\eta - L\eta^2 - \frac{27648}{(1-\sigma)^4} \eta \cdot (\eta L)^{5/3}\\
&= \eta (\frac{1}{2} - L\eta - \frac{27648}{(1-\sigma)^4}  (\eta L)^{5/3})\geq \frac{1}{4}\eta 
\end{align*}
{\color{black}where in the last inequality, we have used $\eta L < \frac{1}{8}$, and $\frac{27648}{(1-\sigma)^4}  (\eta L)^{5/3} < \frac{1}{8}$ ($\Leftarrow \eta L< \frac{(1-\sigma)^{2.4}}{1611}$), all following from our step size condition. }

Next, since $\alpha_t<\alpha_0\leq \frac{1}{2}$ and $\chi_2(\eta) < \frac{2 \eta^{1/3}}{L^{2/3}}$, we have, 
\begin{align*}
& (1-\alpha_t) \frac{\mu}{4} - 3 \kappa^2\chi_2(\eta_t)^2L^3 \geq \frac{\mu}{8} - \frac{432}{(1-\sigma)^2} \eta^{2/3} L^{5/3} \geq \frac{\mu}{16}
\end{align*}
where the last inequality (equivalent to $ \eta^{2/3}< \frac{\mu}{L^{5/3}} \frac{(1-\sigma)^2}{6912}$) follows from the step size condition. Hence, returning to (\ref{eq:fixed:phi_f}), we get 
\begin{align}
\phi_{t+1}^* &\geq f(\bar{x}(t+1)) + \frac{1}{4}\eta \Vert g(t)\Vert^2 + \frac{\mu}{16} \Vert\bar{x}(t) - \bar{y}(t)\Vert^2 \label{eq:fixed:phi_f_ind_err}\\
&\geq f(\bar{x}(t+1)) \nonumber
\end{align}}
Therefore the induction is finished and (\ref{eq:fixed:induction}) is true for all $t$.  Hence, by (\ref{eq:nsc:phi_upperbound}),
\begin{align*}
f(\bar{x}(t))&\leq \phi_t^*   \leq \Phi_t(x^*)  \leq f^* +\lambda_t(\Phi_0(x^*) - f^*).
\end{align*}
Therefore $f(\bar{x}(t)) - f^* = O(\lambda_t)$. Following an argument similar to Lemma~\ref{lem:nsc:alpha} (or simply let $\beta\rightarrow 0$ in Lemma~\ref{lem:nsc:alpha}), we will have $\lambda_t = O(1/t^2)$. As a result, $f(\bar{x}(t)) - f^* = O(\frac{1}{t^2})$, i.e. part (a) of the Theorem.

{\color{black} Finally, we prove part (b) of the Theorem, i.e. upper bounding $f(y_i(t)) - f^*$. We start from \eqref{eq:fixed:phi_f_ind_err}, and using (\ref{eq:nsc:phi_upperbound}), we have
	\begin{align}
	\frac{1}{4}\eta\Vert g(t) \Vert^2 + \frac{\mu}{16} \Vert\bar{x}(t) - \bar{y}(t)\Vert^2 &\leq \phi_{t+1}^* - f(\bar{x}(t+1)) \nonumber \\
	&\leq \Phi_{t+1}(x^*)  - f^* \nonumber \\
	&\leq \lambda_{t+1}  (\Phi_0(x^*) - f^*) \nonumber.
	\end{align}
Therefore, $\Vert g(t)\Vert = O(\frac{1}{t})$, and $\Vert \bar{x}(t)-\bar{y}(t) \Vert =  O(\frac{1}{t})$. Using Lemma~\ref{lem:nsc:rel_con_err}, we then have $\Vert y(t) - \one \bar{y}(t)\Vert= O(\Vert \bar{x}(t)-\bar{y}(t) \Vert ) + O(\Vert g(t)\Vert) = O(\frac{1}{t})$. Then, we have, 
\begin{align*}
\Vert \nabla f(\bar{y}(t)) \Vert &  \leq  \Vert   \nabla f(\bar{y}(t)) - g(t)\Vert + \Vert g(t)\Vert \leq \frac{L}{\sqrt{n}} \Vert y(t) - \one \bar{y}(t)\Vert + \Vert g(t)\Vert = O(\frac{1}{t})
\end{align*}
And similarly,
\begin{align*}
\Vert \nabla f(\bar{x}(t)) \Vert &  \leq  \Vert   \nabla f(\bar{x}(t)) - \nabla f(\bar{y}(t))\Vert + \Vert \nabla f(\bar{y}(t))\Vert 
\leq  L \Vert    \bar{x}(t)) - \bar{y}(t) \Vert + \Vert \nabla f(\bar{y}(t))\Vert 
= O(\frac{1}{t})
\end{align*}
We next upper bound $f(\bar{y}(t))$,
{\small
	\begin{align*}
	f(\bar{y}(t)) &\leq f(\bar{x}(t)) + \langle \nabla f(\bar{x}(t)), \bar{y}(t) - \bar{x}(t)\rangle + \frac{L}{2}  \Vert \bar{y}(t) - \bar{x}(t)\Vert^2\\
	&\leq f(\bar{x}(t)) + \Vert \nabla f(\bar{x}(t))\Vert \cdot \Vert \bar{y}(t) - \bar{x}(t)\Vert + \frac{L}{2}  \Vert \bar{y}(t) - \bar{x}(t)\Vert^2\\
	&= f(\bar{x}(t)) + O(\frac{1}{t} \cdot \frac{1}{t} ) + O(\frac{1}{t^2}) \\
	&= f^* + O(\frac{1}{t^{2}}  ) .
	\end{align*}}
Finally, for any $i$, we have 
{\small
	\begin{align*}
	f(y_i(t)) &\leq f(\bar{y}(t)) + \langle \nabla f(\bar{y}(t)) , y_i(t) - \bar{y}(t) \rangle + \frac{L}{2} \Vert y_i(t) - \bar{y}(t)\Vert^2\\
	&\leq f(\bar{y}(t)) + \Vert \nabla f(\bar{y}(t)) \Vert \cdot \Vert y_i(t) - \bar{y}(t) \Vert + \frac{L}{2} \Vert y_i(t) - \bar{y}(t)\Vert^2\\
	&\leq f(\bar{y}(t)) + O(\frac{1}{t} \cdot \frac{1}{t} ) + O( \frac{1}{t^{2}})\\
	&= f(\bar{y}(t)) + O(\frac{1}{t^{2}}) = f^* + O(\frac{1}{t^{2}}).
	\end{align*}}So we have $f(y_i(t)) - f^* =  O(\frac{1}{t^{2}})$, the desired result of the Theorem. 
}
\subsection{Supplementary Materials for the Simulation in Section~\ref{sec:numerical}}\label{appendix:stepsize}
In this section, we provide the step size parameters for the five figures in Section~\ref{sec:numerical}. We also provide the individual objective errors $f(y_i(t)) - f^*$ achieved by our algorithm for Figure~\ref{fig:case1}, Figure~\ref{fig:case2} and Figure \ref{fig:case3}. 

\textbf{Step sizes for Fig.~\ref{fig:case1}:} (a) For random graph (left plot), $\frac{L}{\mu} = 793.1463, \sigma = 0.59052$ and step sizes are, Acc-DNGD-SC: $\eta = 0.00017,\alpha = 0.011821$; D-NG: $\eta_t = \frac{0.00076687}{t+1}$; DGD: $\eta_t = \frac{0.0015337}{\sqrt{t}}$; EXTRA: $\eta = 0.00092025$; Acc-DGD: $\eta_t = 0.00030675$; CGD: $\eta = 0.0015337$; CNGD: $\eta =0.0015337$, $\alpha = 0.035508$.
(b) For $k$-cycle (middle plot), $\frac{L}{\mu} = 793.1463, \sigma = 0.74566$ and step sizes are, Acc-DNGD-SC: $\eta = 0.00013,\alpha = 0.010338$; D-NG: $\eta_t = \frac{0.00076687}{t+1}$; DGD: $\eta_t = \frac{0.0015337}{\sqrt{t}}$; EXTRA: $\eta = 0.00092025$; Acc-DGD: $\eta_t = 0.00015337$; CGD: $\eta = 0.0015337$; CNGD: $\eta =0.0015337$, $\alpha = 0.035508$.
(c) For 2D grid (right plot), $\frac{L}{\mu} = 772.5792, \sigma = 0.92361$ and step sizes are, Acc-DNGD-SC: $\eta = 0.00005,\alpha = 0.0064959$; D-NG: $\eta_t = \frac{0.00076687}{t+1}$; DGD: $\eta_t = \frac{0.0015337}{\sqrt{t}}$; EXTRA: $\eta = 0.00092025$; Acc-DGD: $\eta_t = 0.00010736$; CGD: $\eta = 0.0015337$; CNGD: $\eta =0.0015337$, $\alpha = 0.035977$.

\textbf{Step sizes for Fig.~\ref{fig:case2}:} (a) For random graph (left plot), $\frac{L}{\mu} = 658.0205, \sigma = 0.59052$ and step sizes are, Acc-DNGD-SC: $\eta = 0.03,\alpha = 0.016707$; D-NG: $\eta_t = \frac{0.16334}{t+1}$; DGD: $\eta_t = \frac{0.32667}{\sqrt{t}}$; EXTRA: $\eta = 0.16334$; Acc-DGD: $\eta_t = 0.081669$; CGD: $\eta = 0.32667$; CNGD: $\eta =0.32667$, $\alpha = 0.055131$.
(b) For $k$-cycle (middle plot), $\frac{L}{\mu} = 658.0205, \sigma = 0.74566$ and step sizes are, Acc-DNGD-SC: $\eta = 0.015,\alpha = 0.011814$; D-NG: $\eta_t = \frac{0.16334}{t+1}$; DGD: $\eta_t = \frac{0.32667}{\sqrt{t}}$; EXTRA: $\eta = 0.081669$; Acc-DGD: $\eta_t = 0.032667$; CGD: $\eta = 0.32667$; CNGD: $\eta =0.32667$, $\alpha = 0.055131$.
(c) For 2D grid (right plot), $\frac{L}{\mu} = 344.5099, \sigma = 0.92361$ and step sizes are, Acc-DNGD-SC: $\eta = 0.015,\alpha = 0.014345$; D-NG: $\eta_t = \frac{0.2116}{t+1}$; DGD: $\eta_t = \frac{0.42319}{\sqrt{t}}$; EXTRA: $\eta = 0.2116$; Acc-DGD: $\eta_t = 0.063479$; CGD: $\eta = 0.42319$; CNGD: $\eta =0.42319$, $\alpha = 0.076193$.

\textbf{Step sizes for Fig.~\ref{fig:case3}}: (a) For random graph (left plot), $\sigma = 0.59052$ and step sizes are, Acc-DNGD-NSC with vanishing step size: $\eta = \frac{0.0027642}{(t+1)^{0.61}},\alpha_0= 0.70711$; Acc-DNGD-NSC with fixed step size: $\eta_t = 0.0027642,\alpha_0=0.70711$; D-NG: $\eta_t = \frac{0.0027642}{t+1}$; DGD: $\eta_t = \frac{0.0055283}{\sqrt{t}}$; EXTRA: $\eta = 0.0055283$; Acc-DGD: $\eta_t = 0.0027642$; CGD: $\eta = 0.0055283$; CNGD: $\eta =0.0055283$, $\alpha_0 = 0.5$. 
	(b) For $k$-cycle (middle plot), $\sigma = 0.74566$ and step sizes are, Acc-DNGD-NSC with vanishing step size: $\eta = \frac{0.0027642}{(t+1)^{0.61}},\alpha_0= 0.70711$; Acc-DNGD-NSC with fixed step size: $\eta_t = 0.0022113,\alpha_0=0.63246$; D-NG: $\eta_t = \frac{0.0027642}{t+1}$; DGD: $\eta_t = \frac{0.0055283}{\sqrt{t}}$; EXTRA: $\eta = 0.0055283$; Acc-DGD: $\eta_t = 0.0022113$; CGD: $\eta = 0.0055283$; CNGD: $\eta =0.0055283$, $\alpha_0 = 0.5$. (c) For 2D grid (right plot), $\sigma = 0.92361$ and step sizes are, Acc-DNGD-NSC with vanishing step size: $\eta = \frac{0.0024928}{(t+1)^{0.61}},\alpha_0= 0.70711$; Acc-DNGD-NSC with fixed step size: $\eta_t = 0.0014957,\alpha_0=0.54772$; D-NG: $\eta_t = \frac{0.0024928}{t+1}$; DGD: $\eta_t = \frac{0.0049856}{\sqrt{t}}$; EXTRA: $\eta = 0.0049856$; Acc-DGD: $\eta_t = 0.0014957$; CGD: $\eta = 0.0049856$; CNGD: $\eta =0.0049856$, $\alpha_0 = 0.5$.

\textbf{Step sizes for Fig.~\ref{fig:individual_err}:} The step size of the upper plot (case I) is identical to that of the Acc-DNGD-SC in Fig.~\ref{fig:case1} left plot. The step size of the lower plot (case III) is identical to that of the Acc-DNGD-NSC with fixed step size in Fig.~\ref{fig:case3} left plot. 

\textbf{Step sizes for Fig.~\ref{fig:timevarying}:}  (a) The upper plot uses the cost function of case I, with $\frac{L}{\mu} = 772.5792$ and the step sizes are $\eta = 0.000011717,\alpha = 0.0031445$; 
(b) The lower plot uses the cost function of case III and the step sizes is $\eta_t = 0.0014957,\alpha_0 = 0.54772$. 

\begin{figure*}[t]
	\begin{subfigure}[b]{0.32\textwidth}
		\centering
		\includegraphics[width=\textwidth]{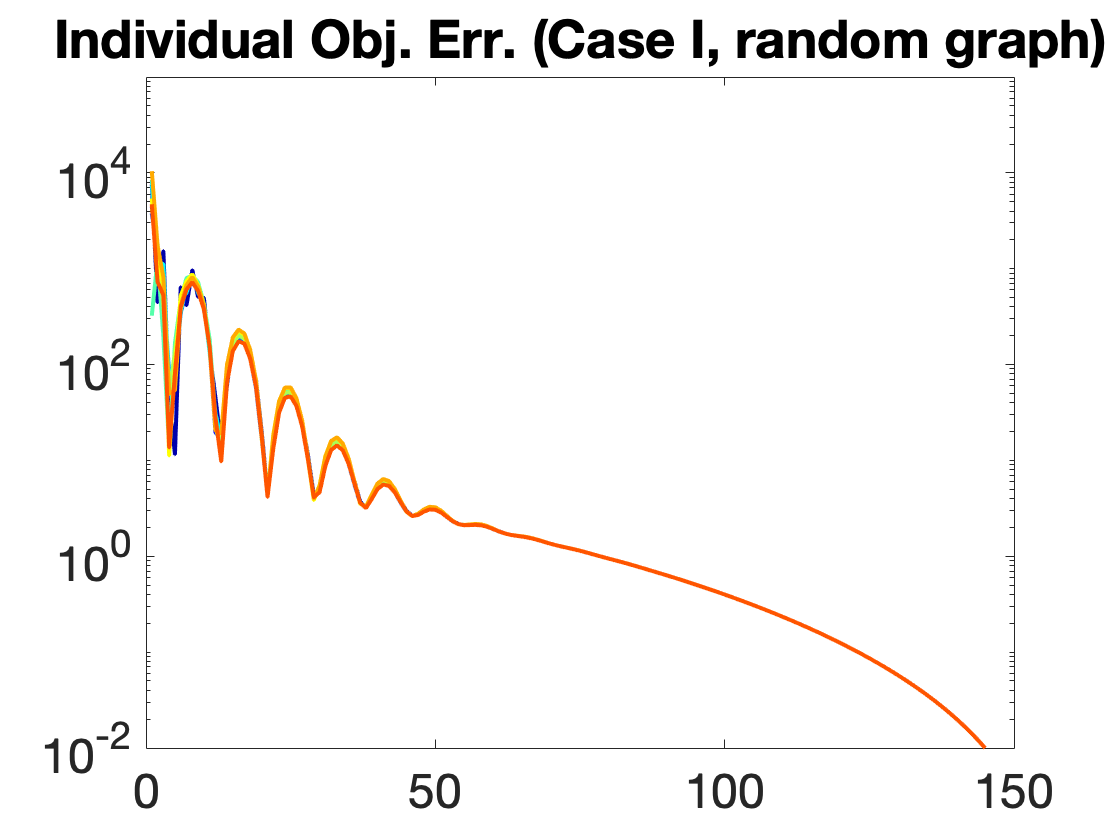} 
	\end{subfigure}
	\begin{subfigure}[b]{0.32\textwidth}
		\centering
		\includegraphics[width=\textwidth]{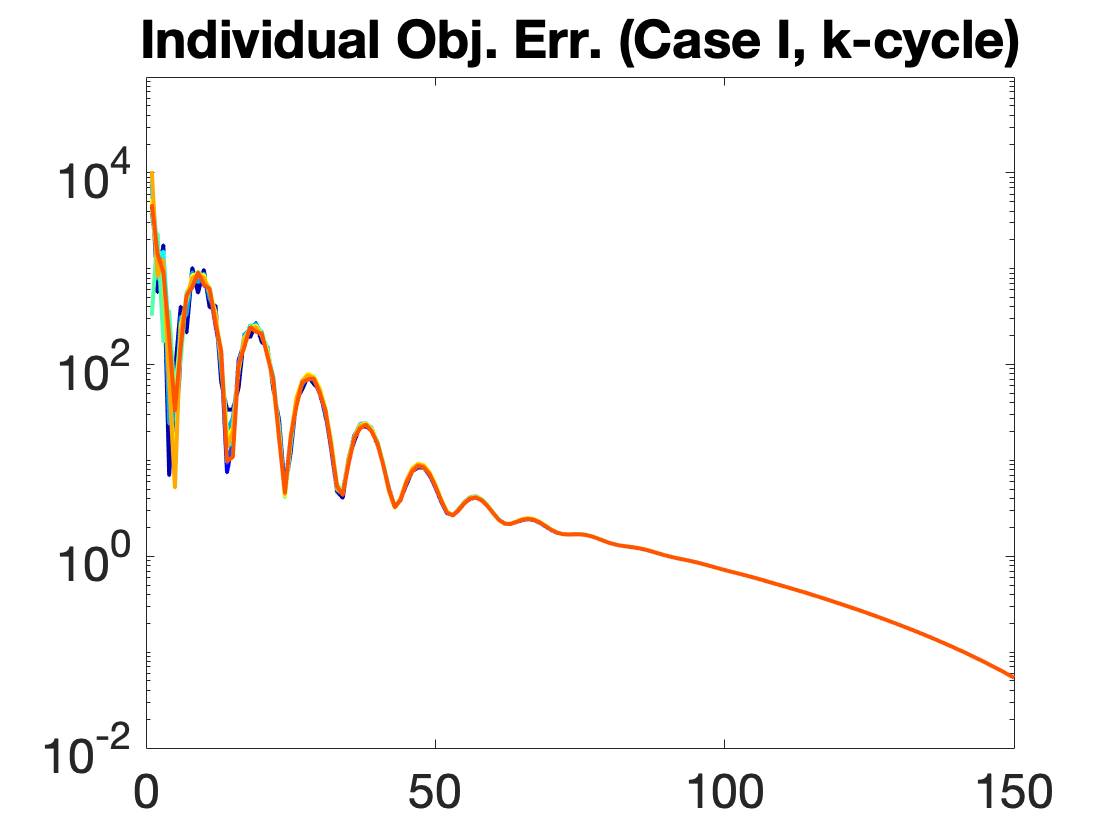} 
	\end{subfigure}
	\begin{subfigure}[b]{0.32\textwidth}
		\centering
		\includegraphics[width=\textwidth]{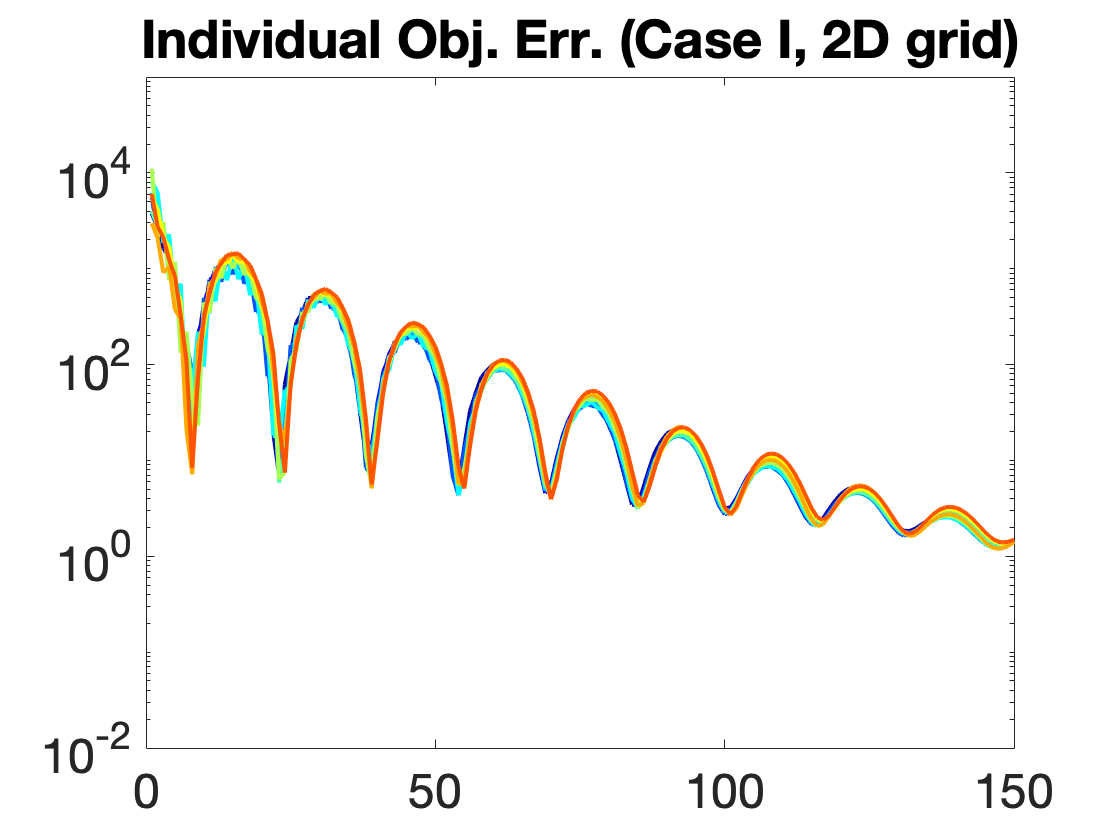} 
	\end{subfigure}
\caption{Individual objective error $f(y_i(t)) - f^*$ achieved by our algorithm Acc-DNGD-SC in Case I (Figure~\ref{fig:case1}), for $10$ choices of $i$ selected evenly spaced from index $1,2,\ldots,n$. }\label{fig:case1_ind}
\end{figure*}

\begin{figure*}[t]
	\begin{subfigure}[b]{0.32\textwidth}
		\centering
		\includegraphics[width=\textwidth]{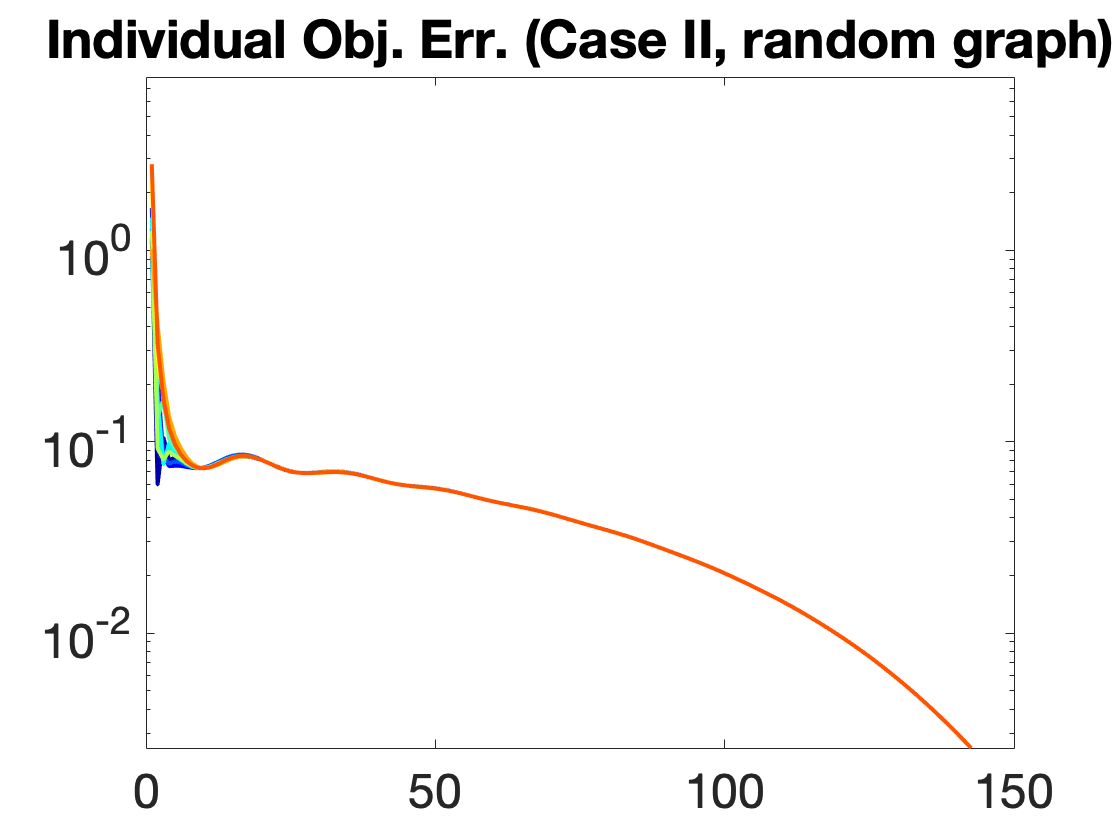} 
	\end{subfigure}
	\begin{subfigure}[b]{0.32\textwidth}
		\centering
		\includegraphics[width=\textwidth]{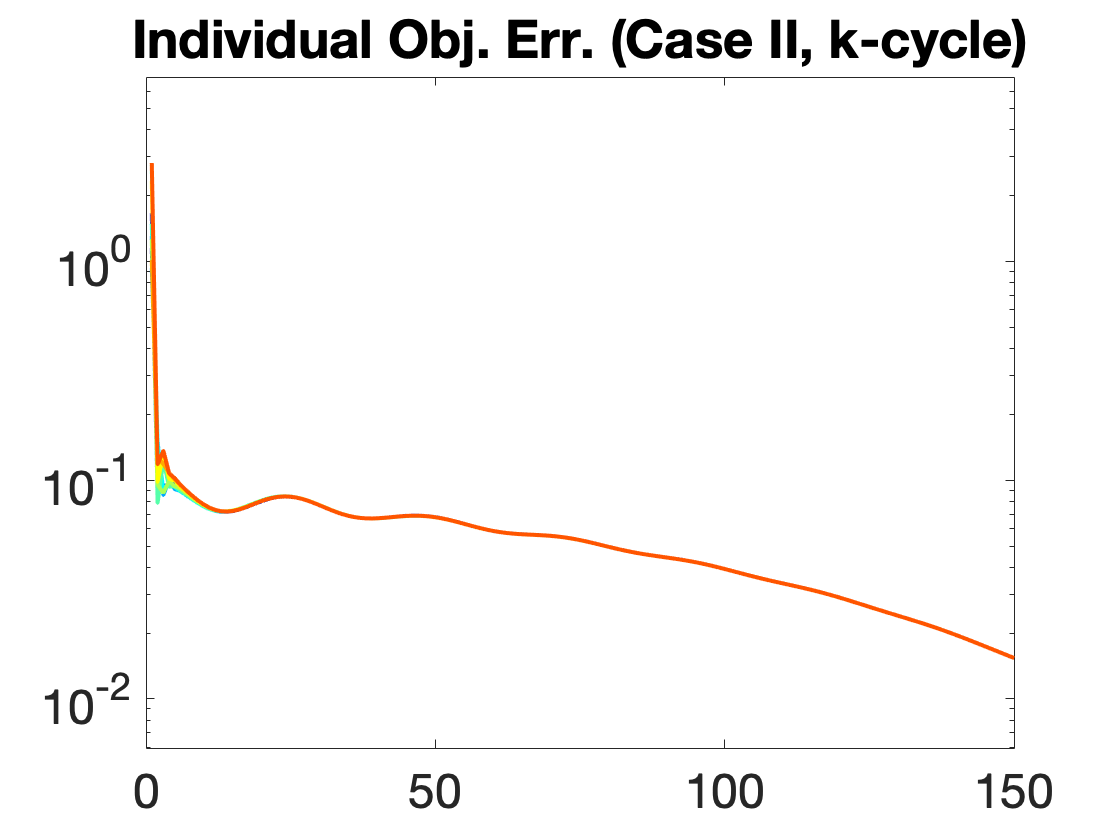} 
	\end{subfigure}
	\begin{subfigure}[b]{0.32\textwidth}
		\centering
		\includegraphics[width=\textwidth]{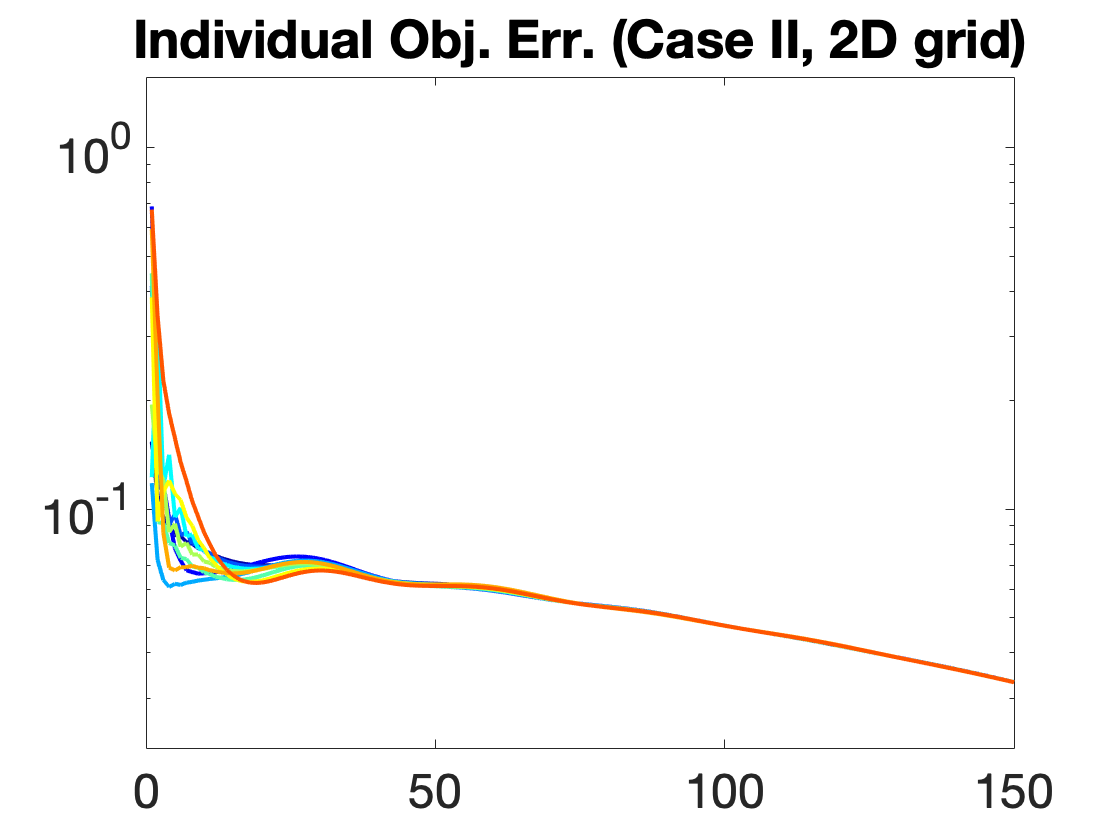} 
	\end{subfigure}
\caption{Individual objective error $f(y_i(t)) - f^*$ achieved by our algorithm Acc-DNGD-SC in Case II (Figure~\ref{fig:case2}), for $10$ choices of $i$ selected evenly spaced from index $1,2,\ldots,n$. }\label{fig:case2_ind}
\end{figure*}

\begin{figure*}[t]
	\begin{subfigure}[b]{0.32\textwidth}
		\centering
		\includegraphics[width=\textwidth]{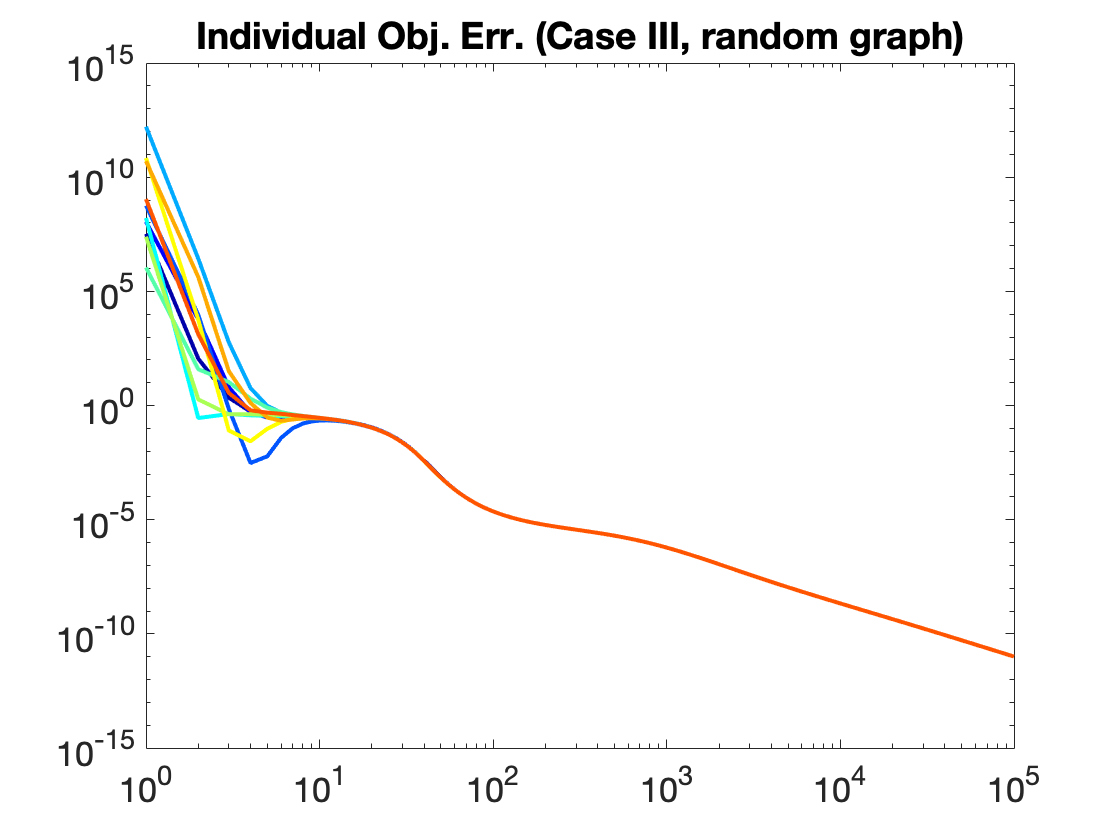} 
	\end{subfigure}
	\begin{subfigure}[b]{0.32\textwidth}
		\centering
		\includegraphics[width=\textwidth]{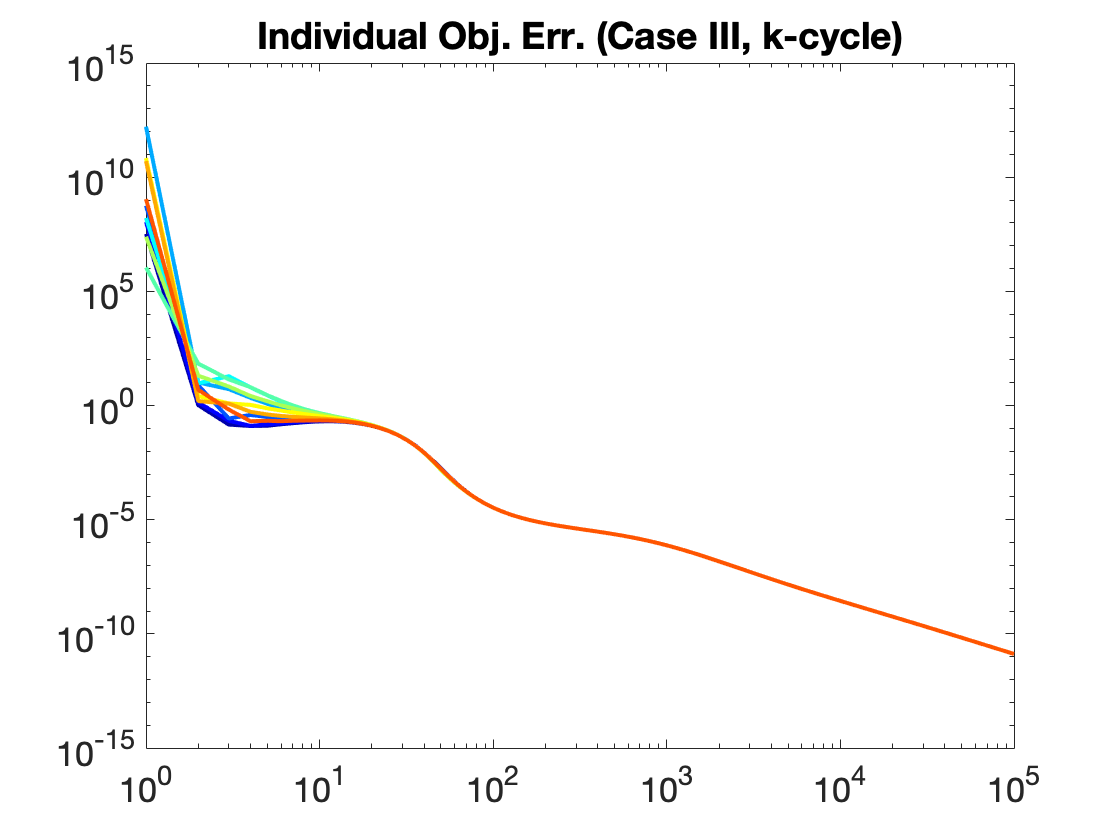} 
	\end{subfigure}
	\begin{subfigure}[b]{0.32\textwidth}
		\centering
		\includegraphics[width=\textwidth]{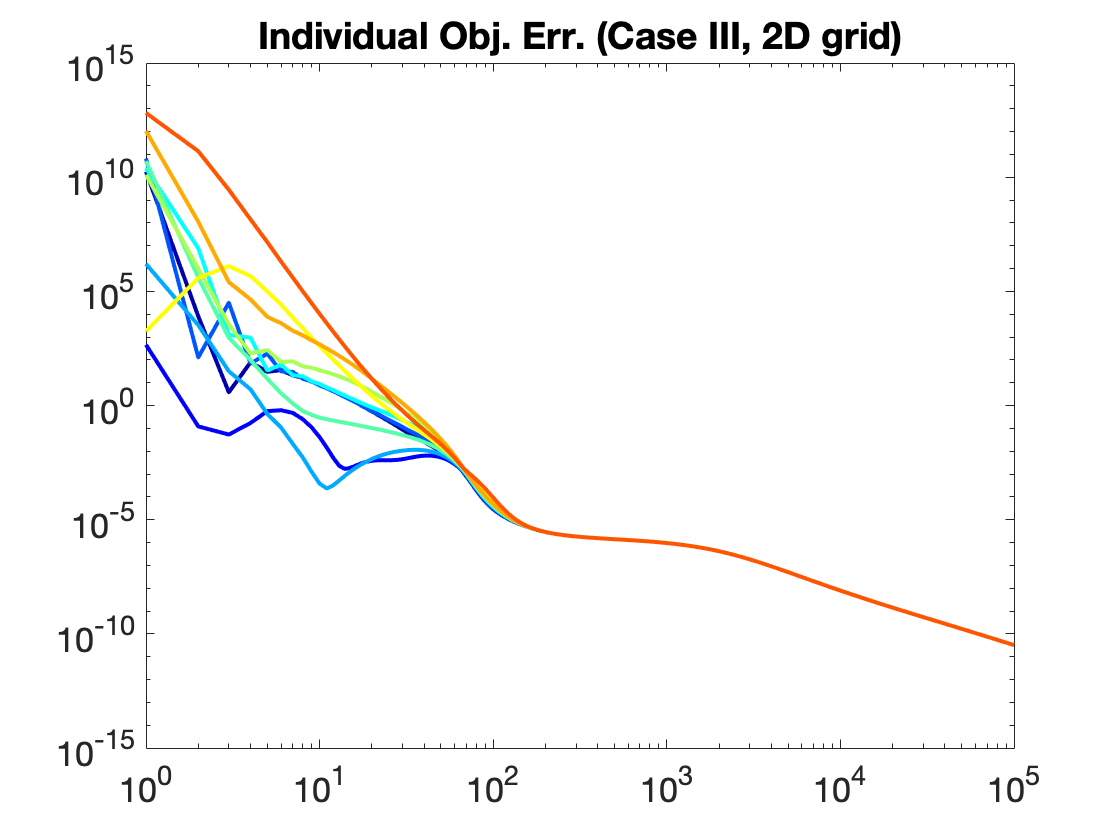} 
	\end{subfigure}
\caption{Individual objective error $f(y_i(t)) - f^*$ achieved by our algorithm Acc-DNGD-NSC (with $\beta=0$) in Case III (Figure~\ref{fig:case3}), for $10$ choices of $i$ selected evenly spaced from index $1,2,\ldots,n$. }\label{fig:case3_ind}
\end{figure*}

\textbf{Individual Objective Errors.} In Figure~\ref{fig:case1_ind}, Figure~\ref{fig:case2_ind} and Figure~\ref{fig:case3_ind}, we provide the individual objective errors $f(y_i(t)) - f^*$ achieved by our algorithm in Case I (Figure~\ref{fig:case1}), Case II (Figure~\ref{fig:case2}) and Case III (Figure \ref{fig:case3}) respectively. 


\subsection{Proof of Lemma \ref{lem:G23_aux}} \label{subsec:G23}
\noindent\textit{Proof of Lemma \ref{lem:G23_aux}:} We have,
$$\left[\begin{array}{ccc}
1 & 1 & 1\\
\gamma_1 & \gamma_2 & \gamma_3 \\
\gamma_1^2 & \gamma_2^2 & \gamma_3^2
\end{array}\right] \left[\begin{array}{c}
\alpha_1\\ \alpha_2 \\ \alpha_3 \end{array}\right]
= \left[\begin{array}{c}
B_0 \\ B_1 \\ B_2
\end{array}\right].$$
Hence we can solve for $\alpha_1,\alpha_2,\alpha_3$, getting
\begin{align*}
\alpha_1 &=\frac{1}{\Delta} \bigg[\gamma_2\gamma_3(\gamma_3-\gamma_2)B_0 + (\gamma_2^2 - \gamma_3^2) B_1  + (\gamma_3-\gamma_2)B_2 \bigg]  \\
\alpha_2 &=\frac{1}{\Delta} \bigg[\gamma_3\gamma_1(\gamma_1-\gamma_3)B_0 + (\gamma_3^2 - \gamma_1^2) B_1  + (\gamma_1-\gamma_3)B_2 \bigg]  \\
\alpha_3 &=\frac{1}{\Delta} \bigg[\gamma_1\gamma_2(\gamma_2-\gamma_1)B_0  +(\gamma_1^2 - \gamma_2^2) B_1  + (\gamma_2-\gamma_1)B_2 \bigg] \\
\end{align*}
where $\Delta = (\gamma_1 - \gamma_2)(\gamma_2 - \gamma_3)(\gamma_3 - \gamma_1)$. 
We now calculate,
\begin{align*}
\vert (\gamma_1 - \gamma_3)\gamma_2^t + (\gamma_2 - \gamma_1)\gamma_3^t\vert & = \vert (\gamma_2 - \gamma_3) \gamma_2^t + (\gamma_2 - \gamma_1)(\gamma_3^t - \gamma_2^t)     \vert\\
&= \vert (\gamma_2 - \gamma_3) \gamma_2^t + (\gamma_2 - \gamma_1)(\gamma_3 - \gamma_2)\sum_{k=0}^{t-1} \gamma_3^{t-1-k}\gamma_2^k     \vert  \\
& \leq \vert \gamma_2 - \gamma_3\vert \Big[ \vert\gamma_2\vert^t + t \vert \gamma_1 - \gamma_2\vert  \vert \gamma_2\vert^{t-1} \Big].\\
\end{align*}
Similarly,
\begin{align*}
\vert (\gamma_3^2 - \gamma_1^2)\gamma_2^t + (\gamma_1^2 - \gamma_2^2)\gamma_3^t\vert & = \vert (\gamma_3^2 - \gamma_2^2) \gamma_2^t + (\gamma_1^2 - \gamma_2^2)(\gamma_3^t - \gamma_2^t)     \vert\\
& \leq 2 \vert \gamma_2 - \gamma_3\vert \vert\gamma_2\vert^t + 2\vert \gamma_1 - \gamma_2\vert \vert\gamma_2 - \gamma_3\vert \vert \sum_{k=0}^{t-1} \gamma_2^{t-1-k}\gamma_3^k\vert \\
&\leq \vert \gamma_2 - \gamma_3\vert \Big[  2 \vert\gamma_2\vert^t + 2t \vert \gamma_1-\gamma_2\vert \vert\gamma_2\vert^{t-1}\Big].
\end{align*} 

\begin{align*}
\vert \gamma_3\gamma_1(\gamma_1 - \gamma_3)\gamma_2^t + \gamma_1\gamma_2(\gamma_2 - \gamma_1)\gamma_3^t\vert & = \vert (\gamma_3-\gamma_2)(\gamma_1^2 - \gamma_1(\gamma_2+\gamma_3)) \gamma_2^t + \gamma_1\gamma_2(\gamma_2 - \gamma_1)(\gamma_3^t - \gamma_2^t)     \vert\\
& \leq 3 \vert \gamma_2 - \gamma_3\vert \vert\gamma_2\vert^t + \vert \gamma_1 - \gamma_2\vert \vert\gamma_2 - \gamma_3\vert \vert \sum_{k=0}^{t-1} \gamma_2^{t-1-k}\gamma_3^k\vert \\
&\leq \vert \gamma_2 - \gamma_3\vert \Big[  3 \vert\gamma_2\vert^t + t \vert \gamma_1-\gamma_2\vert \vert\gamma_2\vert^{t-1}\Big].
\end{align*}
Hence,
\begin{align*}
&|\alpha_2\gamma_2^t + \alpha_3\gamma_3^t| \\
&\leq \frac{1}{\vert\Delta\vert } \Big[ |B_0|\vert \gamma_3\gamma_1(\gamma_1 - \gamma_3)\gamma_2^t + \gamma_1\gamma_2(\gamma_2 - \gamma_1)\gamma_3^t\vert + |B_1|\vert (\gamma_3^2 - \gamma_1^2)\gamma_2^t + (\gamma_1^2 - \gamma_2^2)\gamma_3^t\vert  \\
&\qquad + \vert B_2\vert \vert (\gamma_1 - \gamma_3)\gamma_2^t + (\gamma_2 - \gamma_1)\gamma_3^t\vert   \Big]\\
&\leq \frac{\vert \gamma_2-\gamma_3\vert }{\vert\Delta\vert }\bigg[(3\vert B_0 \vert+ 2\vert B_1\vert + \vert B_2\vert ) \vert\gamma_2\vert^t + (\vert B_0\vert+ 2\vert B_1\vert + \vert B_2\vert ) t \vert \gamma_1 - \gamma_2\vert \vert\gamma_2\vert^{t-1}  \bigg].
\end{align*}
Now let $\min(|\gamma_1 - |\gamma_2| |,|\gamma_1 - |\gamma_3||) = \beta\geq (\sigma\eta L)^{1/3}$. Notice that 
$$\frac{\gamma_1^t}{\vert\gamma_2\vert^t} = ( 1 + \frac{\gamma_1-\vert\gamma_2\vert}{\vert\gamma_2\vert})^t > t  \frac{\gamma_1-\vert\gamma_2\vert}{\vert\gamma_2\vert } \geq \frac{\beta}{\vert\gamma_2\vert} t .$$
Therefore, $t \vert\gamma_2\vert^{t-1} \leq \frac{1}{\beta} \gamma_1^t$. Hence, 
\begin{align*} 
&|\alpha_2\gamma_2^t + \alpha_3\gamma_3^t| \\
&\leq \frac{\vert \gamma_2-\gamma_3\vert }{\vert\Delta\vert } (3\vert B_0 \vert+ 2\vert B_1\vert + \vert B_2\vert ) \vert\gamma_2\vert^t + \frac{\vert \gamma_2-\gamma_3\vert }{\vert\Delta\vert }  (\vert B_0\vert+ 2\vert B_1\vert + \vert B_2\vert ) t \vert \gamma_1 - \gamma_2\vert \vert\gamma_2\vert^{t-1}  \\
&\leq \frac{ (3\vert B_0 \vert+ 2\vert B_1\vert + \vert B_2\vert ) \vert\gamma_2\vert^t}{\vert \gamma_1-\gamma_2\vert \vert \gamma_1-\gamma_3\vert } + \frac{ (\vert B_0\vert+ 2\vert B_1\vert + \vert B_2\vert ) t   \vert\gamma_2\vert^{t-1} }{  \vert \gamma_1-\gamma_3\vert }\\
& \leq \frac{ (3\vert B_0 \vert+ 2\vert B_1\vert + \vert B_2\vert )  \gamma_1^t}{\beta^2 } + \frac{ (\vert B_0\vert+ 2\vert B_1\vert + \vert B_2\vert )    \gamma_1^{t} }{  \beta^2  }\\
&= \frac{ (4\vert B_0 \vert+ 4\vert B_1\vert +2  \vert B_2\vert )  \gamma_1^t}{\beta^2 }. 
\end{align*}
At last, it is easy to check
$$ \vert \alpha_1 \gamma_1^t\vert \leq \frac{(\vert B_0\vert + 2 \vert B_1\vert + \vert B_2\vert) \gamma_1^t }{\beta^2}.$$
Therefore,
$$\vert B_t\vert \leq \frac{ (5\vert B_0 \vert+ 6\vert B_1\vert +3  \vert B_2\vert )  \gamma_1^t}{\beta^2 }\leq \frac{ (5\vert B_0 \vert+ 6\vert B_1\vert +3  \vert B_2\vert )  \gamma_1^t}{(\sigma\eta L)^{2/3} }.  $$\qedd

\end{document}